\documentclass[11pt]{article}
\usepackage{amsfonts,amssymb,amsmath}
\hsize \textwidth
\advance \hsize by -\marginparwidth
\evensidemargin \oddsidemargin

\usepackage{amsfonts,amssymb,amsmath,pdfsync,numbysec,bbm,xcolor}
\usepackage{amsthm}
\usepackage[shortlabels]{enumitem}
\usepackage[normalem]{ulem}
\usepackage{cancel}
\usepackage{scalerel}

\usepackage{mathrsfs}

\advance\hoffset by 5mm
\newcommand{\pdr}[2]{\frac{\partial{#1}}{\partial{#2}}}

\newcommand{\Rm}{{\mathbb R}}
\newcommand{\R}{\mathbb R}
\newcommand{\eps}{\varepsilon}
\newcommand{\commentout}[1]{}
\newcommand{\E}{{\mathbb E}}

\newcommand{\Pm}{{\mathbb P}}

\newcommand{\cE}{\mathcal{E}}
\newcommand{\cL}{\mathcal{L}}

\newcommand{\farc}{\frac}
\newcommand{\be}{\begin{equation}}
\newcommand{\ee}{\end{equation}}
\newcommand{\vphi}{\varphi}
\newcommand{\bal}{\begin{aligned}}
\newcommand{\enbal}{\end{aligned}}

\newcommand{\one}{{\mathbbm{1}}}
\newcommand{\1}{\one}
\newcommand{\ifnty}{\infty}

\newcommand{\breta}{\zeta}
\newcommand{\cnk}[2]{\begin{pmatrix}{#1}\cr {#2}\cr\end{pmatrix}}

\newcommand{\vot}{\hbox{Vote}}
\newcommand{\cT}{{\mathcal T}}

\newcommand{\urd}{u_{\rm rd}}
\newcommand{\urcl}{u_{\rm rcl}}
\newcommand{\wrcl}{w_{\rm rcl}}

\newcommand{\bhat}[1]{\expandafter\hat#1} 
\newcommand{\tinyfrac}[2]{\text{\tiny{$\frac{#1}{#2}$}}}

\newcommand{\smuG}{\mu_{\scaleto{G}{2.8pt}}}
\newcommand{\tmuG}{\mu_{\scaleto{G}{4pt}}}

\newcommand{\doverline}[1]{\overline{\overline{#1}}}

\newcommand{\dy}{\, dy}
\newcommand{\dx}{\, dx}
\newcommand{\dz}{\, dz}
\newcommand{\ds}{\, ds}

\newcommand{\tgammacoefficient}{\omega}

\numberwithin{equation}{section}

\usepackage[colorlinks=true,pdfpagemode=UseNone,urlcolor=blue,linkcolor=blue,citecolor=blue,breaklinks=true]{hyperref}
\usepackage{cleveref}
\usepackage{autonum}

\newtheorem{thm}{Theorem}[section]
\newtheorem{lem}[thm]{Lemma}
\newtheorem{cor}[thm]{Corollary}
\newtheorem{prop}[thm]{Proposition}

\newtheorem{defn}[thm]{Definition}

\renewcommand{\Cref}{\cref}
\crefname{thm}{Theorem}{Theorems}
\crefname{lem}{Lemma}{Lemmas}
\crefname{prop}{Proposition}{Propositions}
\crefname{cor}{Corollary}{Corollaries}
\crefname{appendix}{Appendix}{Appendices}
\crefname{defn}{Definition}{Definitions}
\crefname{section}{Section}{Sections}

\newcommand{\btilde}[1]{\expandafter\tilde#1}

\usepackage{tikz}
\newcommand{\BMDrawer}[7]{
\draw[#4] (#6,#7)
\foreach \x in {1,...,#1}
{   -- ++(#2,rand*#3)
}
node[right] {#5};
}

\newcommand{\BBMTree}{\begin{tikzpicture}
\draw[thick,->] (0,0) -- (7.5,0);
\draw[thick,<->] (0,-3) -- (0,3);
\pgfmathsetseed{1}
\BMDrawer{375}{0.02}{0.2}{black}{}{0}{0}
\pgfmathsetseed{6}
\BMDrawer{290}{0.02}{0.2}{blue}{}{1.7}{-1.93}
\pgfmathsetseed{2}
\BMDrawer{236}{0.02}{0.2}{red}{}{2.78}{-.4}
\pgfmathsetseed{4}
\BMDrawer{163}{0.02}{0.2}{brown}{}{4.25}{1.82}
\pgfmathsetseed{3}
\BMDrawer{105}{0.02}{0.2}{violet}{}{5.4}{-1.2}
\pgfmathsetseed{8}
\BMDrawer{46}{0.02}{0.2}{gray}{}{6.578}{.5}
\draw[thick,orange] (1.7,-.1) -- (1.7, .1) node[anchor=south] {$t_1$};
\draw[thick,orange] (2.78,-.1) -- (2.78, .1) node[anchor=south] {$t_2$};
\draw[thick,orange] (4.25,-.1) -- (4.25, .1) node[anchor=south] {$t_3$};
\draw[thick,orange] (5.4,-.1) -- (5.4, .1) node[anchor=south] {$t_4$};
\draw[thick,orange] (6.578,-.1) -- (6.578, .1) node[anchor=north] {$t_5$};

\draw[thick,->] (9,-3) -- (15,-3);
\draw[thick,black] (9,0) -- (15,2.5);
\draw[thick, blue] (9,0) -- (15,-2.5);
\draw[thick,red] (10.08,-.45) -- (15,0);
\draw[thick,brown] (11.55, -1.0625) -- (15,-.95);
\draw[thick,violet] (12.7,1.5416666) -- (15, 1);
\draw[thick,gray] (13.878, -.98) -- (15, -1.75);
\draw[thick,orange] (9,-3.1) -- (9, -2.9) node[anchor=south] {$t_1$};
\draw[thick,orange] (10.08,-3.1) -- (10.08, -2.9) node[anchor=south] {$t_2$};
\draw[thick,orange] (11.55,-3.1) -- (11.55, -2.9) node[anchor=south] {$t_3$};
\draw[thick,orange] (12.7,-3.1) -- (12.7, -2.9) node[anchor=south] {$t_4$};
\draw[thick,orange] (13.878,-3.1) -- (13.878, -2.9) node[anchor=south] {$t_5$};
\end{tikzpicture}}

\begin{document}

\author{Jing An\footnote{Department of Mathematics, Duke University, Durham, NC 27708, USA;
jing.an@duke.edu}  
\and Christopher Henderson\footnote{Department of Mathematics, University of Arizona, 
Tucson, AZ 85721, USA;
ckhenderson@math.arizona.edu} \and 
Lenya Ryzhik\footnote{Department of Mathematics, Stanford University, Stanford, CA 94305, USA;
ryzhik@stanford.edu}}

\title{Quantitative steepness, semi-FKPP reactions, and pushmi-pullyu fronts}

\maketitle

\begin{abstract}
We uncover a seemingly previously unnoticed algebraic structure of a large class of reaction-diffusion equations  and use it, in particular, to study the long time behavior of the solutions and their convergence to traveling waves in the pulled and pushed regimes, as well as at the pushmi-pullyu boundary. One such new object introduced in this paper is the shape defect function, which, indirectly, measures the difference between the profiles of the solution and the traveling wave. While one can recast the classical notion of `steepness' in terms of the positivity of the shape defect function, its positivity  can, surprisingly,  be used in numerous quantitative ways. In particular, the positivity is used in a new weighted Hopf-Cole transform and in a relative entropy approach  that play a key role in the stability arguments. The shape defect function also gives a new connection between reaction-diffusion equations and reaction conservation laws at the pulled-pushed transition. Other simple but seemingly new algebraic constructions in the present paper   supply various unexpected inequalities sprinkled throughout the paper. Of note is a new variational formulation that applies equally to pulled and pushed fronts, opening the door to an as-yet-elusive variational analysis in the pulled case.
\end{abstract}


\section{Introduction}

%
%

We consider the long time behavior of the solutions to reaction-diffusion equations of the form
\be\label{e.rde}
u_t=u_{xx}+f(u).
\ee
Throughout the paper, we assume that the nonlinearity $f(u)$ is non-negative on $[0,1]$ and satisfies
\be\label{dec2804bis}
f(0)=f(1)=0,~~f'(0)>0,~~f(u)>0\hbox{ for all $u\in(0,1)$.}
\ee
For simplicity, we also assume that $f(u)$ is as smooth as needed. 

Reaction-diffusion equations of this type have been extensively studied, dating back nearly a century to the pioneering works of Fisher~\cite{Fisher}
and Kolmogorov, Petrovskii and Piskunov~\cite{kpp}. 
The main goal of the present paper is to uncover novel algebraic structure and 
properties of these equations that we found to be quite surprising. Indeed, they are somewhat `concealed' and seem to be of an independent interest.  As an application of this theory, we analyze the long time asymptotics of the solutions to~\eqref{e.rde}; that is, we show that the convergence of $u(t,\cdot)$ to the minimal speed traveling wave and identify the precise moving frame in which this convergence occurs. 
To make this more precise, we first recall some well-established notions.

\subsubsection*{Traveling waves and their speed}
 
Under the positivity assumption (\ref{dec2804bis}), 
there exists a minimal speed~$c_*>0$ so that for all $c\ge c_*$ the reaction-diffusion equation (\ref{e.rde}) admits traveling wave solutions  
of the form~$u(t,x)=U_c(x-ct)$.
The profiles $U_c(x)$ satisfy
\be\label{dec2806bis}
	-cU_c'=U_c''+f(U_c),
		\quad
	U_c(-\infty)=1,
		\quad\text{ and }\quad
	U_c(+\infty)=0.
\ee
The minimal speed $c_*$ is characterized by the variational formula
of~\cite{HadelerRothe}:
\be\label{dec2808bis}
c_*[f]=\inf_{p\in{\mathcal K}}\sup_{u\in[0,1]}\left(p'(u)+\farc{f(u)}{p(u)}\right).
\ee
Here, $\mathcal K$ is the class of continuously differentiable functions $p(u)$ such that 
\be\label{dec2810bis}
	p(0)=0,
	\quad p'(0)>0,
	\quad\text{ and }
	\quad p(u)>0 ~\text{ for } u\in(0,1).
\ee
Existence of traveling waves was first established  in the original papers~\cite{Fisher,kpp} 
for the Fisher-KPP type nonlinearities. Recall that $f(u)$ is of the Fisher-KPP type if it satisfies, in addition to~\eqref{dec2804bis},
\be\label{mar408bis}
f(u)\leq f'(0)u,~~\hbox{ for all $0\le u\le 1$.} 
\ee
For the Fisher-KPP nonlinearities, the minimal speed given by (\ref{dec2808bis}) is explicit:
\be\label{jul1402}
c_*=2\sqrt{f'(0)}.
\ee
Such nonlinearities have attracted enormous attention throughout the past few decades.  This is perhaps due to the fact that, under this assumption, 
the linearized model (around $u = 0$) acts as a reasonably faithful approximation of~\eqref{e.rde} under the Fisher-KPP type assumption and, hence, 
many computations can be done semi-explicitly.  This allows one to study~\eqref{e.rde} in impressive details with highly precise results. We refer to a very
recent paper~\cite{Zlatos-lin} for an elegant formulation of this linearizability property of Fisher-KPP reactions. 
Moreover, there is a large class of nonlinearities, the `McKean nonlinearities', all satisfying the Fisher-KPP condition, 
that connect solutions to (\ref{e.rde})  to branching Brownian motion, as discovered by McKean in~\cite{McK}.  
We explain in Section~\ref{sec:voting} below how this connection can be
extended to a much larger class of nonlinearities $f(u)$, not only of the Fisher-KPP type.

\subsubsection*{Convergence in shape} 
Another seminal result of the original 
KPP paper~\cite{kpp} is that the  solution $u(t,x)$   to (\ref{e.rde}) with the step-function initial condition
$u(0,x)=\one(x\le 0)$ `converges in shape' to a minimal speed traveling wave. That is, there exists a reference frame $m(t)$ such that
\be\label{mar808bis}
|u(t,x+m(t))-U_*(x)|\to 0~~~~\hbox{ as $t\to+\infty$, uniformly in $x\in\Rm$.}
\ee
Here, $U_*(x)$ is the traveling wave solution to (\ref{dec2806bis}) with speed $c=c_*$, normalized, for example, so that~$U_*(0)=1/2$.
While (\ref{mar808bis}) was stated in~\cite{kpp} only for the Fisher-KPP type nonlinearities, 
with an extra assumption $f'(u)\le f'(0)$ for $u\in[0,1]$, 
the original proof can be easily adapted to a large class of nonlinearities $f(u)$. The KPP paper
also showed that 
\be\label{jul1402bis}
m(t)=c_*t+o(t),~~\hbox{as $t\to+\infty$}.
\ee
Fisher formally argued in~\cite{Fisher} that 
\be\label{jul1404}
m(t)=c_*t+O(\log t),~~\hbox{as $t\to+\infty$},
\ee
although his prediction for the coefficient in front of the logarithmic correction
turned out to be incorrect.

%
%

\subsubsection*{Pushed and pulled fronts in reaction-diffusion equations} 

As noted above, convergence in shape to a traveling wave in (\ref{mar808bis})  
holds for a much larger class of nonlinearities
$f(u)$ than, say, the Fisher-KPP nonlinearities. On the other hand, the rate of 
convergence and the asymptotics of the front location $m(t)$ depend strongly on the nature of the spreading. 
One needs
to distinguish between the so-called `pulled' regime in which propagation is dominated by the behavior far ahead of the front 
where $u\approx 0$ and the `pushed' regime in which spreading is governed by the behavior of the solution near the front where~$u $ is close
neither to~$0$ nor~$1$.  In a sense, this is a question of whether the propagation is linearly determined (pulled fronts) or nonlinearly determined (pushed fronts).  

In the pushed case, convergence in  (\ref{mar808bis})  is exponential in time.  This is, roughly, due to the fact that the important behavior occurs 
in the {\em compact} region around the front.   In addition, in the pushed case
the front location $m(t)$ has the asymptotics 
\be\label{jul1306}
m(t)=c_*t+x_0+o(1)~~\hbox{ as $t\to+\infty$.}
\ee
These results go back to the classical paper~\cite{fife1977approach} by Fife and McLeod.

In the pulled case, the convergence rate in  (\ref{mar808bis}) is algebraic in time, roughly due to the fact that the important behavior occurs on the {\em non-compact} half-line to the right of the front.  Moreover, when the Fisher-KPP condition (\ref{mar408bis}), which guarantees the pulled nature of the front, is  
satisfied,~$m(t)$ has the asymptotics 
\be\label{jul1308}
m(t)=c_*t-\farc{3}{2\lambda_*}\log t+x_0+o(1)~~\hbox{ as $t\to+\infty$.}
\ee
It is notable for the unbounded `delay' between $m(t)$ and the moving frame  $m_{TW}(t)=c_*t+x_0$  of the traveling wave.  
The formula~\eqref{jul1308} was first established by Bramson in~\cite{Bramson1,Bramson2} with probabilisitic
tools for the aforementioned `McKean sub-class' of the Fisher-KPP  
nonlinearities $f(u)$ for which the reaction-diffusion equation (\ref{e.rde}) is 
directly connected to branching Brownian motion.  
More recent analytical proofs allowed for general Fisher-KPP type nonlinearities~\cite{Graham,NRR1,NRR2}.
The algebraic rates of convergence for pulled fronts have been investigated 
in~\cite{AHS,AS1,berestycki2017exact,berestycki2018new,Ebert-vanSaarlos,Graham,NRR2,vanSaarlos}.  
%

Another distinction between the pushed and pulled cases is in the shape of the traveling wave.  Pushed traveling waves have purely exponential asymptotics:
\be\label{jul1310}
U_*(x)\sim B e^{-\lambda _*x},~~\hbox{as $x\to+\infty$},
\ee
while pulled traveling waves have an extra linear factor:
\be\label{jul1312}
U_*(x)\sim (Dx+B)e^{-\lambda _*x},~~\hbox{as $x\to+\infty$}.
\ee
Note that the exponential decay rate $\lambda_*$  appears both in (\ref{jul1312}) and in the pulled front location asymptoptics~(\ref{jul1308}). 

It is well known that the fronts associated to Fisher-KPP type nonlinearities, that is, those satisfying~(\ref{mar408bis}),
are pulled. However, the Fisher-KPP criterion does not describe 
the boundary of the pushed-pulled transition. 
While the transition from the pulled to pushed behavior has been extensively studied in the applied literature, see~\cite{Ebert-vanSaarlos,Leach-Needham,vanSaarlos} and references therein, 
including the recent study of the stochastic effects~\cite{BHK1,BHK2},  
a true mathematical understanding is still lacking. 
We mention~\cite{garnier2012inside} for a general criterion for pulled fronts, \cite{Giletti} for a preliminary investigation into the asymptotics of $m$, the
very recent papers~\cite{AHS,AS1,AS2} for  
a spectral approach to this question, and~\cite{AHR} for the study of the pushed-pulled 
transition in the context of the Burgers-FKPP equation.  We also mention~\cite{Crooks} for an investigation into some explicitly solvable cases similar to our setting.  Shedding new light on this old problem is one motivation in this work.

\subsection{New objects and main results}


%
%
%
%
%
%

\subsubsection*{The shape defect function and an energy functional} 

The proof of (\ref{mar808bis}) in~\cite{kpp} and many of the later references relies 
on a `steepness comparison' between the solution $u(t,x)$ to (\ref{e.rde}) with the initial condition $u(0,x)=\one(x\le 0)$ and the minimal speed 
traveling wave $U_*(x)$. Roughly, if $u$ starts `steeper' than $U_*$, then it remains so.  We refer to~\cite{GM} for a recent beautiful exposition of these ideas
and to~\cite{Bachmann} for its adaptation to a time-discrete setting. We simply note that $u(t,\cdot)$ being `steeper' than $U_*$ means that any shift of $U_*$ intersects the profile of $u(t,\cdot)$ precisely once.

We would like to introduce a seemingly
new way to quantify these elegant arguments. 
It is 
well known that, for all $c\ge c_*$, the traveling wave profiles $U_c(x)$  
are strictly  decreasing in~$x$. This can be seen, for instance,  
by the sliding method~\cite{BerNir}.   Thus, there exists
a `traveling wave profile function'~$\eta_c(u)$, so that~$U_c(x)$ satisfies a {\em first order} 
ordinary differential equation
\be\label{dec2812bis}
	-U_c'=\eta_c(U_c),~~~
	U_c(-\infty)=1,
		~~~
	U_c(+\infty)=0,
\ee
in addition to (\ref{dec2806bis}).  
Thus, to quantify the difference between the shape of the solution $u(t,x)$ to~\eqref{e.rde}, we propose the following object:
\begin{defn}\label{d.sdf}
	The  shape defect function is
	\be\label{dec2121bis}
		w(t,x)=-u_x(t,x)-\eta_c(u(t,x)).
	\ee
\end{defn}
To the best of our knowledge,
this notion is new. 
Additionally, introduce the energy associated to the shape defect function:
\be\label{jul1302}
	\cE_c(t)=\int e^{cx}|w(t,x+ct)|^2dx.
\ee
We should stress that the traveling wave profile function is rarely explicit,
except in the pushmi-pullyu and some pushed cases pointed out below, but its explicit form is not needed.  A generalization to higher dimensions is discussed in \Cref{s.sdf} below.

We mention two related concepts that have been considered in the past.  Recently, Matano and Pol\'a\v{c}ik    leveraged 
the trajectories $\tau_t = \{(u(t,x),u_x(t,x)): x \in \R\} \subset \R^2$ to access phase plane methods and deduce strong results in the context of propagating terraces~\cite{MatanoPolacik,Polacik}.  Even earlier, Fife and McLeod~\cite{FifeMcLeod_PP} studied the convergence to pushed traveling waves by deriving a degenerate partial differential equation for $p(t,z) = u_x(t, u(t)^{-1}(z))$.  In our notation, these approaches have, as goals, to prove the convergence of $\tau_t \to \{(z, -\eta_*(z))\}$ and $p(t,z) \to -\eta_*(z)$, respectively, as $t\to\infty$.

The early method of the original KPP paper~\cite{kpp} and  some subsequent work used the intersection number
to compare the steepness of $u(t,\cdot)$ and $U_*$ and relied on the fact that it is decreasing~\cite{Angenent,Matano-78}.  The shape defect function  
provides a new way to capture and {\em quantify} this insightful notion.  Indeed, a simple calculation yields:
\begin{prop}
For any time $t>0$, the profile $u(t,\cdot)$ is as steep as $U_*$ if and only if the shape defect function satisfies:
\be
	w(t,x) = -u_x(t,x) - \eta_*(u(t,x)) \geq 0
		\qquad\text{ for all } x \in \R.
\ee
\end{prop}
\noindent Further, the larger $w$ is, roughly, the more severe the difference in steepness.

The key property noticed in the original KPP paper~\cite{kpp} is that steepness ordering is preserved; that is, if $u_0$ is steeper than $U_*$, then so is $u(t,\cdot)$.  Again, their argument is based on the intersection number.  
By recasting this in terms of the shape defect function, this property is:
\be\label{aug202}
	\text{if } ~~
	w(0,\cdot) \geq 0, ~~
	\text{ then } ~~
	w(t,\cdot) \geq 0 ~~
	\text{ for all }
	t\geq 0.
\ee
This can be established, for instance, using the parabolic equation \eqref{dec2120} below
satisfied by $w$ and applying the maximum principle (see~\eqref{e.c92201}).  
Surprisingly, the positivity of the shape defect has numerous other, 
much more quantitative consequences than (\ref{aug202}), and this is one of the main
points of this paper.  In particular, the equation~\eqref{dec2120} for $w$ 
can be analyzed directly, for example, in order to obtain bounds on $w$, as in
\Cref{l.w_bounds}.  In addition, the positivity of $w$
leads to, for example, the new weighted Hopf-Cole transform and entropy inequalities
discussed below.
Such bounds and inequalities  
play a key role in the arguments throughout the paper,
and while they use crucially the positivity of $w$, they 
are of a different nature than the intersection number arguments.

Importantly, the shape defect function also measures, in a sense, the distance between $u$ and the `nearest' traveling wave.  Indeed, if $u(t,x)=U_c(x-ct-x_0)$ with some shift $x_0\in\Rm$, then~$w(t,x)\equiv 0$ because of (\ref{dec2812bis}). 
An elementary computation in Section~\ref{s.sdf} shows that for any $c\ge c_*$ 
the reaction-diffusion equation~(\ref{e.rde}) is a gradient flow for $\cE_c$ and
\be\label{jul1304}
\farc{d\cE_c(t)}{dt}\le 0.
\ee
This is a natural and, as far as we are aware, new way to quantify the convergence in 
shape result in (\ref{mar808bis}). 
Alternative variational formulations for reaction-diffusion equations have been previously introduced,
for instance, in~\cite{GR,LMN,MN1,MN2,Risler}. However, again, to the best of our knowledge, this is the first kind of energy
that can be used for equations of the Fisher-KPP type at the minimal speed~$c=c_*$ or, more broadly, for pulled type reaction-diffusion equations.  
Previous variational formulations were restricted   to pushed type reaction-diffusion equations. 
Moreover,~\eqref{jul1304} gives
a quantitative reason behind the convergence to traveling waves in reaction-diffusion equations.


Let us stress that the definition of the shape defect function
implicitly depends on the choice of the 
speed~$c\ge c_*$ via the function $\eta_c(u)$. Unless otherwise specified,
we use it with~$c=c_*$ and also use the notation
\be\label{jul1502}
\eta_*(u)=\eta_{c_*}(u).
\ee

\subsubsection*{The front location asymptotics}

To be concrete, in this paper we consider a family of nonlinearities of the form 
\be\label{jul1314}
f(u)=f'(0)(u-A(u))(1+\chi A'(u)).
\ee
We assume that 
\be\label{jul1802}
	A(0) = 0 = A'(0),\quad
	A(1)=1,\quad
	\hbox{and}\quad
	\hbox{$\alpha(u) := \dfrac{A(u)}{u} \in C^2$ is increasing and convex.}
\ee
As a consequence of (\ref{jul1802}), we know that $A(u)$ itself is increasing
and convex. 
This class of nonlinearities appears to be fairly general: for example,
it includes the nonlinearities suggested in~\cite{Ebert-vanSaarlos,HadelerRothe,murray2007mathematical}
\be\label{e.c110601}
	f(u) = u(1-u^{n-1})(1 + n \chi u^{n-1})
		\qquad n \geq 2.
\ee
While the Fisher-KPP nonlinearity has an extremely elegant interpretation in terms of the 
position of the maximal particle of branching Brownian motion, due to McKean~\cite{McK}, this 
is not possible with other simple models like~\eqref{e.c110601}, even when $\chi=0$.  
In \Cref{sec:voting} we investigate a connection to branching Brownian motion 
``voting models,'' which applies to, e.g.,~\eqref{e.c110601} 
and which provides a motivation for the convexity assumption on $\alpha(u)$.

Let us comment on the form of $f(u)$ in~\eqref{jul1314} and the assumptions on $\alpha(u)$ in~\eqref{jul1802}.  The main advantage of~\eqref{jul1314}-\eqref{jul1802} is that it helps highlight a concealed algebraic structure underlying~\eqref{e.rde}.  As we discuss in \Cref{s.A}, for any $f(u)$ one can find $\chi$ and~$A(u)$ satisfying~\eqref{jul1314}, so the only assumption we are making is that $\alpha(u)$ is $C^2$ and convex.  At the same time, writing $f(u$) in the form of~\eqref{jul1314} helps to make certain aspects of the structure more apparent. 
For example, the form~\eqref{jul1314} allows one to see explicitly the connection between~\eqref{e.rde} and the reactive conservation law~\eqref{e.rcl}, as well as derive the equation satisfied by the weighted Hopf-Cole transform~\eqref{jul1902} of $u$ (see~\eqref{jul1904}). 
On the other hand, the assumption~\eqref{jul1802} on the convexity of $\alpha$ is clearly not always satisfied for all~$f(u)$.  It, however, has the advantage of allowing for elegant proofs that bypass extra technicalities and highlight the role played by the algebraic structure.  Indeed, the convexity of $\alpha$ allows us to discard several errors terms that have a ``good sign.'' Otherwise, we believe that these terms could be
estimated by the smallness of $w(t,x)$ at the price of loss of elegance (see below).  
Thus, the convexity assumption
on $\alpha(u)$ can be relaxed.   Nonetheless, those $f$ satisfying~\eqref{jul1314}-\eqref{jul1802}  
present a rich class of nonlinearities in which to perform our investigation.

Nonlinearities of the form (\ref{jul1314})-(\ref{jul1802}) are of the Fisher-KPP type
when $\chi=0$ and also when~$\chi$ is positive but sufficiently small.
However, they are not of Fisher-KPP type for $\chi$  close to $\chi=1$ (see~\eqref{mar804} and \Cref{l.semiFKPP_range}). 
We show in Section~\ref{sec:semi-push} that traveling 
waves for such nonlinearities have asymptotics~(\ref{jul1312}) with~$D\neq 0$  when $0\leq \chi<1$. We refer to such nonlinearities as {\bf semi-FKPP} type, and one can show that they have speed $c_* = 2\sqrt{f'(0)}$ (see \Cref{prop-dec2804}).  Thus, these remain pulled waves, despite $f$ not satisfying the Fisher-KPP condition.

When $\chi = 1$, the traveling wave decay becomes purely exponential~\eqref{jul1310}, as in the pushed case.
Nevertheless, the speed remains linearly determined $c_*  = 2\sqrt{f'(0)}$, as in the pulled case.  This is the
boundary between the pulled and pushed regimes.  As such, we call this the {\bf pushmi-pullyu} case, following~\cite{AHR}.  
Below, we discuss various remarkable algebraic properties of the solutions to (\ref{e.rde}) in the pushmi-pullyu case.
Here, we simply note that 
the traveling wave profile function~$\eta_*(u)$ is explicit when $\chi = 1$ (see \Cref{prop-dec2802}):
\be\label{jul1506}
\eta_*(u)=u-A(u). 
\ee
An immediate consequence is the purely exponential decay of the minimal speed traveling wave $U_*$ mentioned above
that follows from the regularity assumptions on $A(u)$.  The identity~\eqref{jul1506}, however, has many further consequences in the
analysis  of the pushmi-pullyu case.  In the semi-FKPP case, the traveling wave profile $\eta_*(u)$ does not have the simple form
(\ref{jul1506}) and is not explicit but satisfies the bounds
\begin{equation}\label{aug204}
\sqrt{\chi}(z-A(z))\leq \eta_*(z)\leq z-A(z).
\end{equation}
proved in Lemma~\ref{lem-jun1402}.

Let us also comment that in this paper we do not treat the pushed regime $\chi>1$ as it can be handled by the standard 
existing methods originating in~\cite{fife1977approach}.  For completeness, we provide a short proof that these fronts are pushed -- see \Cref{p.large_chi}.

We now state  the main theorem of this paper on the asymptotics of the front location for the reaction-diffusion equation (\ref{e.rde}) with nonlinearities of the form
\eqref{jul1314}-(\ref{jul1802}).
We also assume that the initial condition $u_0(x)=u(0,x)$ satisfies
\be\label{mar1614}
0\le u_0(x)\le 1,~~\hbox{ for all $x\in\Rm$},
\ee
and is compactly supported
on the right: there is $L_0\in\Rm$ so that
\be\label{mar1616}
u_0(x)= 0,~~\hbox{ for all $x\ge L_0$}.
\ee

\begin{thm}\label{t.rde-intro}
Under the above assumptions, suppose that $u(t,x)$ solves~\eqref{e.rde} and that initially the shape defect function is non-negative: $w(0,x)\ge 0$ for all $x\in\Rm$. \\
(i) If $0\le\chi<1$ so that  $f$ is a semi-FKPP type nonlinearity, the front location has the asymptotics:
\be\label{mar1620bis3}
m(t)=2t-\farc{3}{2}\log t+x_0+o(1),~~\hbox{ as $t\to+\infty$.}
\ee
(ii) If  $\chi=1$ so that $f$ is a pushmi-pullyu type nonlinearity,  
then $m(t)$ has the asymptotics
\be\label{mar1620bis2}
m(t)=2t-\farc{1}{2}\log t+x_1+o(1),~~\hbox{ as $t\to+\infty$.}
\ee
The constants $x_0$ and $x_1$ in (\ref{mar1620bis3}) and (\ref{mar1620bis2}) depend on the  initial condition $u_0$ for (\ref{e.rde}) and the
nonlinearity~$f$. 
\end{thm}

The asymptotics (\ref{mar1620bis3}) for semi-FKPP type nonlinearities is exactly the same as (\ref{jul1308}) for the
Fisher-KPP nonlinearities. However, in the pushmi-pullyu
case, the logarithmic correction in (\ref{mar1620bis2}) is different. This change has been predicted 
in~\cite{Ebert-vanSaarlos,Leach-Needham,vanSaarlos} using formal matched asymptotics for the situations when the minimal speed traveling wave has purely exponential
decay as in (\ref{jul1310}).
To the best of our knowledge, the only rigorous result in this direction is the expansion
\be\label{mar1623bis}
m(t)=2t-\farc{1}{2}\log t+o(\log t),~~\hbox{ as $t\to+\infty$,}
\ee
obtained in~\cite{Giletti} by a careful gluing of sub- and supersolutions, a very different
approach from the present paper.

Let us point out that a simple consequence of our analysis is a  way to identify 
some pulled fronts: 
waves associated to~\eqref{e.rde} are pulled if $f\leq f_\chi$ for some $f_\chi$ of the form~\eqref{jul1314}-\eqref{jul1802} with $\chi < 1$.  
This is more general than the Fisher-KPP condition and, unlike the commonly used condition $c_*=2\sqrt{f'(0)}$ for pulled fronts, 
allows one to distinguish between 
the pulled and pushmi-pullyu cases. Moreover, Corollary~\ref{p.large_chi} in 
Appendix~\ref{sec:appendix}
shows that if $f\ge f_\chi$ with $\chi>1$ then the fronts for (\ref{e.rde}) are pushed.
Actually, if $f(u)=f_\chi(u)$ and $\chi>1$, then both the minimal speed $c_*$ and the traveling wave profile function
$\eta_*(u)$ are explicit, see Proposition~\ref{p.large_chi-explicit}. \

One technical note is that \Cref{t.rde-intro}   does not explicitly
mention the convergence of $u$ to $U_*$ as this follows from previous results (see the discussion surrounding~\eqref{jul1402bis}).  
Due to the type of soft arguments used, these previous convergence results are not quantitative, and, thus, it is not possible to use them in any way in the present 
proof of \Cref{t.rde-intro}.  As a result, we provide an alternative proof of convergence as a part of establishing the precise asymptotics 
for~$m(t)$.

It seems that the assumption $w(0,\cdot) \geq 0$ can be relaxed.  Indeed, the lack of positivity of the shape defect function can be compensated by its smallness that, as we have mentioned, can be proved independently by directly analyzing its dynamics in~\eqref{dec2120} (see \Cref{l.w_bounds} and the forthcoming work~\cite{AHR3}).  As with the convexity assumption on $\alpha$, the positivity of $w$ allows for an elegant proof highlighting the subtle algebraic structure that we uncover in this work.

\subsubsection*{Connection to the reactive conservation laws and the precise front asymptotics}

A first example of the special structure in the pushmi-pullyu case 
is the connection to 
 the reactive conservation law:
\be\label{e.rcl}
	\mu_t + (A(\mu))_x = \mu_{xx} + \mu - A(\mu).
\ee
In particular, we show in Section~\ref{sec:cons-law-push} that the reaction-diffusion equation (\ref{e.rde}) with $f(u)$ given 
by 
\be\label{jul1508}
f(u)=(u-A(u))(1+A'(u))
\ee 
and the reactive conservation law (\ref{e.rcl}) have exactly the same   minimal speed traveling 
wave solutions.  Thus, the traveling wave profile function $\eta_*(u)$ is the same for the two equations,
and the shape defect function for (\ref{e.rcl}) is still defined 
as\footnote{Although with a notational change from $w$ to $\wrcl$ in order to not confuse it with the shape defect function for~\eqref{e.rde}.} in \Cref{d.sdf}:
\be\label{aug216}
	\wrcl(t,x)=-\mu_x(t,x)-\eta_*(\mu(t,x)).
\ee

The connection between the two equations (\ref{e.rde}) and (\ref{e.rcl}) 
extends further: if the shape defect function $w(t,x)\ge 0$, defined by (\ref{dec2121bis}) then the  solution to~(\ref{e.rde}) is a
subsolution to (\ref{e.rcl}):
\be\label{jul1316bis}
u_t+(A(u))_x\le u_{xx}+u-A(u).
\ee
This algebraic  miracle allows us to bound the solutions to the pushmi-pullyu reaction-diffusion equation and the conservation law (\ref{e.rcl})
in terms of each other:
\be\label{e.rde_rcl_compare}
	u \leq \mu.
\ee
Further novel algebraic properties of the semi-FKPP and pushmi-pullyu nonlinearities are discussed in Section~\ref{sec:intro-proof} below, as well as in
Sections~\ref{sec:semi-push} and~\ref{sec:alg}.

We also obtain the analogous result of \Cref{t.rde-intro}.(ii) for the reactive conservation law (\ref{e.rcl}). 
\begin{thm}\label{t.rcl-intro}
Let $\mu$ be a solution to~\eqref{e.rcl}, with an initial condition such that $\wrcl(0,x)\ge 0$  
for all~$x\in\Rm$.  
Then there is $x_2 \in \R$ so that, for all $L > 0$,
\be\label{e.c101301}
	\lim_{t\to\infty} \sup_{|x| \leq L} |\mu(t,x + m(t)) - U_*(x)|
		= 0
	\qquad\text{ with }~
	m(t) = 2t - \frac{1}{2}\log t + x_0.
\ee
\end{thm}
The positivity of the initial shape defect function does not appear to be necessary for the conclusions of Theorems~\ref{t.rde-intro} and~\ref{t.rcl-intro}
to hold. We have opted to use this extra assumption not simply because it shortens the proof  
but also because it makes many steps in the proofs elegant rather than technical and reveals
a number of algebraic properties that are not seen without this assumption.

\subsection{Key elements of the proofs of \Cref{t.rde-intro} and  \Cref{t.rcl-intro}}\label{sec:intro-proof} 

As the proofs are quite intricate, a detailed summary would be too long to contain here.  Instead we discuss the major aspects of the proof and,
in particular, the new tools that become available due to the positivity of the shape defect function.  A more detailed outline of the proofs can be found in \Cref{sec:outline}.

\subsubsection*{The weighted Hopf-Cole transform}

The crucial tool to compensate for the lack of the Fisher-KPP condition for $f$ is the weighted Hopf-Cole transform:
letting $\hat u(t,x) = u(t,x+2t)$, we define
\be\label{jul1902}
	v(t,x)=\exp\Big(x+\sqrt{\chi}\int_x^\infty \alpha(\hat u(t,y))dy\Big)\hat u(t,x).
\ee
One should stress a key  difference between the semi-FKPP ($\chi \in [0,1)$) and pushmi-pullyu ($\chi = 1$) cases: 
\be
	v(t,-\infty) = 0 \quad\text{ when } \chi < 1
	\qquad\text{ and }\qquad
	v(t,-\infty) > 0 \quad\text{ when } \chi = 1.
\ee


Due to the remarkable algebraic structure discussed in \Cref{sec:alg}, whenever $u$ solves~\eqref{e.rde} and the shape defect function is nonnegative, we obtain the differential inequality
\be\label{jul1904}
	v_t\le v_{xx}.
\ee
This is slightly easier to see in the pushmi-pullyu case ($\chi = 1$) where, after an intricate computation in the proof of \Cref{prop-mar902}, we obtain
\be
	v_t - v_{xx}
		\leq - 2 \alpha(\hat u) w
		\leq 0.
\ee
Intuitively, the importance of~\eqref{jul1904} is the following. When $\chi < 1$, it is clear that $v(t,x)$ is `small' for~$x< 0$.  We can then, roughly, think of $v$ as solving the heat equation on the half-line with {\em Dirichlet} boundary conditions, which implies that $v$ decays like $O(t^{-3/2})$.  This $3/2$ is the same as the one in~\eqref{mar1620bis3}.  Indeed,
\be\label{e.c101402}
	O(t^{-3/2})
		= v(t, m(t)- 2t)
		= e^{m(t) - 2t + O(1)} u(t,m(t))
		= O(e^{m(t) - 2t}).
\ee
When $\chi = 1$, we have, on the other hand,
\be
	v_x(t,x) = - \exp\Big(x+\int_x^\infty \alpha(\hat u(t,y))dy\Big) w(t,x),
\ee
which is `small' for $x<0$. We can then, roughly, think of $v$ as solving the heat equation on the half-line with {\em Neumann} boundary conditions, which implies that $v$ decays like $O(t^{-1/2})$.  This $1/2$ is the same as the one in~\eqref{mar1620bis2} by a similar computation to~\eqref{e.c101402}.

We note that this intuition can only be turned into a proof in the semi-FKPP case.  For technical reasons, it does not go through in 
the pushmi-pullyu case, and one must first obtain the `rough' front asymptotics 
\be\label{nov1002}
m(t)=2t - (1/2) \log t + O(1),
\ee  
without the use of the weighted Hopf-Cole transform.  The transform, however, does play a crucial role in upgrading the front asymptotics to 
the precision of~\eqref{mar1620bis2} (resp.~\eqref{e.c101301}) in the pushmi-pullyu case. 


\subsubsection*{Relative entropy and a weighted Nash inequality in the pushmi-pullyu case}

Let us now explain how the `rough' front asymptotics (\ref{nov1002}) is obtained. We focus on the conservation law (\ref{e.rcl}). 
Along the lines of~\eqref{e.c101402}, this asymptotics
can be recast as the~$L^\infty$ upper and lower estimates
\be\label{e.c101302}
	\quad\|p(t,x)\|_\infty = O(1/\sqrt t)
\ee
for the function
\be\label{nov1004}
			p(t,x) = e^x \hat \mu(t,x)),
\ee
where $\hat\mu(t,x) = \mu(t, x+2t)$.  It turns out that the upper bound in~\eqref{e.c101302} is much more difficult to establish, so we discuss that now.  

The intuitive reason behind~\eqref{e.c101302} is that $p(t,x)$ satisfies an inhomogeneous conservation law
\be\label{e.c101304}
	p_t + (\alpha(\hat \mu) p)_x = p_{xx}.
\ee
Heuristically, $\alpha(\hat \mu)$ pushes all of the mass of $p$ to the right, where~\eqref{e.c101304} is essentially the heat equation as $\alpha(\hat \mu) \approx 0$
as $x\to+\infty$.  This indicates that, after a boundary layer in time, we should have heat-equation-like decay $O(1/\sqrt t)$, as desired.

In order to create a proof out of this simple idea, we use a relative entropy approach.  Another surprising consequence of the steepness comparison 
of the solution and the traveling wave is that, if~$\wrcl(0,x)\ge 0$ for all $x\in\Rm$, then
the function
\be\label{jul1908}
	\rho(t,x) = \exp \Big( - \int_0^{\hat \mu(t,x)} \frac{A(u')}{u'(u'-A(u'))} \, du'\Big)
\ee
is a supersolution to (\ref{e.c101304})
\be\label{jul1910}
	\rho_t + (\alpha(\hat \mu) \rho)_x \geq \rho_{xx}.
\ee
This follows from a rather involved computation in the proof of Lemma~\ref{lem.rho}.
An observation coming from~\cite{AHR} is that if $\rho(t,x)$ is a supersolution and  $p(t,x)$ is a solution to (\ref{e.c101304})
then the function
\be
\varphi(t,x) =\farc{p(t,x)}{\rho(t,x)}= \farc{e^x \hat \mu(t,x)}{ \rho(t,x)},
\ee
obeys a dissipation inequality 
\be
	\frac{1}{2}\frac{d}{dt} \int  \varphi(t,x)^2 \rho(t,x) dx
		\leq - \int \varphi_x(t,x)^2 \rho(t,x) dx.
\ee
Were $\rho\equiv 1$, this would be the key differential inequality for the heat equation that, along with the Nash inequality and a self-adjointness `trick,' yields the $O(1/\sqrt t)$-decay desired in~\eqref{e.c101302}.  Unfortunately, $\rho\not\equiv 1$ (notice that $\rho(t,-\infty) = 0$ and $\rho(t,+\infty) = 1$) and our problem is not self-adjoint.  By developing a suitable Nash-type inequality with time-dependent dynamic weights 
to yield suitable weighted~$L^1 \to L^2$ estimates and an additional 
bootstrapping procedure to yield $L^2 \to L^\infty$ estimates, we obtain the analogous bound $\|\varphi \|_\infty \leq O(1/\sqrt t)$ that is the key step in obtaining~\eqref{e.c101302} and, thus, the asymptotics~\eqref{mar1620bis2} and~\eqref{e.c101301}.

\subsubsection*{Decay of $w$ and the reactive conservation law}

As we discussed above, the weighted Hopf-Cole transform is not used in the pushmi-pullyu case to prove the `rough' 
$2t - (1/2) \log t + O(1)$ front asymptotics.  It is, however, needed to upgrade this estimate to~\eqref{mar1620bis2} and~\eqref{e.c101301}; that is, that 
the $O(1)$ term is actually of the form $x_0 + o(1)$ for some $x_0 \in \R$.  In the reaction-diffusion case~\eqref{e.rde}, we use~\eqref{jul1904} and the arguments of~\cite{AHR} to obtain this improved precision.

Unfortunately, the Hopf-Cole transform does not yield~\eqref{jul1904} in the reactive conservation law case.  Instead, analogous computations yield
\be
	v_t
		- v_{xx}
			\leq e^\gamma (\hat \mu \alpha'(\hat \mu) - \alpha(\hat\mu)) \wrcl, 
\ee
where we abuse notation and use $v$ for the same quantity~\eqref{jul1902} with $\hat \mu$ replacing $\hat u$.  The right hand side is {\em positive} for convex $\alpha$  except in the special case $\alpha(u) = u$ considered in~\cite{AHR}.

As we do not have the miraculous Hopf-Cole cancelation so we require an additional 
argument.  In particular, we need smallness of $w_{\rm rcl}$.  As can be seen by a change to self-similar coordinates, it turns out that the required bounds are of the type 
$\|w_{\rm rcl}\|_\infty = O(1/t)$ and suitable integrability.  To obtain these bounds, we leverage the fact that $w_{\rm rcl}$ satisfies the equation
\be\label{e.c101401}
	\partial_t w_{\rm rcl} + A'(\hat \mu) \partial_x w_{\rm rcl}
		= \partial_{xx} w_{\rm rcl} + w_{\rm rcl}(1 - A'(\hat \mu)).
\ee
It is not immediately obvious why $w_{\rm rcl}$ should decay like $O(1/t)$ from the above.  Indeed, $w_{\rm rcl}$ has a positive growth rate for $x \geq 2t - (1/2) \log t + O(1)$.

Roughly, $w_{\rm rcl}$ decays because the region where it has a positive growth rate `moves too fast' to the right.  Indeed, the production rate $1- A'(\hat \mu)$ 
in (\ref{e.c101401}) 
is only positive when $\hat \mu \approx 1$, which, by the work outlined above, corresponds to $x \geq 2t - (1/2) \log t + O(1)$.  On the other hand, the linearization of~\eqref{e.c101401} is the same as that of Fisher-KPP, indicating that $w_{\rm rcl}$ `wants' to have a `front' at $2t - (3/2) \log t + O(1)$.  The $\log t$ discrepancy between this and where $w_{\rm rcl}$ has a positive growth rate, along with the natural $O(e^{-x})$ decay of the problem (see~\eqref{jul1310}) indicates that
\be
	w_{\rm rcl} \leq O(e^{-\log t}) = O\Big(\frac{1}{t}\Big).
\ee

Interestingly, it seems that analyzing the equations for $w$ and $w_{\rm rcl}$ yields sharp bounds for convergence rates of $u$ to the traveling wave by using the heuristic ideas outlined above.  We plan to explore this in a future work.

%
%
%
%
%
%
%
%
%
%

{\bf A note on the relationship between the present work and~\cite{AHR}.}
Let us note that some parts of the analysis here are motivated by the  considerations in~\cite{AHR} for the Burgers-FKPP equation,
which is the reactive conservation law (\ref{e.rcl}) in the special case $A(u)=u^2$.
We take care to present the work here to highlight {\em only} the novel elements of the proofs as compared to \cite{AHR}.  We point out that~\cite{AHR} examined 
only the reactive conservation law equation and, moreover, only the particular case $\alpha(u) = u$.  In this simple setting, there are additional cancellations and 
many explicit calculations are possible.  Hence, a key aspect of our work here is to explore the algebraic structure of the reactive conservation laws~\eqref{e.rcl} and the reaction diffusion equations~\eqref{e.rde} and to show that the seemingly {\em ad hoc} techniques in~\cite{AHR} are actually a part of this larger structure.  One example of this is the choice of $\rho$: in~\cite{AHR} a `lucky guess' led to the choice $\rho = 1- \hat u$, but it is now clear that this follows from~\eqref{jul1908} where, in that particular case, $A(u) = u^2$.

\subsection{Organization of the paper}

In \Cref{sec:voting}, we discuss  a connection between  branching Brownian motion and   reaction-diffusion
equations~\eqref{e.rde} using probabilistic voting models.  While some examples of Fisher-KPP type nonlinearities are McKean nonlinearities,
those that, roughly, give the statistics of the maximal particle in a branching Brownian motion, many simple choices fall outside this class; for example, $f(u) = u - u^n$ for any $n \geq 3$ is {\em not} McKean type.  We present a class of voting models that 
allows to go beyond the McKean class of nonlinearities to a much larger class of  nonlinearities that includes Fisher-KPP type nonlinearities such as $u - u^n$ and  many non-Fisher-KPP type nonlinearities.  
This provides an additional motivation for the setting of \Cref{t.rde-intro}.

We then discuss the basic properties of the semi-FKPP and pushmi-pullyu nonlinearities in Section~\ref{sec:semi-push}. 
In particular, 
we discuss
the properties of the traveling wave profile function $\eta_c(u)$ 
and explain why the minimal traveling wave speed is given by (\ref{jul1402}) not only for the Fisher-KPP but also for semi-FKPP nonlinearities.   

Section~\ref{sec:alg} discusses some of remarkable algebraic properties of reaction-diffusion equations that are elucidated by the
use of the shape defect and traveling wave profile functions.  There, we  consider 
in greater detail 
the shape defect function,
establish the variational formulation (\ref{jul1302}) 
for reaction-diffusion equations,
and explain the aforementioned natural connection between reaction-diffusion equations  
and reactive conservation laws. We also introduce the weighted Hopf-Cole transform (see~\eqref{jul1902}) and establish the key differential inequality~\eqref{jul1904}.

Section~\ref{sec:outline} gives the outlines of the proofs of Theorems~\ref{t.rde-intro} and~\ref{t.rcl-intro}.  These results are quite intricate, and so we provide a high level discussion of its major difficulties and how the elements of the proof introduced above fit together to overcome these obstructions.  The full proofs are contained in \Cref{sec:pushmi-pf} (pushmi-pullyu case) and \Cref{sec:sfkpp} (semi-FKPP case).

Finally, Appendix~\ref{sec:appendix} contains some auxiliary results used in the paper as well as a discussion of the generality of the model~\eqref{jul1314}.

\medskip

{\bf Acknowledgments.} CH was supported by NSF grants DMS-2003110 and DMS-2204615. 
LR was supported by NSF grants DMS-1910023 and DMS-2205497 and by ONR grant N00014-22-1-2174.  JA and CH acknowledge support of the Institut Henri Poincar\'e (UAR 839 CNRS-Sorbonne Université), and LabEx CARMIN (ANR-10-LABX-59-01).

\section{Probabilistic interpretations of reaction-diffusion equations} \label{sec:voting}

\subsection{The McKean nonlinearities}

A special class of the Fisher-KPP nonlinearities arises naturally in the context of the branching Brownian
motion (BBM), as originally discovered by McKean in~\cite{McK}. 
They have the form
\be\label{mar410}
f(u)=\gamma\Big(1-u-\sum_{k=2}^\infty p_k(1-u)^k\Big),
\ee
with $\gamma>0$ and $p_k\ge 0$ such that
\be\label{mar710}
\sum_{k=2}^\infty p_k=1.
\ee
Here, $p_k$ is the probability that a BBM particle branches into $k$ offspring
particles at a  branching event, and~$\gamma>0$ is the exponential rate of branching.  
Specifically, when $f(u)$ has the form (\ref{mar410}), the solution to (\ref{e.rde}) is given by
\be
\label{jul1912}
u(t,x)=1-\E_x\Big(\prod_{k=1}^{N_t}(1-u_0(X_k(t))\Big).
\ee
Here, $X_1(t),\dots,X_{N_t}(t)$ are the locations of the BBM particles at the time $t$. 
In the absence of branching, this is simply the standard interpretation of the solutions of the 
heat equation in terms of a Brownian motion. 
We refer to~\cite{JB,Bovier} for excellent introductions to BBM and the McKean connection
between BBM and reaction-diffusion equations. 
 
%
%

\subsection{Voting  nonlinearities}

However,  the nonlinearities that admit the McKean interpretation form just a sub-class of the Fisher-KPP nonlinearities. 
Maybe the most basic example of an FKPP
nonlinearity not in the McKean class is $f(u)=u-u^n$ with $n\ge 3$, as it can be 
easily   checked that it both satisfies (\ref{mar408bis}) 
and can not be written in the form (\ref{mar410}).

\begin{figure}
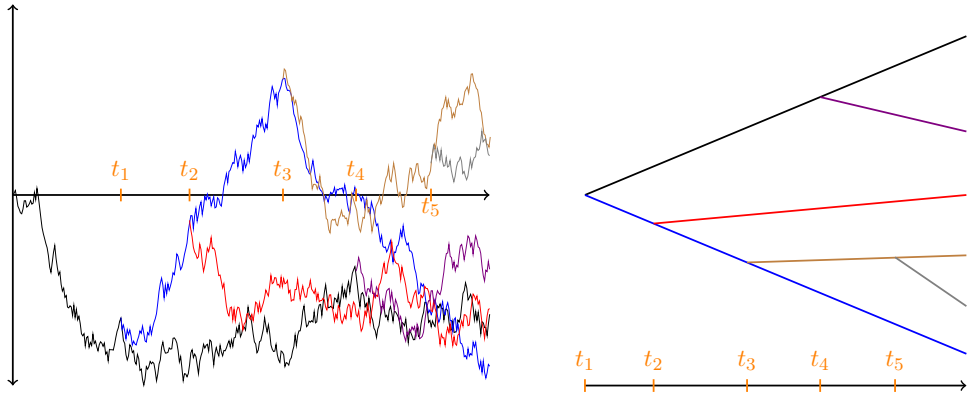

\begin{center}
\resizebox{5 in}{!}{\BBMTree}
\end{center}
\caption{On the left, a sample (binary) branching Brownian motion, which branches at times $t_i$.  On the right, the associated tree.}
\label{f.bbm}
\end{figure}

Let us now briefly describe the construction in~\cite{AHR-Vot}, originating in the  beautiful and insightful
ideas of~\cite{EFP}, that provides
a connection between a much larger class of semilinear parabolic equations and BBM than that of McKean. For the sake of concreteness
and simplicity,
we consider $f(u)=u-u^n$ and refer to~\cite{AHR-Vot} for a description of this connection in  more generality.  
Let us start with a BBM running at an exponential rate $\beta>0$: the branching times have the law
\be\label{jun706}
\Pm(\tau>t)=e^{-\beta t}.
\ee
We assume that at each branching event, the parent particle produces exactly $n$ offspring.  
 There is a natural way to associate a random genealogical tree $\cT$ to each realization of
the BBM, with each vertex of the tree corresponding to a branching event.  This is illustrated in Figure~\ref{f.bbm}. 
Each of the edges
coming out of a vertex represents an offspring particle born at that branching event. The root of the tree $\cT$ is the original particle
that started at the time~$t=0$ at the position $x$.

We now describe the voting procedure. Let us run the above BBM
until a final time $t>0$. At that time, each of the particles $X_1(t),\dots,X_{N_t}(t)$ 
that are present at the time $t$ 
votes~$1$ or $0$, with the probability
\be\label{jun904}
 \Pm(\vot(X_k(t))=1)=	1 - \Pm(\vot(X_k(t))=0) =g(X_k(t)).
\ee
Here, the function $g(x)$ is fixed and takes values in $[0,1]$.  

Given the votes of the last generation of particles that are present at the time $t$, we propagate the vote up the genealogical tree $\cT$ as follows. Let us fix $\gamma>0$ sufficiently small and define the probabilities~$\mu_{kn}$,~$k=0,\dots,n$, as 
\be\label{jun902}
	\mu_{0n}=0,
		\qquad \mu_{nn}=1,
		\qquad\text{and}\qquad
		\mu_{kn}=\farc{(1+\gamma)k}{n}~
		\text{ when } 0 < k < n.
\ee
As we need to have $0\le\mu_{kn}\le 1$ for all $k$, the parameter $\gamma>0$ must satisfy
\be\label{jun704}
	0 < \gamma \le\frac{1}{n-1}.
\ee
With the probabilities $\mu_{kn}$ in hand, given a parent  
particle on the genealogical tree $\cT$, if $k$ out of its~$n$ children voted~$1$, then the parent particle votes $1$ with the 
probability $\mu_{k}$ given by (\ref{jun902}).
Using this rule iteratively to go up the tree all the way to the root  produces the random vote~$\vot_{\rm orig}$ of the original ancestor particle, and we can define
\be\label{mar1704}
u(t,x)=\Pm_x(\vot_{\rm orig}=1).
\ee
Here, the probability is taken both with respect to the randomness in the original 
voting in (\ref{jun904}), and with respect to the randomness in the vote of each parent. 
If there was no branching event until the time $t$, so that $N_t=1$, then the vote of the original particle is $1$ with the probability $g(X_1(t))$. 

An elementary computation in~\cite{AHR-Vot} shows that the~function $u(t,x)$ defined in (\ref{mar1704}) satisfies the initial value problem
\be\label{mar1708}
\bal
&u_t=\Delta u+f(u),\\
&u(0,x)=g(x),
\enbal
\ee
with the nonlinearity
\be\label{mar1710}
\bal
f(u)&=\beta\gamma (u-u^n).
\enbal
\ee
Note that  
the range of $\gamma$ is restricted by
(\ref{jun704}) but $\beta>0$ can be arbitrary. 

We may consider the above voting model with
a  more general branching Brownian motion, with the probability $\alpha_n$ to branch into $n$ children at each branching event. If we keep the probabilities~$\mu_{kn}$
for the parent with $n$ total children to vote $1$ if $k$ of its $n$ children voted $1$,
we would obtain a convex combinations of the nonlinearities in~(\ref{mar1710}):
\be\label{jun710}
f(u)=\beta\gamma\sum_{k=1}^N\alpha_k(u-u^k),~~\sum_{k=2}^N\alpha_k=1.
\ee
They are still within the Fisher-KPP class. 
We refer to nonlinearities of the form (\ref{jun710}) as voting-FKPP nonlinearities. 
Such nonlinearities have the form 
\be\label{mar708}
f(u)=\lambda (u-A(u)),
	\quad\text{ with } \lambda = f'(0)>0.
\ee
The functions  $A(u)$ are non-negative, convex on $[0,1]$, and satisfy (\ref{jul1802}):
\be\label{mar718}
	A(0)=0,
	\quad A(1)=1,
	\quad\text{ and }\quad
	 A'(0)=0.
\ee 
It follows that the function $A(u)$ is increasing since $A'(0)=0$ and $A(u)$ is convex on $[0,1]$. 

It is sometimes convenient for us to write $f(u)$ in the form 
\be\label{mar702}
f(u)=\lambda u(1-\alpha(u)),
\ee
with 
\be\label{jun916}
\alpha(u)=\farc{A(u)}{u}.
\ee
For nonlinearities of the form (\ref{jun710}), the functions $A(u)$ and $\alpha(u)$ have the forms
 \be\label{mar712}
	 A(u)=\sum_{k=2}^N\alpha_ku^k
 		\quad\text{ and }\quad
	\alpha(u)=\sum_{k=2}^\infty \alpha_ku^{k-1}.
\ee
Therefore, in these examples both $A(u)$ and $\alpha(u)$ are increasing and convex.
We should mention that for the McKean nonlinearities (\ref{mar410}) if we were to write them in the form (\ref{mar708}), 
the function $A(u)$ is increasing and convex but $\alpha(u)$ is necessarily concave. That is, the voting Fisher-KPP nonlinearities
represent a complementary  class to the McKean type. 

Let us also briefly comment that the simple voting procedure described above can be generalized in many ways, and, 
unlike the McKean interpretation, voting models can lead 
to reaction-diffusion equations (\ref{e.rde}) with nonlinearities $f(u)$
that need not  be of the Fisher-KPP 
type. This is done simply by changing the voting rules, without changing the underlying branching Brownian motion. 
That is, one considers the same BBM, with exactly the same genealogical tree as above.
However, we  modify the probabilities $\mu_{kn}$ in~(\ref{jun902}) 
for the parent to vote~$1$ if $k$ out of its $n$ children voted~$1$.   
Let us assume for simplicity that the BBM has a fixed number $n$ of children at each 
branching event.
Then, the function
$u(t,x)$ defined by (\ref{mar1704}), which is 
the probability for the original ancestor particle 
to vote~$1$, satisfies
the reaction-diffusion equation (\ref{mar1708}), although with the nonlinearity
\be\label{aug302}
f(u)=\beta\Big(\sum_{k=0}^{2n-1}\cnk{2n-1}{k}\mu_{kn} u^k(1-u)^{2n-1-k}-u\Big).
\ee
As discussed in~\cite{AHR-Vot}, with a suitable choice of $\mu_{kn}$ one can obtain nonlinearities not of the Fisher-KPP type.
Indeed, the original example in~\cite{EFP} is the Allen-Cahn equation 
\be\label{aug304}
u_t=u_{xx}+f(u), 
\ee
with the nonlinearity 
\be\label{aug306}
f(u)=u(1-u)(2u-1),
\ee
that is not of the Fisher-KPP type. To obtain (\ref{aug304})-(\ref{aug306}), 
one considers ternary BBM,~$n=3$, and the voting probabilities $\mu_{03}=\mu_{13}=0$ 
and $\mu_{23}=\mu_{33}=1$ that come from the simple majority voting
rule. More examples can be found in~\cite{AHR-Vot}.

 \section{The algebra of semi-FKPP and pushmi-pullyu traveling waves}\label{sec:semi-push}

In this section, we explore the identities and properties of the traveling wave profile function $\eta_c(u)$ (see Propositions~\ref{prop-dec2802} and~\ref{dec28-prop06}) as well as the relationship between its regularity and the decay of the traveling wave $U_*(x)$ as $x\to+\infty$ (see \Cref{dec28-prop06}).  
We also show in Proposition~\ref{prop-dec2804} that the minimal speed for traveling waves is still given by the FKPP formula $c_*=2\sqrt{f'(0)}$ for the
nonlinearities of the semi-FKPP type.

%
%
%
%
%
%

\subsection{Semi-FKPP and pushmi-pullyu nonlinearities}  
 
We first introduce the (seemingly) new types of nonlinearities: semi-FKPP and pushmi-pullyu.
Let us start with a Fisher-KPP nonlinearity of the form 
\be\label{mar720}
\breta(u)=\lambda(u-A(u)),
\ee
with some $\lambda>0$ and an increasing convex function $A(u)$ that satisfies (\ref{mar718}). 
Note that 
\be\label{nar802}
\lambda=\breta'(0),
\ee
since $A'(0)=0$ by (\ref{mar718}). 
This class includes both the McKean and voting-FKPP nonlinearities. 

We say that a function $f(u)$ is a pushmi-pullyu nonlinearity if it has the form 
\be\label{mar721}
f(u)=\breta(u)(2\lambda-\breta'(u))=\lambda^2(u-A(u))(1+A'(u)),
\ee
with $\breta(u)$ as above. As in (\ref{nar802}), since $A'(0)=0$, it follows that $f'(0)=\lambda^2$. Thus, we
may represent a pushmi-pullyu nonlinearity in the form 
\be\label{jun906}
f(u)=f'(0)(u-A(u))(1+A'(u)),
\ee
As we shall see, this class of nonlinearities represents the boundary between those of pushed and pulled type.  Any nonlinearities smaller than $f$ should be pulled; however, they may not be Fisher-KPP type.  It is, thus, natural to say that to say that a function $f(u)$ is of the semi-FKPP type if, in contrast to~\eqref{jun906}, it satisfies 
\be\label{mar722}
f(u)\le f'(0)(u-A(u))(1+\chi A'(u)),
\ee
with $0\le\chi<1$ and an increasing convex function $A(u)$ that satisfies  (\ref{mar718}).

Let us note that a semi-FKPP nonlinearity satisfies the Fisher-KPP condition
(\ref{mar408bis}) if
\be
(u-A(u))(1+\chi A'(u))\le u,~~\hbox{ for all $0\le u\le 1$,}
\ee
or, equivalently,
\be\label{mar804}
0\le \chi \le\chi_{FKPP}=\min_{u\in[0,1]}\farc{A(u)}{A'(u)[u-A(u)]}. 
\ee
We are mostly interested in the range $\chi_{FKPP}\le\chi\le 1$, where
$f(u)$ is not of the Fisher-KPP type.  As Theorem~\ref{t.rde-intro} shows, 
the solutions to~(\ref{e.rde}) still
exhibit the Fisher-KPP type behavior. Let us comment that \Cref{l.semiFKPP_range}, below, shows that,
as
long as $A(u)$ satisfies (\ref{mar718}), we have~$\chi_{FKPP}\leq 1/2$, so that
the semi-FKPP range is `uniformly nontrivial.' 
However, in the pushmi-pullyu case~$\chi=1$ the behavior of the solutions changes drastically, as seen from the second statement in
Theorem~\ref{t.rde-intro}.

A well-known example from~\cite{HadelerRothe}, also discussed in 
detail in~\cite{murray2007mathematical}, 
is the nonlinearity 
\be\label{dec2831}
f(u)=u(1-u)(1+au),
\ee
with $a>0$. This nonlinearity satisfies the FKPP property for all $0\le a\le 1$. In the terminology of the present paper, with $\breta(u) = u - u^2$, it is  
semi-FKPP type in the larger range~$0\le a< 2$ and pushmi-pullyu type when $a=2$. 
A generalization of this example:
\be\label{dec2840}
f(u)=u(1-u^n)(1+au^n),
\ee
was considered in~\cite{Ebert-vanSaarlos}. This corresponds to (\ref{mar722}) with   $\breta(u)=u-u^n$.  To put this in the form of~\eqref{mar722}, we set $A(u) = u^{n+1}$ and $\chi = a/(n+1)$. 
This nonlinearity also has the 
Fisher-KPP property for all~$0\le a\le 1$, as can be seen from (\ref{mar804}), but is of the 
semi-FKPP type in the range~$0\le a< n+1$. It is pushmi-pullyu type when $a = n+1$.

%

\subsection{The traveling wave profile and the nonlinearity}

In order to further explain the particular form of the nonlinearity~\eqref{jul1314}, we need to discuss some basic facts about traveling waves
for semi-linear parabolic equations. 
%
The main result of this section is an expression for the nonlinearity $f(u)$
in terms of the traveling wave profile function $\eta_c(u)$
with speeds $c\ge c_*$, defined by (\ref{dec2812bis}). 
It is elementary,
but is quite interesting in its own right and used frequently in the sequel, so we state this connection 
as a standalone result.

\begin{prop}\label{prop-dec2802}
For each $c\ge c_*$, the function $\eta_c(u)$ is continuously differentiable for $u\in[0,1]$ and satisfies
\be\label{dec2825}
	\eta_c(0)=\eta_c(1)=0,
		~~
	\eta_c'(0)>0,
		~~
	\eta_c(u)>0
	~~ \hbox{and}~~
	\eta_c'(u)<c,
\ee
for all $u \in (0,1)$. 
Moreover, for each $c\ge c_*$,
the function $f(u)$ can be expressed in terms of $\eta_c(u)$ by
\be\label{dec2826}
f(u)=\eta_c(u)(c-\eta_c'(u)),~~\hbox{ for all $u\in(0,1)$.}
\ee
\end{prop}
\medskip\noindent{\bf Proof.} The positivity of $\eta_c(u)$ for $u\in(0,1)$ follows from the aforementioned 
strict negativity of~$U_c'$, which also implies that $\eta_c(0)=\eta_c(1)=0$.
To check the continuous differentiability of~$\eta_c(u)$ we only need to analyze
the behavior near~$u=0$ and~$u=1$; indeed, the case $u \in (0,1)$ follows by the negativity of $U_c'$ and the inverse function theorem.  We consider only the behavior near~$u=0$ as 
the other case can be handled similarly. 
Let $U_c(x)$ be a solution to (\ref{dec2806bis}) with some $c\ge c_*$. 
Recall that 
traveling waves have the asymptotics (\ref{jul1312}) 
\be\label{dec2816}
U_c(x)\sim (D_cx+B_c)e^{-\lambda_cx},
\ee
with $\lambda_c$ being a positive root of 
\be\label{dec2827}
\lambda_c^2-c\lambda_c+f'(0)=0,
\ee
given by
\be\label{dec2828}
\lambda_c=\frac{c\pm \sqrt{c^2-4f'(0)}}{2}.
\ee
Let us make two remarks about~\eqref{dec2816}-\eqref{dec2828}.  First, the `$+$' sign 
in (\ref{dec2828}) appears only in the case~$c = c_*>2\sqrt{f'(0)}$,  
and the `$-$' sign corresponds to all the other cases: either~$c=c_*=2\sqrt{f'(0)}$, or~$c>c_*$; 
see \cite[Proposition~4.4]{AronsonWeinberger}.
Second, the 
coefficient $D_c$ may be non-zero only if~$c=2\sqrt{f'(0)}$, so that $\lambda_c$ is a
double root of (\ref{dec2827}).

We deduce from
(\ref{dec2816}) that
\be\label{dec2818}
\lim_{x\to+\infty}\farc{U_c'(x)}{U_c(x)}=-\lambda_c.
\ee
Using (\ref{dec2818}) and that $\eta_c'(U_c(x)) U_c'(x) = - U_c''(x)$, we obtain
\be\label{dec2820}
\bal
\eta_c'(0)&=-\lim_{x\to+\infty}\frac{U_c''(x)}{U_c'(x)}
=c+\lim_{x\to+\infty}\frac{f(U_c(x))}{U_c'(x)}=
c+\lim_{x\to+\infty}\frac{f'(0)U_c(x)}{U_c'(x)}=c-\farc{f'(0)}{\lambda_c}\\
&=\farc{c\lambda_c-f'(0)}{\lambda_c}=\lambda_c.
\enbal
\ee
We used (\ref{dec2827}) in the last step. 
In particular, it follows from (\ref{dec2820}) that $\eta_c'(0)>0$. To see that~$\eta_c'(u)<c$,
we use the monotonicity of $U_c(x)$ to write
\be\label{dec2835}
\eta'(U_c(x))=-\farc{U_c''(x)}{U_c'(x)}=c+\frac{f(U_c(x))}{U_c'(x)}<c,
\ee

In order to establish (\ref{dec2826}), insert~(\ref{dec2812bis}) into (\ref{dec2806bis}) to find
\be\label{dec2822}
	c\eta(U_c)
		=U_c''+f(U_c)
		=(-\eta(U_c))' + f(U_c)
		=-\eta'(U_c)U_c'+f(U_c)
		=\eta(U_c)\eta'(U_c)+f(U_c).
\ee
As (\ref{dec2822}) holds for all $x\in\Rm$ and $U_c(x)$ is monotonically decreasing and obeys the limits
at infinity in (\ref{dec2806bis}),
this gives expression (\ref{dec2826}) for $f(u)$:
\be\label{dec2823}
f(u)=\eta_c(u)(c-\eta_c'(u)),~~\hbox{ for all $0\le u\le 1$,}
\ee 
finishing the proof.~$\Box$

\medskip
 
One consequence of Proposition~\ref{prop-dec2802} is that
$p(u)=\eta_{c_*}(x)$ is an admissible test function in the Hadeler-Rothe variational principle (\ref{dec2808bis}) 
as it is continuously differentiable.
A simple observation, also going essentially back to~\cite{HadelerRothe} is 
that~$\eta_*(u)$ is actually the optimizer in (\ref{dec2808bis}). 
This is a consequence of (\ref{dec2826}) with $c=c_*$: if~$p(u)=\eta_*(u)$, the expression
inside the $\inf\sup$ in (\ref{dec2808bis}) becomes 
\be\label{dec2824}
	\eta_*'(u)+\farc{f(u)}{\eta_*(u)}=\farc{\eta_*(u)\eta_*'(u)+f(u)}{\eta_*(u)}=
\farc{\eta_*(u)\eta_*'(u)+\eta_*(u)(c_*-\eta_*'(u))}{\eta_*(u)}=c_*.
\ee

A kind of converse to Proposition~\ref{prop-dec2802} is also true. Let $f(u)$ be 
a $C^1([0,1])$ function of the form~(\ref{dec2826}) with $\eta_c(u)\in C^1([0,1])$, and 
consider a traveling wave solution to 
\be\label{dec931}
-cU_c'=U_c''+\eta_c(U_c)(c-\eta_c'(U_c)),~~U_c(-\infty)=1,~U_c(+\infty)=0. 
\ee
We claim that $U_c$ solves the first order ODE (\ref{dec2812bis}):
\be\label{dec930}
-U_c'=\eta_c(U_c),~~~~U_c(-\ifnty)=1,~~U_c(+\ifnty)=0.
\ee
To see this, let $U$ be the unique (up to translation) solution of~\eqref{dec930}.  Then, we have 
\be\label{dec932}
-U''-cU'=\eta'(U)U'-cU'=-\eta_c'(U)\eta_c(U)+c\eta_c(U)=\eta_c(U)(c-\eta_c'(U)),
\ee
which is (\ref{dec931}). Since the traveling wave profiles for both (\ref{dec931}) and (\ref{dec930}) are unique, this shows that, up to translation, $U_c=U$.  
Hence, $U_c$ satisfies~\eqref{dec930}.

Another simple but important comment is that if $f(u)$ has the pushmi-pullyu form (\ref{mar721}), then the minimal speed 
traveling wave profile function is
simply
\be\label{jul2002}
\eta_*(u)=\zeta(u) = u - A(u).
\ee
This property is very convenient in the analysis of the pushmi-pullyu case.

\subsection{Purely exponentially decaying waves}

The next statement characterizes the nonlinearities $f(u)$ for which
the decay is purely exponential, so that $D_c=0$. 
It also gives the asymptotics of $\eta_c(u)$ as~$u\to 0$ 
for waves that have an exponential decay with a linear pre-factor.  In particular, it shows that the difference between pulled and pushmi-pullyu waves can be seen in the regularity of the traveling wave profile function $\eta_c$.

 \begin{prop}\label{dec28-prop06}
Let $c\ge c_*$ and  $f\in C^1([0,1])$. If $\eta_c$ is the traveling wave profile function associated to~\eqref{e.rde} then:

\noindent(i) Suppose, for some $p>1$, there exists $C>0$ so that 
		\be\label{e.c41401}
			\eta_c(u) \sim \lambda_c u + O\Big( \frac{u}{(1+ |\log u|)^p}\Big)
				\qquad\text{ as } u \searrow 0,
		\ee
		where $\lambda_c$ is defined as in~\eqref{dec2828}. 
		In particular, this is true if $\eta_c\in C^{1,\delta}([0,1])$ for some $\delta > 0$.  Then, the profile $U_c$ has purely exponential decay:
		\be\label{e.c41403}
			U_c(x) \sim B_c e^{-\lambda_c x}
				\quad\text{ as } x \to \infty.
		\ee
(ii) If $U_c$ has exponential decay with the linear factor, that is,  $D_c\neq 0$ in~\eqref{dec2816}, then
		\be
			\eta_c(u) \sim \lambda_c u \Big(1 + \frac{1}{\log u} + O\Big(\frac{\log(\log (u^{-1}))}{(\log u)^2}\Big)\Big)
				\qquad\text{ as } u \searrow 0.
		\ee
\end{prop}
\medskip\noindent{\bf Proof of (i).}  Using directly~\eqref{dec2812bis} and then~\eqref{e.c41401}, we find
\be\label{jun903}
-1= \frac{U_c'}{\eta_c(U_c)}
\geq \frac{U_c'}{\lambda_c U_c(1 - C (1+|\log U_c |)^{-p})},
\ee
for some constant $C>0$. 
After possibly shifting the traveling wave, we may assume without loss of generality that the denominator in (\ref{jun903}) never vanishes for $x>0$. 
Multiplying both sides by $\lambda_c$, integrating in $x$, and using the monotonicity of $U_c(x)$ to make a change of variables $z = - \log(U_c)$, yields, for~$x>0$:
\be
\begin{split}
	- \lambda_c x
		&\geq \int_0^x \frac{U_c'}{U_c(1 - C (1+|\log U_c|)^{-p})} dx'
		= \int_{-\log{U_c(0)}}^{-\log  {U_c(x)}} \frac{(-1)}{1 - C (1+ z)^{-p}} dz\\
		&\geq \int_{-\log {U_c(0)}}^{-\log  {U_c(x)}} \Big(-1 + \frac{1}{C (1 + z)^p}\Big) dz
		\geq \log U_c(x)
			- C.
\end{split}
\ee
Hence, we have, for all $x>0$,
\be\label{e.c41402}
	U_c(x) \leq C e^{-\lambda_c x}.
\ee

To refine this bound to the asymptotics in (\ref{e.c41403}), we let $\bar U_c (x) = e^{\lambda_c x} U_c(x)$.    
Using~\eqref{dec2812bis} and~\eqref{e.c41401} again, we find
\be
|\bar U_c'(x)|= e^{\lambda_cx}|\eta_c(U_c)(x) - \lambda_c U_c(x)|
		\leq \frac{C \bar U_c(x)}{(1 + |\log(U_c(x))|)^p}.
\ee
Using \eqref{e.c41402} to bound the numerator and the denominator, we obtain
\be
|\bar U_c'(x)|
		\leq \frac{C}{(1 + x)^p},
\ee
from which the claim~\eqref{e.c41403} follows.~$\Box$ 

\medskip\noindent{\bf Proof of (ii).} 
First, by suitably shifting $U_c$, we may assume $B_c = 0$ without loss of generality.  Then, for $x \gg 1$, we have
\be
	U_c(x) = A_c x e^{-\lambda_c x} + O(e^{-(\lambda_c +\delta) x}),
\ee
for some fixed $\delta>0$.  Next, fix any $u>0$ sufficiently small and let $x$ be such that $u = U_c(x)$.   By~\eqref{dec2812bis}, we have
\be\label{e.c41501}
	\begin{split}
	\eta_c(u)
		&= \eta_c(U_c(x))
		= - U_c'(x)
		= A_c (\lambda_c x - 1) e^{-\lambda_c x}
		+ O(e^{-(\lambda_c + \delta) x})
		\\&
		= \lambda_c U_c(x)
			- \frac{U_c(x)}{x}
			+ O\Big(U_c(x)^{1 + \frac{\delta - \eps}{\lambda_c}}\Big),
	\end{split}
\ee
for any $\eps > 0$.  Notice that
\be
	- \lambda_c x
		= \log U_c - \log(A_c x)
			+ O(e^{-\delta x})
		= \log U_c - \log(-\log(U_c))
			+ O(1).
\ee
Hence,~\eqref{e.c41501} becomes
\be
	\begin{split}
	\eta_c(u)
		&= \lambda_c u + \frac{\lambda_c u}{\log u - \log(-\log u)
			+ O(1)} + O(u^{1+ \frac{\delta-\eps}{\lambda_c}})
		\\&
		= \lambda_c u\Big( 1 + \frac{1}{\log u} + O \Big(\frac{\log( - \log u)}{( \log u)^2}\Big)\Big).
	\end{split}
\ee
This concludes the proof.
~$\Box$

\subsection{When the minimal speed is given by the Fisher-KPP formula}

It is well known that for the Fisher-KPP type nonlinearities the minimal speed is given by the
Fisher-KPP formula
\be\label{dec2829}
c_*[f]=2\sqrt{f'(0)}.
\ee
However, the Fisher-KPP condition (\ref{mar408bis}) is not necessary for (\ref{dec2829}) to hold.
A well-known example of a non-FKPP type nonlinearity that satisfies (\ref{dec2829}), is the nonlinearity of the form~(\ref{dec2831}) in the 
range $1\le a\le2$,
as discussed in detail in~\cite{HadelerRothe,murray2007mathematical}.  This is also true for nonlinearities of the form~(\ref{dec2840})
with $1\le a\le n+1$, considered in~\cite{Ebert-vanSaarlos}.   
A natural question is: for which other nonlinearities does the the Fisher-KPP formula for the speed~\eqref{dec2829} hold?  Below, we give a sufficient condition for this.  
\begin{prop}\label{prop-dec2804}
Assume that $f(u)$ satisfies (\ref{dec2804bis}) and, in addition, that 
there is a $C^1([0,1])$-function~$\breta(u)$ that satisfies 
\be\label{dec2830}
\breta(0)=\breta(1)=0,~~\breta'(0)=1,~~\breta(u)>0,~\breta'(u)<2,~\hbox{ for $0<u<1$,}
\ee
and such that 
\be\label{dec925}
f(u)\le  f'(0)\breta(u)(2-\breta'(u)),~~\hbox{ for all $0\le u\le 1$.}
\ee
Then, the minimal speed $c_*[f]$ is 
\be\label{dec2842}
c_*[f]=2\sqrt{f'(0)}.
\ee 
\end{prop}

\medskip\noindent{\bf Proof.} The first observation is that if $f(u)$ satisfies (\ref{dec2804bis}), 
then we may find a Fisher-KPP type nonlinearity~$f_1(u)\le f(u)$ such that
$f_1'(0)=f'(0)$. As $f_1(u)$ is of the Fisher-KPP type, we have 
\be
c_*[f_1]=2\sqrt{f_1'(0)}=2\sqrt{f'(0)}.
\ee
The comparison principle implies that
\be
c_*[f]\ge c_*[f_1]=2\sqrt{f'(0)}.
\ee

%
%
%
To show that 
\be\label{dec2839}
c_*[f]\le 2\sqrt{f'(0)},
\ee
we note that because of~(\ref{dec2830}), the function $p(u)=\lambda_*\breta(u)$ can be used as test function in the Hadeler-Rothe 
variational principle~(\ref{dec2808bis}).  Here, we have set 
\be\label{mar806}
\lambda_*=\sqrt{f'(0)}.
\ee  
Using assumption~(\ref{dec925}), 
this gives
\be\label{dec922}
c_*[f]\le \sup_{u\in[0,1]}\Big( \lambda_*\breta'(u)+\farc{f(u)}{ \lambda_*\breta(u)}\Big)\le 
\sup_{u\in[0,1]}\Big( \lambda_*\breta'(u)+2\lambda_*-\lambda_*\breta'(u)\Big)=2\sqrt{f'(0)},
\ee
finishing the proof.~$\Box$

\medskip

An important consequence of \Cref{prop-dec2804} is that the Fisher-KPP formula holds for the semi-FKPP and pushmi-pullyu nonlinearities (that is,~\eqref{jul1314} with $\chi < 1$ and $\chi = 1$, respectively)  as  can be seen by taking $\breta (u)= u - A(u)$.  
Notice that the pushmi-pullyu nonlinearities are at the boundary of the validity of our condition for 
the Fisher-KPP formula.

\begin{cor}\label{cor-mar802}
If $f(u)$ is a semi-FKPP or pushmi-pullyu nonlinearity then $c_*[f]=2\sqrt{f'(0)}$. 
\end{cor}
Corollary~\ref{cor-mar802} shows that, at the level of the propagation speed, we do not see a difference
between semi-FKPP and pushmi-pullyu nonlinearities -- both behave similarly to the Fisher-KPP type.

Proposition~\ref{prop-dec2804} 
allows to use explicit 
functions~$\breta(u)$ to verify the validity of the Fisher-KPP formula. 
For example, if we take $\breta(u)=u(1-u)$, then assumption (\ref{dec925}) becomes
\be
f(u)\le u(1-u)(1+2u).
\ee
It holds for nonlinearities of the form (\ref{dec2831}) exactly in the semi-FKPP
range~$0\le a\le 2$. 
On the other hand, for 
\be
\breta(u)= u(1-u^n),
\ee
the assumption (\ref{dec925}) becomes
\be\label{dec2841}
f(u)\le  u(1-u^n)(1+(n+1)u^n).
\ee
Nonlinearities of the form (\ref{dec2840}) satisfy (\ref{dec2841}) 
in the range $0\le a\le n+1$, in agreement with the aforementioned
results in~\cite{Ebert-vanSaarlos}.

\section{Algebraic properties of  semi-FKPP and pushmi-pullyu nonlinearities: the Cauchy problem}\label{sec:alg}

In this section, we discuss some special properties of the semi-FKPP and pushmi-pullyu
nonlinearities. In particular, we introduce the shape defect function,
establish a new variational formulation for reaction-diffusion equations,
and explain a  natural connection between the reaction-diffusion equations  
and reactive conservation laws.

\subsection{Convergence in shape and the shape defect function}\label{s.sdf}

As we have discussed in the introduction, it was proved already in the original KPP paper~\cite{kpp} that the  solution $u(t,x)$   
to (\ref{e.rde}) with the initial condition
$u(0,x)=\one(x\le 0)$ converges in shape to a minimal speed traveling wave, in the sense that (\ref{mar808bis}) holds: 
there exists a reference frame $m(t)$ such that
\be\label{mar808}
|u(t,x+m(t))-U_*(x)|\to 0,~~\hbox{ as $t\to+\infty$, uniformly in $x\in\Rm$.}
\ee
Recall that we have defined in the introduction the traveling wave shape defect function (or shape defect function, for short) as 
\be\label{dec2121}
w(t,x)=-u_x(t,x)-\eta_c(u(t,x)).
\ee
In particular, if $u(t,x)=U_c(t,x)$ then $w(t,x)\equiv 0$ because the traveling wave $U_c(x)$
satisfies (\ref{dec2812bis}). 
However, we also have $w(t,x)\equiv 0$ if $u(t,x)\equiv 0$ or  $u(t,x)\equiv 1$.

Let us note that if $w(t,x)>0$ for all $x\in\Rm$ then $u(t,x)$ is steeper than the traveling wave profile~$U_c(x)$. Here, we use the notion of steepness from~\cite{AHR,GM}: if $u_{1,2}(x)$
are two monotonically decreasing
functions, then $u_1$ is steeper than $u_2$ if for every $u\in\Rm$ such that there exist
$x_{1,2}\in\Rm$ such that $u_1(x_1)=u_2(x_2)=u$, we have $|u_1'(x)|>|u_2'(x)|$.   Actually, a stronger relationship holds: among functions $u$ connecting $1$ at $x=-\infty$ and $0$ at $x=+\infty$, $u$ is as steep as $U_*$ if and only if $w \geq 0$.

A direct computation, using (\ref{e.rde}) and the representation (\ref{dec2826}) for $f(u)$ shows that the shape defect function  satisfies 
\be\label{dec2120}
\bal
&w_t-w_{xx}=-u_{xt}+u_{xxx}-\eta_c'(u)u_t+\eta_c'(u)u_{xx}+\eta_c''(u)u_x^2
\\
&~~=-(\eta_c(u)(c-\eta_c'(u)))_x-\eta_c'(u)(u_{xx}+\eta_c(u)(c-
\eta_c'(u)))+\eta_c'(u)u_{xx}+\eta_c''(u)u_x^2
\\&~~=-\eta_c'(u)(c-\eta_c'(u))u_x
+\eta_c(u)\eta_c''(u)u_x-\eta_c(u)\eta_c'(u)(c-\eta_c'(u))+\eta_c''(u)u_x^2\\
&~~= \eta_c''(u)u_x(u_x+\eta_c(u))-\eta_c'(u)(c-\eta_c'(u))(u_x+\eta_c(u))
=-(\eta_c''(u)u_x-\eta_c'(u)(c-\eta_c'(u)))w\\
&~~=\big(\eta_c''(u)(w+\eta_c(u))+\eta_c'(u)(c-\eta_c'(u))\big)w.
\enbal
\ee
The maximum principle then implies that
\be\label{e.c92201}
	\text{if } w(0,\cdot) \geq 0
	\quad\text{ then }\quad
	w(t,\cdot) \geq 0
	\quad\text{ for all } t>0.
\ee
This is the preservation of steepness property of KPP: if the initial condition $u(0,x)$ is steeper than
a traveling wave, it remains steeper than the wave for all $t>0$. We use 
this property extensively throughout the paper.

\subsubsection*{A variational formulation in terms of the shape defect function}

An interesting observation we have mentioned in the introduction 
is that the shape defect function provides an energy for the reaction-diffusion equation (\ref{e.rde}).
Consider this equation   in the moving frame in view of the identity~\eqref{dec2826}:
\be\label{dec1806}
u_t-cu_x=u_{xx}+\eta_c(u)(c-\eta_c'(u)),
\ee
and define the energy functional
\be\label{dec1802}
\cE_c(u)=\farc{1}{2}\int_\Rm e^{cx}(u_x+\eta_c(u))^2dx.
\ee 
Let us compute
\be\label{dec1804}
\bal
\farc{\delta\cE_c}{\delta u}&=-\pdr{}{x}\Big(e^{cx}(u_x+\eta_c(u))\Big)+e^{cx}(u_x+\eta_c(u))\eta_c'(u)\\
&=e^{cx}\big(-u_{xx}-\eta_c'(u)u_x-cu_x-c\eta_c(u)+u_x\eta_c'(u)+\eta_c(u)\eta_c'(u)\big)\\
&=-e^{cx}\big(u_{xx}+cu_x+\eta_c(u)(c-\eta_c'(u))\big).
 \enbal
\ee
Therefore, equation (\ref{dec1806})  has a variational formulation
\be\label{dec1808}
\pdr{u}{t}=-e^{-cx}\farc{\delta\cE_c}{\delta u}.
\ee
As a consequence, it follows that if $u(t,x)$ is a solution to (\ref{dec1806}), then
\be\label{mar816}
\farc{d\cE_c(t)}{dt}\le 0,
\ee
with a strict inequality unless $u=U_c(x)$, $u\equiv 0$ or $u\equiv 1$.

Other variational formulations for reaction-diffusion equations have been considered; see, 
for example, 
in~\cite{GR,LMN,MN1,MN2,Risler}. The energy functional considered in those papers is
\be\label{mar814}
\tilde\cE_c[u]=\int_\Rm e^{cx}
\Big(\farc{1}{2}u_x^2+F(u)\Big)dx.
\ee
Here, $F$ is the anti-derivative of $-f$: $F'(u)=-f(u)$. One issue with the functional~$\tilde\cE_c[u]$ is that it is not defined
for $u(t,x)=U_c(x)$ unless $2\lambda_c>c$.
This essentially restricts its use to pushed fronts, that is, those that propagate at speed $c=c_*$ where $c_* > 2 \sqrt{f'(0)}$, as these pushed waves have
decay rate given by
\be\label{mar812}
\lambda_c=\farc{c+\sqrt{c^2-4f'(0)}}{2}
	> \frac{c}{2}.
\ee
On the other hand, the functional $\cE_c[u]$ vanishes if $u(t,x)=U_c(x)$  and is thus well-defined,
regardless of the pushed or pulled nature of the traveling wave. 
In addition,  it coincides with $\tilde\cE_c[u]$ for sufficiently rapidly decaying 
solutions. To see this, let us set
\be
N_c(u)=\int_0^u\eta_c(u')du',
\ee
and write
\be
\bal
\cE_c[u]&=\farc{1}{2}\int_\Rm e^{cx}(u_x+\eta_c(u))^2dx=\farc{1}{2}\int_\Rm e^{cx}
(u_x^2+2u_x\eta_c(u)+\eta_c^2(u))dx\\
&=\farc{1}{2}\int_\Rm e^{cx}
(u_x^2+2(N_c(u))_x+\eta_c^2(u))dx=\farc{1}{2}\int_\Rm e^{cx}
(u_x^2-2cN_c(u)+\eta_c^2(u))dx\\
&= \int_\Rm e^{2x}
\Big(\farc{1}{2}u_x^2+V_c(u)\Big)dx,
\enbal
\ee
where, we have defined
\be
V_c(u)= - cN_c(u)+\frac{1}{2}\eta_c^2(u).
\ee
However, $V_c(u)$ is an anti-derivative of $(-f(u))$ because  
\be
V_c'(u)= -c\eta_c(u)+\eta_c(u)\eta_c'(u)=-\eta_c(u)(c-\eta_c'(u))=-f(u).
\ee
This agrees with (\ref{mar814}), so that $\cE_c[u]$ coincides with $\tilde\cE_c[u]$
when both are defined. 

Let us comment that the energy (\ref{dec1802}) can be
generalized in a natural way to 
dimensions~$d>1$. Let $c\in\Rm$ be a speed such that a traveling wave solution to
the one-dimensional problem
(\ref{e.rde}) exists.  Consider the corresponding reaction-diffusion equation in
$\Rm^d$, in a moving 
frame, going in a direction~$e\in{\mathbb S}^{n-1}$, with $\|e\|=1$, at the speed $c$:
\be\label{oct302}
u_t-c e\cdot\nabla u=\Delta u+f(u),~~t>0,~x\in\Rm^d.
\ee
Given  $e\in{\mathbb S}^{n-1}$ and $c$, consider the energy
\be\label{oct304}
\cE_{c,e}[u]=\farc{1}{2}
\int e^{c(e\cdot x)}|\nabla u+e\eta_c(u)|^2dx.
\ee
Then, an essentially verbatim computation to (\ref{dec1804}) shows that
\be\label{oct306}
\bal
\farc{\delta\cE_{c,e}}{\delta u}
	&=-\sum_{k=1}^d\pdr{}{x_k}\Big(e^{c(e\cdot x)} 
\Big(\pdr{u}{x_k}+e_k\eta_c(u)\Big)\Big)+e^{c(e\cdot x)}
\sum_{k=1}^d e_k\Big(\pdr{u}{x_k}+e_k\eta_c(u)\Big)\eta_c'(u)
	\\&
	=e^{c(e\cdot x)}\big(-\Delta u-\eta_c'(u)(e\cdot\nabla u)-ce\cdot\nabla u
-c\eta_c(u)+\eta_c'(u)(e\cdot\nabla u)+\eta_c(u)\eta_c'(u)\big)
	\\&
	=-e^{c(e\cdot x)}\big(\Delta u+c e\cdot\nabla u+\eta_c(u)(c-\eta_c'(u))\big)
	=-e^{c(e\cdot x)}(\Delta u+ce\cdot\nabla u+f(u)).
 \enbal
\ee
It follows that (\ref{oct302})  has a variational formulation
\be\label{oct308}
\pdr{u}{t}=-e^{-c(e\cdot x)}\farc{\delta\cE_{c,e}}{\delta u}.
\ee
In particular, we have
\be\label{oct310}
\farc{d\cE_c(t)}{dt}\le 0.
\ee
A strict inequality holds in (\ref{oct310}) unless $u(t,x)=U_c(x\cdot e-ct)$, 
$u\equiv 0$ or $u\equiv 1$.

We are not going to pursue this variational direction in the present paper, but this approach
seems to make the variational tools of~\cite{GR,Risler} available for a larger class
 than the bistable equations considered in the aforementioned papers.  In particular, it opens the door to a variational analysis of 
Fisher-KPP type equations.

\subsection{Connection to reactive conservation laws for pushmi-pullyu nonlinearities} \label{sec:cons-law-push}

\subsubsection*{The common traveling wave profiles}

We describe a new connection between reaction-diffusion equations and reactive conservation laws provided by the shape defect
function. 
Let us assume that the nonlinearity $f(u)$ is of pushmi-pullyu type (recall~\eqref{jun906}): 
\be\label{jan510}
	f(u)
		= \lambda_*^2 (u-A(u))(1 + A'(u)),
\ee
with some $\lambda_*>0$. Recall that in this case the wave profile function is $\eta_*(u)=\zeta(u)$, as in (\ref{jul2002}).
%
To make the notation less heavy we assume without loss of generality that   
$\lambda_*=1$ and $c_*=2$. 
We also let $U(x)$ be the corresponding minimal speed traveling wave profile, the solution to (\ref{dec2812bis}):
\be\label{jan506}
	-U'= U - A(U),
		\quad U(-\infty)=1,
		\quad U(+\infty)=0,
\ee
and
\be\label{jan514}
-2U'=U''+f(U).
\ee
Let us 
write
\be\label{dec933}
	-2U'+(A(U))'-U''
		=-2U'+U'-(\breta(U))'-U''
		=- U'
		=U-A(U).
\ee
Thus, apart from (\ref{jan514}), the solution to (\ref{jan506}) is also a traveling wave solution to the
reactive conservation law
\be\label{dec934}
u_t+(A(u))_x=u_{xx}+u-A(u).
\ee
In other words, if $f(u)$ is a pushmi-pullyu type nonlinearity, then the reactive conservation 
law~(\ref{dec934}) and
the reaction-diffusion equation  
\be\label{mar818}
u_t=u_{xx}+f(u),
\ee
with 
\be\label{mar820}
	f(u)
		=(u-A(u))(1+A'(u)),
\ee
have exactly the same minimal speed traveling wave profiles.

\subsubsection*{Comparison to reactive conservation laws}

The connection between the reactive conservation law (\ref{dec934}) and
the reaction-diffusion equation   (\ref{mar818}) goes beyond the common traveling wave profile. 
Let $u$ be a solution to~\eqref{mar818}-\eqref{mar820}
and, as usual, assume that $A(u)$ satisfies (\ref{mar718}) and is increasing:
\be\label{dec2110}
A'(u)\ge 0
	\quad\text{ for all $u\in[0,1]$.}
\ee
We claim that if the shape defect function is non-negative:
\be\label{dec2108}
	w(t,x)
		=-u_x(t,x)-\eta_*(u(t,x))
		\ge 0
	\quad\text{ for all $x\in\Rm$ and $t\ge 0$,}
\ee
then
$u(t,x)$ is a subsolution to the reactive conservation law (\ref{dec934}):
\be\label{dec2102}
u_t+(A(u))_x\le u_{xx}+u-A(u).
\ee
Therefore, we can use the comparison
principle to bound the solution to the reaction-diffusion equation (\ref{mar818}) from above by the solution
to the reactive conservation law (\ref{dec934}).  This is used extensively below. 
Let us recall that we have shown that (\ref{dec2108}) holds at all times $t>0$ as long as it is satisfied 
at $t=0$. Thus, (\ref{dec2108}) is at most a restriction on the initial condition. 

To show that (\ref{dec2102}) holds,  
we use (\ref{jul2002}) to write
\be\label{dec2106}
\bal
	u_t&+A'(u)u_x- u_{xx}-u+A(u)
		= (u-A(u))(1 + A'(u)) + A'(u) u_x - u + A(u)
		\\
		&=(u-A(u))(1+A'(u))+A'(u)(-w-\eta_*(u))-u+A(u)
		\\
		&=u-A(u)+A'(u)(u-A(u))-wA'(u)-A'(u)(u-A(u))-u+A(u)
		\\
		&=-wA'(u)
		\le 0,
\enbal
\ee
because of (\ref{dec2110}) and (\ref{dec2108}). 
   
In the general case, if  (\ref{dec2108}) does not hold, so that the shape defect function is not positive everywhere,  
a solution~$u(t,x)$  to (\ref{mar818})-\eqref{mar820}   satisfies the 
forced reactive conservation law
\be\label{dec2114}
u_t+(A(u))_x= u_{xx}+u-A(u)-A'(u)w,~~~~w=-u_x-u+A(u).
\ee

\subsubsection*{The shape defect function for reactive conservation laws}

We note that the shape defect function can be defined for solutions of this reactive conservation law in the same manner.  
For any solution $\urcl$ to~\eqref{dec934}, let
\be\label{jun702}
	w_{\rm rcl}(t,x)
		= - \partial_x\urcl(t,x) - \eta_*(\urcl(t,x)),
\ee
where we have used the `rcl' subscript to distinguish $w_{\rm rcl}$ from the shape defect function for solutions to~\eqref{mar818} and we have changed to the $\partial$ notation in order to avoid the awkward double subscript.  
Let us recall that if $f(u)$ is a pushmi-pullyu nonlinearity, then the minimal speed traveling wave profiles for the reaction-diffusion equation (\ref{mar818})  
and the reactive conservation law 
(\ref{dec934}) are the same. This allows us to use the same wave profile function $\eta_*(u)$ both in the definition (\ref{dec2108}) 
of the shape defect function for the reaction-diffusion equation and in (\ref{jun702}). 
In this case, we find
\be\label{e.wrcl}
	\begin{split}
		\partial_t &w_{\rm rcl}
			- \partial_{x}^2 w_{\rm rcl}
			=
			\left[- \partial_{xt}\urcl - \partial_t \urcl (1- A'(\urcl))\right]
				+ \left[\partial_{x}^3\urcl + \partial_{x}^2\urcl - A(\urcl)_{xx}\right]
			\\&=
			- \partial_{x}^3\urcl + A(\urcl)_{xx} - \partial_x\urcl + A(\urcl)_x
			+ (- \partial_{x}^2\urcl + A(\urcl)_{x}
			\\&\qquad
			- \urcl + A(\urcl))(1 - A'(\urcl)) + \partial^3_{x}\urcl + \partial^2_{x}\urcl - A(\urcl)_{xx}
			\\&=
			- \partial_x\urcl + A(\urcl)_x
			+ (- \partial_{x}\urcl + A(\urcl)_{x} - \urcl + A(\urcl))(1 - A'(\urcl))
			+ \partial^2_{x}\urcl
			\\&=
			- \partial_x\urcl + A(\urcl)_x
			+ \partial^2_{x}\urcl A'(v)
			+ A'(\urcl) \partial_x\urcl (1 - A'(\urcl))
			\\&\qquad
			- \urcl(1 - A'(\urcl))
			+ A(\urcl) (1 - A'(\urcl))
			\\&=
			w_{\rm rcl}(1 - A'(\urcl))
			- A'(\urcl) \partial_x w_{\rm rcl}.
	\end{split}
\ee
While the form of the equation~\eqref{e.wrcl} for $w_{\rm rcl}$ is different from that of~\eqref{dec2120} for $w$, it still preserves positivity. We conclude that
\be
	\text{if }
	w_{\rm rcl}(0,\cdot) \geq 0
	\quad\text{ then }\quad
	w_{\rm rcl}(t,\cdot) \geq 0 \qquad \text{ for all } t\geq 0.
\ee

\subsection{The weighted Hopf-Cole transform} \label{s.hopf_cole}

To finish this section, we introduce the weighted Hopf-Cole transform for pushmi-pullyu and semi-FKPP nonlinearities. 
It generalizes a similar transformation for the Burgers-FKPP equation considered in~\cite{AHR}. 
Let us consider a reaction-diffusion equation in the frame $x\to x-2t$, with a semi-FKPP or pushmi-pullyu
nonlinearity of the form
\be\label{dec1431}
u_t-2u_x=u_{xx}+ f(u),
\ee
where $f$ has the form
\be\label{jun914}
f(u)=\breta(u)(1+\chi A'(u)),
\ee
with  $0\le \chi\le 1$ and $\breta(u)$ related to $A(u)$ by $\breta(u)=u-A(u)$.  
The function $A(u)$ is convex and satisfies the familiar 
assumptions~(\ref{mar718}) and (\ref{dec2110}).  Thus, the nonlinearity
in (\ref{dec1431}) is of semi-FKPP type if~$0\le \chi<1$ and of pushi-pullyu type if $\chi=1$.

As in (\ref{jun916}), we set
\be\label{mar902}
\alpha(u)=\farc{A(u)}{u},
\ee
so that 
\be\label{mar912}
\breta(u)=u(1-\alpha(u)).
\ee
We define the weighted Hopf-Cole transform via
\be\label{mar906}
v(t,x)=\exp\Big(x+\sqrt{\chi}\int_x^\infty \alpha(u(t,y))dy\Big)u(t,x).
\ee
Our goal is to show the following. 
\begin{prop}\label{prop-mar902}
Let $u(t,x)$ be a solution to (\ref{dec1431}), with the above assumptions
on the functions~$\breta(u)$ and~$A(u)$, and $0\le\chi\le 1$. Assume, in addition,
that $\alpha(u)$ is convex and increasing on~$[0,1]$. Suppose also that the shape defect
function satisfies $w(0,x)\ge 0$ for all  
$x\in\Rm$. 
Then, the function $v(t,x)$ defined by (\ref{mar906}) is a subsolution to the heat equation:
\be
v_t-v_{xx}\le 0.
\ee
\end{prop}
\medskip\noindent{\bf Proof.}  We use the notation 
\be
\Gamma(t,x)= x+\sqrt{\chi}\int_x^\infty \alpha(u(t,y))dy
\ee
for short, and utilize the following computations
\be
v_x=e^\Gamma u_x+(1-\sqrt{\chi}\alpha(u))e^\Gamma u,
\ee
and
\be
\bal
v_{xx}=e^\Gamma u_{xx}+2(1-\sqrt{\chi}\alpha(u))e^\Gamma u_x-\sqrt{\chi}\alpha'(u)e^\Gamma uu_x+(1-\sqrt{\chi}\alpha(u))^2e^\Gamma u,
\enbal
\ee
as well as
\be
\bal
v_t&=e^\Gamma u_t+\Big(\sqrt{\chi}\int_x^\infty \alpha'(u(t,y))u_t(t,y)dy\Big)e^\Gamma u\\
&=e^\Gamma u_t+\Big(\sqrt{\chi}\int_x^\infty \alpha'(u(t,y))
(u_{yy}+\zeta(u)(1+\chi A'(u))+2u_y) dy\Big)e^\Gamma u
		\\
		&=e^\Gamma u_t+\sqrt{\chi}e^\Gamma u\Big(- \alpha'(u)u_x- \int_x^\infty \alpha''(u)u_y^2dy-2 \alpha(u)
		 +
		 \int_x^\infty \alpha'(u)\zeta(u)(1+\chi A'(u))dy\Big).
\enbal
\ee  
Here we used that $A\in C^2$.  Next, we write
\be
\bal
	&e^{-\Gamma}\big(v_t-v_{xx}\big)
		=u_t-u_{xx}-2u_x
			+2\sqrt{\chi}\alpha(u) u_x
			+\sqrt{\chi}\alpha'(u)  uu_x-(1-\sqrt{\chi}\alpha(u))^2u
		\\&\qquad
			+\sqrt\chi u
			\Big(
			-\alpha'(u)u_x
			-\int_x^\infty \alpha''(u)u_y^2dy
			-2\alpha(u)
			+\int_x^\infty \alpha'(u)\zeta(u)(1+\chi A'(u))dy
			\Big),
\enbal
\ee
which is
\be\label{dec1504}
\bal
	e^{-\Gamma}\big(v_t-v_{xx}\big)
		&=
			f(u) - u - \chi \alpha(u)^2 u + 2 \sqrt\chi \alpha(u) u_x
		\\&\qquad+\sqrt{\chi} u\Big(- \int_x^\infty \alpha''(u)u_y^2dy +
			 \int_x^\infty \alpha'(u)\zeta(u)(1+\chi A'(u))dy\Big).
\enbal
\ee
By assumption, the function $\alpha(u)$ is increasing and convex. 
In addition, as $w(0,x)\ge 0$ for all~$x\in\Rm$, we
know that $w(t,x)\ge 0$ for all $t>0$ and $x\in\Rm$. 
Positivity of $w(t,x)$ implies that~$u_x(t,x)<0$ for all $t>0$ and $x\in\Rm$.
We also note that, after using \Cref{lem-jun1402},
\be\label{e.c72601}
	\sqrt{\chi}\zeta(u)\le \eta_*(u)=-u_x-w\le -u_x.
\ee

With these ingredients in hand, we can estimate the first integral in the right side of~(\ref{dec1504}) as 
\be\label{dec1506}
\bal
-\sqrt{\chi}& \int_x^\infty \alpha''(u)u_y^2dy
=-\sqrt{\chi}\int_x^\infty \alpha''(u)(-u_y)(-u_y)dy\\
&	\leq -\chi\int_x^\infty \alpha''(u)\zeta(u)(-u_y)dy.
\enbal
\ee
Therefore, after integrating by parts, we have 
\be
\bal		
-\sqrt{\chi}& \int_x^\infty \alpha''(u)u_y^2dy\le - {\chi}\int_x^\infty \alpha''(u)\zeta(u)(-u_y)dy
= -\chi\alpha'(u)\zeta(u) - \chi\int_x^{\infty} \alpha'(u)\zeta'(u) u_y dy.
\enbal
\ee
To estimate the second integral in the right side of~(\ref{dec1504}), we use~\eqref{e.c72601} again:
\be\label{dec1510}
\bal
\sqrt{\chi}\int_x^\infty &\alpha'(u)\zeta(u)(1+\chi A'(u))dy 
\leq -\int_x^{\infty}\alpha'(u)(1+\chi A'(u))u_y dy\\& = 
\alpha(u) - \chi\int_x^{\infty}\alpha'(u)(1-\zeta'(u))
 u_y dy =\alpha(u) + \chi\alpha(u) + \chi\int_x^{\infty} \alpha'(u)\zeta'(u) u_y dy.
 \enbal
\ee
Combining (\ref{dec1506}) and (\ref{dec1510}) gives
\be
\bal
\sqrt{\chi} u\Big(- \int_x^\infty \alpha''(u)u_y^2dy &+
 \int_x^\infty \alpha'(u)\big(\zeta(u)(1+\chi A'(u))\big)dy\Big)\leq \chi\alpha(u)u + \alpha(u)u-\chi\alpha'(u)\zeta(u)u.
\enbal
\ee
Going back to (\ref{dec1504}) and recalling that $\zeta(u)=u(1-\alpha(u))$, and 
\[
f(u) = \zeta(u)(1 + \chi (u \alpha'(u) + \alpha(u))),
\]
we obtain 
\be\label{dec1516}
\bal
	e^{-\Gamma}\big(v_t-v_{xx}\big)
		&\leq f(u) - u - \chi \alpha(u)^2 u + 2 \sqrt \chi \alpha(u) u_x
			+ \chi \alpha(u)u + \alpha(u) u - \chi u \alpha'(u) \zeta(u)
		\\&
		= f(u) - u(1-\alpha(u))+ \chi \alpha(u) u(1 - \alpha(u))
			- \chi u \alpha'(u) \zeta(u)
			+ 2 \sqrt \chi \alpha(u) u_x
		\\&
		= f(u) - \zeta(u) - \chi \alpha(u) \zeta(u)
			+ 2 \chi \alpha(u) \zeta(u) - \chi u\alpha'(u) \zeta(u) + 2 \sqrt \chi \alpha(u) u_x
		\\&
		= 2 \sqrt \chi \alpha(u) (\sqrt \chi \zeta(u) + u_x)
		\leq 0.
\enbal
\ee
We used~\eqref{e.c72601} once again, as well as the positivity of $w(t,x)$ in the last line above.~$\Box$

\section{Discussion of the proofs of Theorems~\ref{t.rde-intro} and \ref{t.rcl-intro}} \label{sec:outline}

\subsection{Generalities}

The first step in the proof of both of the convergence results in Theorem~\ref{t.rde-intro} is to identify the candidates for the
coefficient $r_\chi$ and the shift $x_0$ in first three terms in 
the expansion 
\be\label{jun1002}
m(t)=2t - r_\chi \log t + x_0+o(1),~~\hbox{ as $t\to+\infty$.}
\ee
This is achieved by analyzing the precise tail behavior of the solution at the position $x = m(t) + t^\gamma$, with a small~$\gamma>0$. 
The second step is to argue that, after shifting the traveling wave so that it  matches
the solution $u(t,x)$  at the position
\be\label{jun1004}
x_\gamma(t)=2t - r_\chi \log t + t^\gamma,
\ee
the shifted wave and the solution are close not just at $x_\gamma(t)$ 
but everywhere to the left of it as well. This general strategy has been introduced for the Fisher-KPP nonlinearities in~\cite{Graham,HNRR,NRR1,NRR2},
and is a manifestation of the pulled nature of the propagation: behavior at the front is controlled by the behavior 
far ahead of the front. In the Fisher-KPP case considered in the above references, 
the second step is relatively straightforward,
and relies heavily on the Fisher-KPP property of the nonlinearity. 
More precisely, it allowed to control the sign
of the principal eigenvalue for a certain half-line problem that ultimately induces the convergence. 
This gives a quantitative way to look at the pulled nature of the Fisher-KPP type equations. 
 
As both the semi-FKPP and pushmu-pullyu nonlinearities do not satisfy the Fisher-KPP property, genuinely new ingredients are required in both steps of this approach. 
The types of difficulties that arise, and where they arise, are quite different in the semi-FKPP and pushmi-pullyu cases.  
We now briefly discuss the particulars of each case in greater detail.

\subsection{Semi-FKPP fronts}

\noindent{\bf First step: identifying a candidate expansion (\ref{jun1002}).}
The first step is to identify the coefficient $r_\chi=3/2$ in~(\ref{jun1002}) 
and the constant $x_0$ for semi-FKPP nonlinearities. 
Here, the lack of the Fisher-KPP property of the nonlinearity is compensated by a use of the miracle of the weighted Hopf-Cole transform~\eqref{mar906}. 
While somewhat mysterious, it makes this step relatively straightforward, allowing to use the modification of the 
strategy of~\cite{HNRR,NRR1,NRR2}, introduced in~\cite{AHR}, 
even though the nonlinearity is not of the Fisher-KPP type. 
That is, the weighted Hopf-Cole transform overcomes the absence of the Fisher-KPP
property, to an extent, as it did in~\cite{AHR} for the Burgers-FKPP equation. 

\smallskip

\noindent{\bf Second step: matching the traveling wave far ahead of the front and tracing back.} 
Next, we use the precise information about $u(t,x)$ at the position $x_\gamma(t)$ given 
by (\ref{jun1004}) with $r_\chi=3/2$ obtained in the first step, to match $u(t,x)$ to a shift
of the traveling wave $U_*$ at $x=x_\gamma(t)$.  
The argument of~\cite{HNRR,NRR1,NRR2} proceeds by taking the difference $s=u-U_*$ and using the 
Fisher-KPP property of the nonlinearity to show that $s(t,x)$ satisfies a Dirichlet problem
with a time-decaying solution. The time decay comes about 
because the zero-order term has a good sign and $x_\gamma(t)$ is not 
too large. In~ the Burgers-FKPP situation of~\cite{AHR}, this was adapted to the  
difference $s= \tilde u - \psi$, where~$\tilde u$ 
is the appropriate shift of the weighted Hopf-Cole transform of $u$ and $\psi$ is the weighted Hopf-Cole transform of $U_*$. This modification still relied heavily on the fact that the nonlinearity in the Burgers-FKPP equation   is of the
Fisher-KPP type. 
Here, this `Hopf-Cole modified' argument fails, 
due to a term that does not have a sign because of the lack of the Fisher-KPP property. 
%
However, by a series of intricate manipulations, surprisingly, we manage to extract 
a term that is sufficiently positive, as long as the nonlinearity is precisely of the semi-FKPP type.
With this, 
we obtain an upper bound for $s$ that tends to zero. The main novelty here is in the proofs of Lemmas~\ref{l.sfkpp_s} and~\ref{l.s_eqn} that require extremely delicate cancellations. These computations quantify the pulled nature of semi-FKPP fronts.

The proof of \Cref{t.rde-intro}.(i)  is presented in \Cref{sec:sfkpp} with the first step contained in \Cref{sec:sfkpp_tgamma} and the second step in \Cref{sec:sfkpp_conv}.

\subsection{Pushmi-pullyu fronts}

{\bf First step: identifying the front location to precision $O(1)$.}
For the pushmi-pullyu fronts, 
we need to separate the identification of the constants $r_\chi$ and $x_0$ in~(\ref{jun1002})
into two steps; that is, we first obtain the asymptotics of the front location as 
\be\label{jun1008}
m(t)=2t - \tfrac{1}{2}\log t + O(1),~~\hbox{as $t\to+\infty$}
\ee
and then we bootstrap (\ref{jun1008}) to find the precise $O(1)$ term.  
%

%
%
%
%
For the proof of (\ref{jun1008}), the connection between the reaction-diffusion equation and the reactive conservation law is crucial.  By~\eqref{dec2102}, it is enough to obtain a lower bound on the reaction diffusion equation and an upper bound on the reactive conservation law.  The lower bound 
in (\ref{jun1008}) can be obtained via a suitable estimate on an exponential moment, inspired by
the Fabes-Stroock proof of the heat kernel bounds in~\cite{FS}.

For the upper bound, we apply the methods developed in~\cite{AHR}, which we briefly outline.  First, after passing to the moving frame $x \mapsto 2t + x$, a change of function $p = e^x u$ is made, resulting in a seemingly simple inhomogeneous conservation law
\be\label{e.c41601}
	p_t + (\alpha(u) p)_x = p_{xx},
\ee
where we recall that $\alpha(u) = A(u)/u$.  At this point, it is enough to show that
\be
	\|p\|_\infty \leq O(1/\sqrt t)
\ee
in order to conclude the front is behind $2t - (1/2) \log t$ due to the definition of $p$.  When $A(u) = u^2$, one can show that $\rho = 1-u$ is a supersolution to~\eqref{e.c41601} so that the relative entropy methods of~\cite{Const, MMP} can be extended and applied to find
\be
	\frac{d}{dt} \int \left( \frac{p}{\rho} \right)^2 \rho\,  dx
		\leq - C \int \left(\frac{p}{\rho} \right)_x^2 \rho \, dx.
\ee
A suitable Nash-type inequality with a dynamic weight is developed in \cite{AHR} for the measure $d\mu_\rho = \rho\, dx$ allowing us to obtain~$O(1/t^{1/4})$-decay of the $L^2_{\mu_\rho}$-norm of $p/\rho$.  Bootstrapping this to $O(1/\sqrt t)$-decay of the $L^\infty$-norm of $p$ requires an intricate argument in which one finds a ``good'' time $t_{\rm good} \approx t$ where  
\[
\|p(t_{\rm good},\cdot)\|_\infty \leq O(1/\sqrt t),
\]
and then ``trapping'' that norm using that~\eqref{e.c41601} conserves mass and enjoys a comparison principle.

In the context where $A(u) \neq u^2$, one must first determine an appropriate $\rho$ that is a supersolution 
to~(\ref{e.c41601}) and  satisfies the assumptions of the aforementioned Nash-type inequality.  As we detail below, the appropriate choice, via yet another 
surprising piece of algebra, turns out to be 
\be\label{e.c41602}
	\rho(t,x)
		= \exp \Big( - \int_0^{u(t,x)} \frac{\alpha(u')}{u'-A(u')}du'\Big). 
\ee
One can check that, in the special case
\be
	A(u) = u \alpha(u) = u^2,
\ee
the formula~\eqref{e.c41602} simplifies to $\rho=1-u$, while in the case $A(u) = u^n$ (recall~\eqref{dec2840}, which was introduced in~\cite{Ebert-vanSaarlos}),
it has the form
\be
\rho
=(1-u^{n-1})^{1/(n-1)}. 
\ee
However, it does not, in general, have such a simple form.

One can then check that such $\rho(t,x)$ satisfies the assumptions of the Nash-type inequality \cite[Proposition~5.9]{AHR} associated to the measure $d\mu_\rho = \rho\, dx$ in a uniform way (that is, independent of $t$).  From here, the proof proceeds as in~\cite{AHR}.
The proof of this step is contained in \Cref{sec:pp_first}.

\smallskip

\noindent{\bf Second step: identifying the front location to precision $o(1)$.} 
 The approach of~\cite{AHR}  to identifying the constant term $x_0$ in $m(t)$ is to apply the weighted Hopf-Cole transform, change to self-similar variables, and show convergence to a constant multiple of the principal eigenfunction of the resulting operator.  This proceeds in a straightforward way for the solution to the reaction-diffusion equation~\eqref{e.rde} because the resulting function $v$ is a subsolution to the heat equation (see \Cref{prop-mar902}).  Unfortunately, this is not true for the solution to the reactive conservation law~\eqref{e.rcl}; indeed, passing to the moving frame $x \mapsto 2t + x$ and letting
\be
	v(t,x) = \exp\Big(x + \int_x^\infty \alpha(u(t,y)) dy\Big)u(t,x),
\ee
we have  
\be\label{jul708}
	v_t - v_{xx}
		\leq \frac{v}{u} (u\alpha'(u)-\alpha(u) ) w_{\rm rcl}.
\ee
We note that the Burgers-FKPP case considered in \cite{AHR}, in which $\alpha(u) = u$, is the unique case in which the right hand side vanishes.  Changing to self-similar variables $t = e^\tau$, $y = e^{\tau/2}x$, the right hand side accumulates an $e^\tau$ multiplicative factor, making the right hand side not only positive but, potentially, large.
The proof is saved, however, by showing that 
\be\label{jul706}
w _{\rm rcl}\leq \farc{Cxe^{-x}}{t},
\ee
in Lemma~\ref{l.w_bounds}, controlling the right side of (\ref{jul708}) by, in self-similar coordinates, by a bounded function that tends to zero in $L^p(0,\infty)$ for all $1\leq p < +\infty$.  
The proof of this bound is achieved by the construction of a somewhat complicated supersolution (see \Cref{l.w_bounds}).

These estimates provide a sharp estimate of the behavior of $u$ at $x = O(t^\gamma)$ for $\gamma > 0$.    The precise shift of the front $x_0$ is then identified via:
\be\label{e.c70701}
U_*(- x_0 + t^\gamma)= u\big(t, 2t - \tfrac12 \log t + t^\gamma\big)
+ o(e^{-t^\gamma}).
\ee
The step is treated  in \Cref{sec:pp_second}.

\smallskip

\noindent {\bf Third step: using the pulled nature of the problem.}  To upgrade the sharp 
estimate~\eqref{e.c70701} at the position $x\sim t^\gamma$ 
to full convergence 
we can use the smallness and decay of $w$. 
This is done as follows. Let $\delta(x) = u(t,x+m(t)) - U(x)$ and notice that
\be\label{e.c41701}
	- \delta_x - \delta \eta_*'(u)
		\approx - \delta_x - \delta \frac{\eta_*(u(t,\cdot+m(t))) - \eta_*(U)}{u(t,\cdot+m(t)) - U}
		= w.
\ee
Thus, ahead of the front we have 
\be
	-\delta_x \approx\delta \eta'(0) + w = \delta + w,
\ee
which can be approximately solved:
\be\label{jul710}
	e^x \delta(x) \approx \delta(t^\gamma) e^{t^\gamma} + \int_x^{t^\gamma} e^y w(y) dy.
\ee
The first term on the right is small due to~\eqref{e.c70701}.  This reflects the pulled nature 
of the pushmi-pullyu case: $u$ and $U$ must be ``very'' close at far ahead of the 
front (at $x=m(t) + t^\gamma$) in order to be close everywhere.  The second term is small 
due to~\eqref{jul706}.  This step is contained in \Cref{sec:pp_third}.

\section{Pushmi-pullyu fronts: proofs of Theorems~\ref{t.rde-intro}.(ii) and~\ref{t.rcl-intro}}\label{sec:pushmi-pf}

\subsection{Outline of the proofs and notation} 

In this section, we prove part (ii) of Theorems~\ref{t.rde-intro} and~\ref{t.rcl-intro}, 
both of which concern the pushmi-pullyu fronts. We work simultaneously with the solutions $\urd$ 
and $\urcl$ to the reaction-diffusion equation~\eqref{e.rde} and the reactive conservation law~\eqref{e.rcl}, respectively.  Due to the comparison~\eqref{dec2102}, we have immediately that
\be\label{e.c41702}
	\urd \leq \urcl,
\ee
as long as $\urd(0,x)=\urcl(0,x)$ for all $x\in\Rm$.  There is, thus, 
some efficiency in considering both at the same time.
The notation $\urd$ and $\urcl$ is somewhat overwrought, so  
we use simply $u$ when it is possible to do so without risk of confusion or when the argument applies to 
both~$\urd$ and $\urcl$.  We also use this subscript notation  for derived quantities, such as
the shape defect functions $w_{\rm rd}$ and~$w_{\rm rcl}$.

Additionally, we need 
to work in two different moving frames.    
We use $\hat u$ to denote $u$ in the~$x \mapsto 2t + x$ moving frame and $\tilde u$ for $u$ in the $x \mapsto 2t - (1/2)\log (t+1) + x$ moving frame:
\be\label{e.moving_frame}
	\hat u(t,x)
		= u(t, x+ 2 t)
	\quad\text{ and }\quad
	\tilde u(t,x)
		= u\big(t, x + 2t - \tfrac{1}{2}\log(t+1)\big).
\ee
To formalize the idea of the `location of the front', we define the position
$(-\mu(t))$ by:
\be\label{e.normalization}
	\alpha(\tilde u(t,-\mu(t))) = \frac12.
\ee
This specific normalization is chosen just for technical convenience (indeed, as $w \geq 0$, any two level sets remain a bounded distance away from each other at all times). 
Since $\urd \leq \urcl$, and the function~$\alpha(u)$ is increasing, we immediately get
\be
	\mu_{\rm rcl} \leq \mu_{\rm rd}.
\ee
Let us recall that Proposition~\ref{prop-dec2802} implies that in the pushmi-pullyu case we have
\be\label{jul702}
\eta_*(u)=\zeta(u)=u(1-\alpha(u)),
\ee
a relation that we use repeatedly in this section.

The proof proceeds as follows. First, we establish the location of the front
with the precision~$O(1)$, in Section~\ref{sec:pp_first}. The main result in 
that section is Proposition~\ref{prop.pp_O(1)}. At the heart of the proof is a relative
entropy computation with respect to a supersolution and a dynamically weighted
Nash inequality. Both of them are originating in the analysis of~\cite{AHR} for the Burgers-FKPP equation
but the details
in the general case are very different and involve some additional fortunate pieces of algebra
in the construction of the supersolution in Lemma~\ref{lem.rho}. 
 This proof occupies the bulk of this section.

The $O(1)$ term in the front location is identified in Section~\ref{sec:pp_second}.
Its main result is Proposition~\ref{prop.pp_o(1)}. Its proof for the reaction-diffusion case
relies on the weighted Hopf-Cole transform in Proposition~\ref{prop-mar902}. For the reactive
conservation law, the weighted Hopf-Cole transform gives not a subsolution to the heat equation
but only an approximate subsolution far on the right. 
This requires additional estimates on the error term.
This part of the argument is discussed
in Section~\ref{sec:rcl-hopf}.

The final step in the proof of Theorems~\ref{t.rde-intro}.(ii) 
and~\ref{t.rcl-intro},
convergence to a single wave, is presented in Section~\ref{sec:pp_third}. It is very
also different from the corresponding step in~\cite{AHR,NRR1,NRR2} and uses the algebra of the
pushmi-pullyu case.

\subsection[The first step: the expansion]
{The first step: the expansion $m(t) = 2t - \frac{1}{2}\protect{\rm log}(t) + O(1)$ of the front location}
\label{sec:pp_first}

As a preliminary observation, we note two bounds on $\tilde{u}$ viewed from the location $(-\mu(t))$
that follow immediately from the positivity of the shape defect function $w$:
\be\label{e.u_asymptotics}
	1 - C e^{\alpha'(1) x}
		\leq \tilde u(t,x - \mu(t))
		\leq C e^{-x}
			\qquad\text{ for all } x\in \R.
\ee
To see this, we argue as follows. Let us normalize a translate of the minimal speed wave so  that~$\alpha(U_*(0)) = 1/2$, and set
\[
s(t,x) = \tilde u(t,x) - U_*(x+\mu(t)).
\]
Notice that
\be\label{e.c42604}
	-s_x
		= \xi s
			+ w
	\quad\text{ and }\quad s(t,-\mu(t)) = 0,
\ee
with
\be
	\xi(t,x) =\frac{\eta_*(\tilde u(t,x)) - \eta_*(U_*(x + \mu(t)))}{\tilde u(t,x) - U_*(x+\mu(t))}.
\ee
In view of~\eqref{e.c42604} and the positivity of $w$, we have 
\be
	s > 0
		\quad\text{ for } x < -\mu(t)
	\qquad\text{and}\qquad
	s < 0
		\quad\text{ for }  x > -\mu(t). 
\ee
From here,~\eqref{e.u_asymptotics} follows directly from the fact that 
\be\label{e.U_asymptotics}
	1 - C e^{\alpha'(1) x}
		\leq U_*(x)
		\leq C e^{-x}
			\qquad\text{ for all } x\in \R,
\ee
because
\be
	- U_*' = \eta_*(U_*) = U_*(1 - \alpha(U_*)),
\ee
as seen from (\ref{jul702}). 
We note also that the inequalities in~\eqref{e.U_asymptotics} are sharp in the sense that
\be\label{e.U_asymptotics_precise}
	\lim_{x\to\infty} e^x U_*(x)
		\quad\text{ and }\quad
	\lim_{x\to-\infty} e^{-\alpha'(1) x} (1 - U_*(x))
	\qquad\text{ exist and are positive.}
\ee

The goal of this section is to prove the following:
\begin{prop}\label{prop.pp_O(1)}
There exist $\underbar L\le\bar L$ and $T_0>0$ sufficiently large, so that for all $t>T_0$ we have 
\be\label{jun1302}
	\btilde\urd(t, \bar L)\leq \btilde \urcl(t, \bar L)
		\,\stackrel{\mbox{\tiny$(i)$}}{\leq} \,
			\alpha^{-1}\Big(\frac{1}{2}\Big)
		\,\stackrel{\mbox{\tiny$(ii)$}}{\leq}\,
			\btilde\urd(t,  \underbar L)\le\btilde\urcl(t,\underbar L).
\ee
As a consequence, we have the following bounds on $\mu(t)$:
\be
	\underbar L
		\leq \mu_{\rm rcl}(t)
		\leq \mu_{\rm rd}(t)
		\leq \bar L,~~\hbox{for all $t>T_0$.}
\ee
\end{prop}
We begin with the inequality (i), that is, the upper bound on $\btilde\urcl$, 
as it is required for the proof of (ii), the lower bound on $\btilde\urd$.  
Due to the  assumption (\ref{mar1616})  on the initial condition and the comparison principle, we may assume, without loss of generality, that
\be
	u_0(x)
		= \one(x < 0).
\ee

\subsubsection{Proof of \Cref{prop.pp_O(1)}.(i): the upper bound on $\btilde\urcl$}

In this section, we work only with $\urcl$. As there is no possibility of confusion, we simply refer to~$\urcl$ as $u$ below.

Changing to the moving frame $x\to x+2t$, yields
\begin{align}\label{criticalconser}
\hat u_t -2\hat u_x+ (\alpha(\hat u) \hat u)_x = \hat u_{xx}+\hat u(1-\alpha(\hat u)).
\end{align}
Making the change of variable
\be
p(t,x) = e^x \hat u(t,x),
\ee
we see that to prove part (i) of  (\ref{jun1302}), it  is enough to show that
\be\label{e.c41804}
	\|p(t,\cdot)\|_\infty \leq \frac{C}{\sqrt t}.
\ee
This estimate is the key for the proof of both the upper and the lower bound in \Cref{prop.pp_O(1)}.

The function $p$ solves the (inhomogeneous) scalar conservation law
\begin{align}\label{oct141}
p_t + (\alpha(\hat u)p)_x = p_{xx},
\end{align}
with the initial condition $p(0,x) = e^x u_0(x)$. 
It is clear that~\eqref{oct141} conserves mass:
\be
\int p(t,x)dx = \int p(0,x)dx.
\ee

This is a general form of the structure used in~\cite{AHR} in order 
to prove~\eqref{e.c41804} in the special case~$A(u)=u^2$,~$\alpha(u)=u$,
which was considered there.  
While the outline of that proof applies here, there are substantial changes that need to be made in order to suitably generalize it to the case considered here.  
We outline the main steps, making note of the major differences.  To make the outline 
easier to follow, we use the same notation as in~\cite{AHR} as much as possible.

\subsubsection*{\bf Step one: the relative entropy calculation}

The basis of the proof in \cite{AHR} is a generalization of 
the relative entropy ideas introduced in~\cite{MMP} (see also~\cite{Const})
to supersolutions.  Suppose that a function $\rho$ satisfies
\be\label{e.c41806bis}
	\rho_t + (\alpha(\hat u) \rho)_x \geq \rho_{xx},
\ee
and let
\be
	\varphi(t,x)
		= \frac{p(t,x)}{\rho(t,x)}.
\ee
Then, \cite[Proposition~5.8]{AHR} gives a dissipation inequality
\be\label{e.c41805}
	\frac{d}{dt} \int \varphi^2(t,x) \rho(t,x) dx
		\leq - 2 \int \varphi_x^2(t,x) \rho(t,x) dx,
\ee
that holds 
for any initial condition $p_0 = e^x u_0$ for which the quantities above are finite.

Let us  point out that, in the special case of the heat equation, 
where  $\alpha \equiv 0$, and we can take~$\rho \equiv 1$, the inequality~\eqref{e.c41805}, 
paired with the Nash inequality, gives us the $O(t^{1/4})$-decay of the $L^2$-norm.  Our approach is analogous, although we appeal to the weighted Nash inequality with time-dependent weights 
introduced in~\cite{AHR}.

In order to use~\eqref{e.c41805}, we need to find a supersolution $\rho$. 
This is provided by the following lemma, whose proof is postponed to the end of this section.  
Let us define an auxiliary function
\be\label{e.F}
	F(u)
		= \exp \Big( - \int_0^u \frac{\alpha(u')}{\eta_*(u')} \, du'\Big)
		= \exp \Big( - \int_0^u \frac{u' - \zeta(u')}{u'\eta_*(u')} \, du'\Big).
\ee

\begin{lem}\label{lem.rho}
The function $F$ defined in~\eqref{e.F} is $C^2_{\rm loc}([0,1))$, satisfies $F(0) = 1$ and 
$F(1) = 0$, and is decreasing and concave.
Moreover,   for any $\eps>0$, there exists $C_\eps$ such that
\be\label{jul716}
\frac{1}{C_\eps} (1 - \alpha(u))^\frac{1+\eps}{\alpha'(1)}
			\leq F(u)
			\leq C_\eps (1 - \alpha(u))^\frac{1}{(1+\eps)\alpha'(1)}
			\qquad\text{ for all } u \in [0, 1].
\ee
Moreover, $\rho(t,x) := F(\hat u(t,x))$ satisfies~\eqref{e.c41806bis} as long as $\hat w \geq 0$ and $\hat u$ satisfies~\eqref{criticalconser}.
\end{lem}

\subsubsection*{\bf Step two: modification of the initial condition}

An important issue with the step-function initial condition $u(0,x)=\one(x\le 0)$, 
that was also present in~\cite{AHR}, is that, at the time $t=0$, 
\[
\rho(0,x)=F(u(0,x))=\one(x>0),
\]
so that
\be
	\int \varphi^2(0,x) \rho(0,x) \dx
		= \int \Big(\frac{p(0,x)}{\rho(0,x)}\Big)^2 \rho(0,x) \dx
		= \int \frac{{e^{2x}}\1(x\leq 0)}{\1(x > 0)} \dx
		= +\infty.
\ee
Hence, the differential inequality~\eqref{e.c41805} cannot be 
used directly with such initial condition.

We address this by modifying the initial data $\hat u(0,\cdot)$ to not take the value $1$ anywhere, 
so that~$\rho(0,\cdot)$
does not vanish anywhere.  It it is important that the new initial condition 
is steeper than the traveling wave, so that the corresponding
 shape defect function $w$ remains non-negative.
While in~\cite{AHR} it was possible to simply write an explicit approximation of~$\1(x\leq 0)$, here, 
we use the traveling wave: for fixed~$\gamma \in (1, 4/3)$ and any~$a > 0$, we let
\be
	\hat u_a(0, x)
		= U_*(\gamma(x - a)) \1(x \leq 0).
\ee
While $\gamma$ remains fixed throughout the proof, we
eventually take the limit $a \to \infty$.  As a result, all constants  depend on $\gamma$, 
but it is important to  
track the dependence on $a$ throughout.

One can check that 
\[
\hat w_a(t,x) = - \partial_x \hat u_a(t,x) - \eta_*(\hat u_a(t,x))
\]
is positive. Indeed, by~\eqref{dec2120}, it is enough to show that $\hat w_a(0,x) \geq 0$, which follows by:
\be
\begin{split}
	\hat w_a(0,x)
		&= - \partial_x \hat u_a(0,x) - \eta_*(\hat u_a(0,x))
		\\&
		= \left(- \gamma U'_*(\gamma(x-a)) - \eta_*(U_*(\gamma(x-a)))\right) \1(x < 0)
			+ U_*(-\gamma a) \delta_0(x)
		\\&
		\geq \left(- \gamma U'_*(\gamma(x-a)) - \eta_*(U_*(\gamma(x-a)))\right) \1(x < 0)
		\\&
		= - (\gamma - 1) U_*'(\gamma(x-a))\1(x < 0)
		\geq 0. 
	\end{split}
\ee
Above we used that $\gamma >1$ and that $U_*$ is decreasing, positive, and  satisfies $- U_*' = \eta_*(U_*)$.

With $\hat u_a$, we now define
\be
	p_a(t,x) = e^x \hat u_a(t,x),
		\quad
	\rho_a(t,x) = F(\hat u_a(t,x)),
		\quad\text{ and }\quad
	\varphi_a(t,x) = \frac{p_a(t,x)}{\rho_a(t,x)}.
\ee
Here, $F$ is defined by~\eqref{e.F}.  
Using \Cref{lem.rho} with the choice $\eps=1/2$ and then~\eqref{e.U_asymptotics} yields
\be\label{e.c41907}
	\begin{split}
		\int \varphi_a^2(0,x) &\rho_a(0,x)\dx
			= \int \Big(\frac{p_a(0,x)}{\rho_a(0,x)}\Big)^2 \rho_a(0,x) \dx
			\leq C\int_{-\infty}^0 
\frac{e^{2x}\, U_*(\gamma(x-a))^2}{(1 - \alpha(U_*(\gamma(x-a)))^\tinyfrac{3}{2\alpha'(1)}} \dx
		\\&
			\leq C\int_{-\infty}^0 \frac{e^{2x}}{(1 - U_*(\gamma(x-a))^\tinyfrac{3}{2\alpha'(1)}} \dx
			\leq C \int_{-\infty}^0 e^{(2 - \tinyfrac{3\gamma}{2}) x + \tinyfrac{3\gamma a}{2}} \dx
			\leq C e^\tinyfrac{3\gamma a}{2}.
	\end{split}
\ee
In the last step, we used that $3\gamma / 2 < 2$, by assumption.  Thus, due to the finiteness of the quantity in~\eqref{e.c41907}, we have the differential inequality~\eqref{e.c41805} at our disposal.  This is crucial in the following steps.

Before proceeding, we discuss how a bound on $p_a$ yields a bound on $p$.  Let us define 
\be
	h = p - p_a
	\quad\text{ and }\quad
	v = \alpha( \hat u)
		+ \hat u_a \frac{\alpha(\hat u) - \alpha(\hat u_a)}{\hat u - \hat u_a} .
\ee
It is straightforward to check that $h$ is nonnegative at $t=0$ and satisfies, for all $t>0$ and $x\in\Rm$,
\be
	h_t + (v h)_x = h_{xx}.
\ee
This equation conserves mass and preserves nonnegativity.  Thus, $h \geq 0$ and
\be
	\begin{split}
	\int h(t,x) \dx
		&= \int h(0,x) \dx
		= \int e^x \left(\hat u(0,x) - \hat u_a(0,x)\right) \dx
		= \int_{-\infty}^0 e^x \left(1 - U_*(\gamma(x-a))\right) \dx
		\\&
		\leq C \int_{-\infty}^0 e^{(1 + \gamma \alpha'(1))x - \gamma \alpha'(1) a} \dx
		\leq  C e^{-\gamma \alpha'(1) a}.
	\end{split}
\ee
In the second-to-last inequality, we applied~\eqref{e.U_asymptotics} again.  By parabolic regularity theory, it is easy to see that, for $t\geq 0$, we have 
a uniform approximation
\be\label{e.c41904}
0 \leq p(t,x) - p_a(t,x)\leq C e^{-\gamma \alpha'(1) a}.
\ee

%
%

\subsubsection*{Step three: weighted $L^2$-decay of $\varphi_a$}

From the work in the previous two steps, we have the dissipation inequality
\be\label{e.c41906}
	\frac{d}{dt} \int \varphi_a(t,x)^2 \rho_a(t,x) \dx
		\leq - 2 \int (\partial_x \varphi_a(t,x))^2 \rho_a(t,x) \dx.
\ee
We need a Nash-type inequality in order to proceed.  This is given by \cite[Proposition~5.9]{AHR}:
for any~$\theta>0$ and  any smooth non-negative function~$\vphi(x)$ that is sufficiently rapidly decaying as $x\to+\infty$ and bounded as $x\to-\infty$, 
we have
\be\label{e.Nashish}
	\int \vphi^2(x)\rho_a(t,x)\, dx
		\le \frac{2}{\theta} \Big( \int \vphi(x)\rho_a(t,x) \dx\Big)^2
			+ 8C_1\max\{1,\theta^2\} \int |\vphi_x(x)|^2 \rho_a(t,x) \dx.
\ee
For this inequality to apply, the function $\rho_a(t,x)$ needs to be positive, 
bounded and increasing, with the left limit~$\rho_a(t,-\ifnty)=0$. In addition, it should
satisfy
\be\label{e.c41905}
	\doverline\rho_a(t,x)
		\leq C_1 \max\{1,\overline\rho_a^2(t,x)\}\rho_a(t,x).
\ee
Here, we use the notation that, for any $r: \R\to\R$,
\be
	\overline r(x) = \int_{-\infty}^x r(x) \dx.
\ee
\begin{lem}\label{lem.Nash_condition}
Under the setting above, $\rho_a(t,x)$ satisfies~\eqref{e.c41905} with $C_1$ independent of $a$ and~$t$.
\end{lem}
We postpone the proof of this lemma for the moment.
We may now proceed nearly verbatim as in the proof of~\cite[Lemma~5.7]{AHR}  and 
use the dissipation inequality~\eqref{e.c41906} together with (\ref{e.Nashish})
to conclude that
\be
	\int \varphi_a(t,x)^2 \rho_a(t,x) \dx
		\leq \frac{C}{\sqrt t},
		~~\text{ for all } t \geq C\int \varphi_a(0,x)^2 \rho_a(0,x) \dx.
\ee
Then, in view of~\eqref{e.c41907}, we see that
\be\label{e.c41908}
	\int \varphi_a(t,x)^2 \rho_a(t,x) \dx
		\leq \frac{C}{\sqrt t}
		\qquad\text{ for all } t \geq T_a=C e^\frac{3\gamma a}{2}.
\ee

\subsubsection*{Step four: bootstrapping from $L^2_{\rho}$-decay to $L^\infty$-decay}

Fix any $T \geq 4T_a$ as in~\eqref{e.c41908}.  The proof in~\cite{AHR} proceeds then by finding a good time $T_g \in [T/2, 3T/4]$ so that
\be\label{e.c41912}
	(i)~ \max_{x \geq - \mu_a(T_g)} \varphi_a(T_g, x) \leq \frac{C}{\sqrt T}
	\quad\text{ and }\quad
	(ii)~ \mu_a(T_g) \geq \frac{1}{2} \log(T) - C.
\ee
To obtain (i), the arguments of~\cite[Lemma~5.10]{AHR} require only the following ingredients: (1) the dissipation inequality~\eqref{e.c41906}; (2) the $L_\rho^2$-decay given by~\eqref{e.c41908}, (3) the mass conservation of~\eqref{oct141}; and~(4) that $\rho_a(t,x) \in [C^{-1}, 1]$ for all $x > - \mu_a$.  All four ingredients are present here, so the proof can be repeated nearly verbatim.  We omit the details and assert~\eqref{e.c41912}.

We now argue that the bounds in~\eqref{e.c41912} can be, essentially, preserved until time $T$.  To this end, fix $K>0$ to be chosen.
An easy computation shows that
\be
	P(x) = e^x U_*\Big(x + \frac{1}{2}\log(K^2T)\Big)
\ee
is a steady solution to
\be
	P_t + (\alpha( e^{-x} P) P)_x
		= P_{xx}.
\ee
Recalling~\eqref{e.U_asymptotics}, that is, that the minimal speed traveling wave 
$U_*$ has a purely exponential decay, we have
\be\label{e.c41914}
	P(x)\leq \frac{C}{K\sqrt T}.
\ee

Let us decompose $p_a$ into its P-part and its error part:
\be
p_a(t,x) = \psi_P(t,x) + \psi_E(t,x).
\ee
Here, $\psi_P$ is the solution to 
\be
\partial_t \psi_P(t,x) + (\alpha(e^{-x} \psi_P) \psi_P)_x
= \partial_x^2 \psi_P,
~~\text{ for all } t\in (T_g, T],
\ee
with initial condition
\be
\psi_P(T_g,x)= \min\{p_a(T_g,x), P(x)\},
\ee
and $\psi_E$ is the solution to 
\be\label{e.c41913}
	\partial_t \psi_E(t,x) + (v \psi_E)_x
		= \partial_x^2 \psi_E,~~\text{ for all } t\in (T_g, T],
\ee
with initial condition
\be
\psi_E(T_g,x)
= p_a(T_g,x) - \min\{p_a(T_g,x), P(x)\},
\ee
and drift term
\be
v(t,x)= \alpha( \hat u_a)  + \frac{\alpha(\hat u) - \alpha(e^{-x} \psi_P)}{\hat u - e^{-x} \psi_P} e^{-x} \psi_P.
\ee
We point out that $\psi_E\geq 0$, by the choice of initial
condition and the maximum principle.

By the comparison principle and~\eqref{e.c41914}, we clearly have
\be\label{e.c41915}
	\psi_P(T,x)
		\leq P(x)
		\leq \frac{C}{K\sqrt T}.
\ee
On the other hand,~\eqref{e.c41913} conserves mass, so that
\be
	\int \psi_E(T,x) \dx
		= \int \psi_E(T_g,x)\dx.
\ee
Arguing exactly as in~\cite[Proof of Lemma~5.5]{AHR}, we see that, 
possibly after increasing $K$, we have 
\[
\psi_E(T_g,x) = 0,~~\hbox{for } x \geq - \frac{1}{2}\log(T) + C,
\]
due to~\eqref{e.c41912} and~\eqref{e.c41914}.  Hence, we have 
\be
	\int \psi_E(T_g,x)\dx
		\leq \int_{-\infty}^{-\frac{1}{2}\log T + C}
			p_a(t,x) \dx
		\leq \int_{-\infty}^{-\frac{1}{2}\log T + C}
			e^x \dx
		= \frac{C}{\sqrt T}.
\ee
By parabolic regularity theory, we have
\be\label{e.c41916}
	\psi_E(T,x)
		\leq C\int \psi_E(T,x) \dx
		= C\int \psi_E(T_g,x)\dx
		\leq \frac{C}{\sqrt T}.
\ee
Thus, combining~\eqref{e.c41915} and~\eqref{e.c41916}, we conclude
\be\label{e.c41917}
	p_a(T,x)
		\leq \frac{C}{\sqrt T}
		\qquad\text{ for any } T \geq 4 C e^{\frac{3\gamma a}{2}}.
\ee

\subsubsection*{Step five: conclusion of the proof of \Cref{prop.pp_O(1)}.(i)}

Fix any $T$ sufficiently large and let
\be
	a = \frac{1}{2\gamma \alpha'(1)} \log T.
\ee
Recall that $\gamma$ is a fixed number in $(1,4/3)$.  By~\eqref{e.c41904}, we have
\be
	\sup_x p(T,x)
		\leq \sup_x p_a(T,x) + \frac{C}{\sqrt T}.
\ee
Notice that
\be
	4 C e^{\frac{3\gamma a}{2}}
		= 4 C e^{\frac{3}{4 \alpha'(1)} \log T}
		= 4 C T^\frac{3}{4\alpha'(1)}.
\ee
Since $\alpha(0) = 0$, $\alpha(1) = 1$, and $\alpha$ is convex, we know 
that $\alpha'(1) \geq 1$.  It follows that
\be
	4 C e^{\frac{3\gamma a}{2}}
		\leq T
\ee
for $T$ sufficiently large.  We may thus apply~\eqref{e.c41917} to deduce
\be
	\sup_x p(T,x)
		\leq \sup_x p_a(T,x) + \frac{C}{\sqrt T}
		\leq \frac{C}{\sqrt T},
\ee
which concludes the proof of \Cref{prop.pp_O(1)}.(i), 
except for the proof of Lemmas~\ref{lem.rho} and~\ref{lem.Nash_condition}.~$\Box$



 \subsubsection*{\bf Proof of \Cref{lem.rho}}

We first show how the function $F$ in (\ref{e.F}) 
can be obtained, and why it satisfies~\eqref{e.c41806bis}.  We begin with the ansatz
\be
	\rho(t,x)
		= F(\hat u(t,x)),
\ee
leaving $F$ as yet undetermined.  Intuitively, as discussed in~\cite{AHR}, we seek $\rho$ that gives more weight to the right than to the left; hence, we wish it to be $1$ at $x=+\infty$ and $0$ at $x=-\infty$.  This motivates the boundary conditions
\be\label{e.c41806}
	1 = \rho(t,\infty) = F(0)
		\quad\text{ and }\quad
	0 = \rho(t,-\infty) = F(1).
\ee
Next, we compute: 
\be\label{dec1314}
\bal
	\rho_t+&(\alpha(\hat u)\rho)_x-\rho_{xx}
		=F'(\hat u)\hat u_t+[\alpha'(\hat u)F(\hat u)+\alpha(\hat u)F'(\hat u)]\hat u_x
-F'(\hat u)\hat u_{xx}-F''(\hat u)\hat u_x^2
		\\&
		=F'(\hat u)[\hat u_{xx}+\hat u-A(\hat u)+2\hat u_x-A'(\hat u)\hat u_x]+[\alpha'(\hat u)F(\hat u)+\alpha(\hat u)F'(\hat u)]\hat u_x
		\\&
			\qquad -F'(\hat u)\hat u_{xx}-F''(\hat u)\hat u_x^2
		\\&
		=F'(\hat u)\big[\hat u-A(\hat u)\big]
			+ \hat u_x\big[2F'(\hat u)-A'(\hat u)F'(\hat u)+\alpha'(\hat u)F(\hat u)+\alpha(\hat u)F'(\hat u)]-F''(\hat u)\hat u_x^2
		\\&
		=F'(\hat u)\eta_*(\hat u)+\hat u_x\big[2F'(\hat u)-\hat u\alpha'(\hat u)F'(\hat u)+\alpha'(\hat u)F(\hat u)\big]-F''(\hat u)\hat u_x^2.
\enbal
\ee
In the last step we used that in the pushmi-pullyu case we have 
$\eta_*(u) = u- A(u)$.  The last term above makes it clear that we want $F$ to be concave.  It is also natural to expect $F$ to be monotonic.  In view of~\eqref{e.c41806}, this means $F$ is decreasing.  Hence, we require
\be\label{e.c41807}
	F'(u), F''(u) \leq 0
		\quad\text{ for all } u \in (0,1).
\ee

Then, recalling that
\be\label{e.c72602}
	0 \leq \hat w = - \hat u_x - \eta_*(\hat u),
\ee
we find
\be\label{e.c72603}
-F''(u)u_x^2=(-F''(u))(-u_x)(-u_x)\ge (-F''(u))(-u_x)\eta_*(u). 
\ee 
Using~\eqref{e.c72602} and~\eqref{e.c72603} in (\ref{dec1314}) gives 
\be\label{dec1416}
\bal
\rho_t+(\alpha(\hat u)\rho)_x&-\rho_{xx}
=F'(\hat u)\eta_*(\hat u)+\hat u_x\big[2F'(\hat u)-\hat u\alpha'(\hat u)F'(\hat u)+\alpha'(\hat u)F(\hat u)]-F''(\hat u)\hat u_x^2\\
&\ge F'(\hat u)\eta_*(\hat u)+\hat u_x\big[2F'(\hat u)-\hat u\alpha'(\hat u)F'(\hat u)+\alpha'(\hat u)F(\hat u)\big] 
+F''(\hat u)\hat u_x\eta_*(\hat u)
\\&
\geq F'(\hat u)(-\hat u_x)+\hat u_x\big[2F'(\hat u)-\hat u\alpha'(\hat u)F'(\hat u)+\alpha'(\hat u)F(\hat u)+F''(\hat u)\eta_*(\hat u)\big]
\\&
= \hat u_x\big[F'(\hat u)-\hat u\alpha'(\hat u)F'(\hat u)+\alpha'(\hat u)F(\hat u)+F''(\hat u)\eta_*(\hat u)\big].
\enbal
\ee
Hence, we seek a concave function $F(u)$ such that
\be\label{dec1422}
F'(u)-u\alpha'(u)F'(u)+\alpha'(u)F(u)+F''(u)\eta_*(u)=0.
\ee
We construct such an $F$ and then check that it verifies~\eqref{e.c41806} and~\eqref{e.c41807}.

Notice that, using (\ref{jul702}) once again, we have 
\be
\eta_*'(u)=1-\alpha(u)-u\alpha'(u).
\ee
Thus, \eqref{dec1422} can be written as
\be\label{jul714}
\bal
0&=F''(u)\eta_*(u)+F'(u)-u\alpha'(u)F'(u)+\alpha'(u)F(u)\\
&=(F'(u)\eta_*(u))'-F'(u)(1-\alpha(u)-u\alpha'(u))
+F'(u)[1-u\alpha'(u)]+\alpha'(u)F(u)\\
&=(F'(u)\eta_*(u))'-F'(u)(-\alpha(u))
+ \alpha'(u)F(u)=(F'(u)\eta_*(u)+\alpha(u)F(u))'.
\enbal
\ee
Hence, we may take $F(u)$ such that
\be\label{dec1425}
\farc{F'(u)}{F(u)}=-\frac{\alpha(u)}{\eta_*(u)}.
\ee
The boundary condition $F(0)=1$ in~\eqref{e.c41806} implies that 
\be\label{dec1426}
	F(u)
		=\exp\Big(-\int_0^u\farc{\alpha(u')}{\eta_*(u')}du'\Big)
		=
\exp\Big(-\int_0^u\farc{u'-\eta_*(u')}{u'\eta_*(u')}du'\Big),
\ee
which is exactly (\ref{e.F}). 
Note that $\eta_*(1)=0$ and $\eta_*'(1) = -\alpha'(1) <0$ (recall that $\alpha$ is convex and increasing) so that, for some $u_0 \in (0,1)$,
\be\label{dec1429}
	\int_0^1\farc{u'-\eta_*(u')}{u'\eta_*(u')}du'
		\geq \frac{1}{C} \int_{u_0}^1 \frac{1}{1-u'} du'
		=+\infty.
\ee
Hence, $F(1)=0$.  This completes the proof that $F(u)$ satisfies~\eqref{e.c41806}.

In addition, one can see directly that $F' \leq 0$.  It remains to check the concavity condition in~\eqref{e.c41807}. We see from (\ref{jul714}) that
\be\label{dec1424bis}
\bal
\eta_*(u)&F''(u)=-F'(u)+u\alpha'(u)F'(u)-\alpha'(u)F(u)=F(u)\Big[
\frac{\alpha(u)}{\eta_*(u)}-\frac{\alpha(u)}{\eta_*(u)}u\alpha'(u)-\alpha'(u)\Big]\\
&=\farc{F(u)}{\eta_*(u)}\big[\alpha(u)-u\alpha(u)\alpha'(u)-\alpha'(u)u(1-\alpha(u))\big]
=\farc{F(u)}{\eta_*(u)}\big[\alpha(u)- \alpha'(u)u \big]\le 0,
\enbal
\ee 
where the last inequality holds by the assumed convexity of $\alpha$.  Thus~\eqref{e.c41807} holds.  Finally, we note that the fact that $F \in C_{\rm loc}^2$ follows directly from the computations above and the regularity of $\eta_*$.

It remains to prove (\ref{jul716}). 
First notice that, for any $\eps>0$, there is $u_\eps$ so that
\be
	\frac{1}{(1+\eps)} \frac{u \alpha'(u)}{\alpha'(1)}
			\leq \alpha(u)
			\leq \Big( \frac{1+ \eps}{\alpha'(1)}\Big) u \alpha'(u)
		\quad\text{ for all } u \in [u_\eps, 1].
\ee
Indeed, to see this, notice that equality holds when $\eps = 0$ and $u = 1$.  By making $\eps>0$, we can then obtain the above inequalities for some range of $u$ near $1$ due to the regularity of $\alpha$.

Hence, there is $C_\eps$, depending on $\eps$ and changing line-by-line, so that, for all $u \in [u_\eps, 1]$,
\be\label{e.c41902}
	\begin{split}
		F(u)=\exp\Big(-\int_0^u\farc{\alpha(u')}{u'(1-\alpha(u'))}du'\Big)
			\geq \frac{1}{C_\eps}\exp \Big( - \frac{1+\eps}{\alpha'(1)} \int_{u_\eps}^u \frac{\alpha'(u')}{1 -  \alpha(u')} \, du'\Big)
			= \frac{1}{C_\eps}\big(1 - \alpha(u)\big)^\frac{1+\eps}{\alpha'(1)}.
	\end{split}
\ee
The lower bound is clear for $u < u_\eps$ by choosing $C_\eps = 1/ F(u_\eps)$ and using the monotonicity of $F$.  
A similar argument yields the matching bound
\be\label{e.c41902a}
	F(u)
		\leq C_\eps (1 - \alpha(u))^\frac{1}{(1+\eps) \alpha'(1)}
		\quad \text{ for all } u \in [0, 1].
\ee
This completes the proof.~$\Box$

\subsubsection*{Proof of \Cref{lem.Nash_condition}}

From the definition of  $\mu_a(t)$~\eqref{e.normalization} and the relationship between $\hat u$ and $\tilde u$~\eqref{e.moving_frame}, we have
\be
	\alpha(\hat u_a(t, -\tfrac{1}{2}\log(t+1)-\mu_a(t))) = \frac{1}{2}.
\ee
Then, for $y \leq x < - (1/2)\log(t+1) -\mu_a$, 
we have
\be
	\begin{split}
		\rho_a(t,y)
			&= F(\hat u_a(t,y))
			= \exp\Big( - \int_0^{\hat u_a(t,y)} \frac{\alpha(u')}{\eta_*(u')} \, du'\Big)
			= \rho_a(t,x) \exp\Big( - \int_{\hat u_a(t,x)}^{\hat u_a(t,y)} \frac{\alpha(u')}{\eta_*(u')} \, du'\Big)
			\\&
			= \rho_a(t,x) \exp\Big( \int_y^x \frac{\alpha(\hat u_a(t,z))}{\eta_*(\hat u_a(t,z))} \partial_z \hat u_a(t,z)\, dz\Big)
			\leq \rho_a(t,x) \exp\Big( - \int_y^x \alpha(\hat u_a(t,z))\, dz\Big)
			\\&
			\leq \rho_a(t,x) e^{- \frac{1}{2}(x - y)}.
	\end{split}
\ee
In the first inequality, we used that
\be
	\partial_z \hat u_a
		= -\eta_*(\hat u_a) - \hat w_a
		\leq -\eta_*(\hat u_a),
\ee
and in the second inequality, we used that $\alpha(\hat u_a) \geq 1/2$ due to the definition of $\mu_a$ and the monotonicity of $\hat u_a$.  As an aside, we note that this is the motivation for the definition of $\mu_a$ as it guarantees $\alpha$ is positive on all of $(y,x)$.  This is not guaranteed (e.g., the case $\alpha(u) = 4^n (u-(3/4))_+^n$ for $n \geq 3$).

Thus we have, for all $x \leq-(1/2)\log(t+1) -\mu_a$,
\be\label{e.c41910}
	\overline \rho_a(t,x)
		= \int_{-\infty}^x \rho_a(t,y) \dy
		\leq \rho_a(t,x) \int_{-\infty}^x e^\frac{y-x}{2} \dy
		= 2 \rho_a(t,x),
\ee
and, consequentially,
\be\label{e.c41909}
	\doverline \rho_a(t,x)
		\leq 4 \rho_a(t,x)
		\leq 4 \max\{1, \overline \rho_a^2(t,x)\} \rho_a(t,x).
\ee
This is exactly the desired conclusion when $x \leq -(1/2)\log(t+1)- \mu_a$.

To conclude, we need to obtain the desired bound for $x > -(1/2)\log(t+1)- \mu_a$.  
By the definition of $\mu_a$ and the regularity of $\alpha$, it is easy to see that
\be
	\rho_a(t,x)
		\geq \frac{1}{C}
			\quad\text{ for all } x > -\tfrac{1}{2}\log(t+1)-\mu_a.
\ee
Hence,
\be\label{e.c41911}
	\overline \rho_a(t,x)
		\leq C \overline \rho_a(t,x) \rho_a(t,x)
			\quad\text{ for all } x > -\tfrac{1}{2}\log(t+1) -\mu_a.
\ee
When $x>-(1/2)\log(t+1)-\mu_a$, we combine (\ref{e.c41910}) with  (\ref{e.c41911}) and recall that $\rho_a$ and $\overline \rho_a$ are increasing:
\be
	\begin{split}
		\doverline \rho_a(t,x)
			&= \int_{-\infty}^{-\frac{1}{2}\log(t+1)-\mu_a} \overline \rho_a(t,y) \dy
				+ \int_{-\frac{1}{2}\log(t+1)-\mu_a}^x \overline \rho_a(t,y) \dy
			\\&
			\leq 2\int_{-\infty}^{-\frac{1}{2}\log(t+1)-\mu_a} \rho_a(t,y) \dy
				+ \int_{-\frac{1}{2}\log(t+1)-\mu_a}^x C \overline \rho_a(t,y) \rho_a(t,y) \dy
			\\&
			\leq 2 \overline \rho_a\big(t,-\tfrac{1}{2}\log(t+1)-\mu_a\big)
				+ C \overline \rho(t,x) \int_{-\frac{1}{2}\log(t+1) -\mu_a}^x \rho_a(t,y)\dy
			\\&\leq 4  \rho_a(t,-\mu_a)
				+ C \overline \rho_a^2(t,x)
			\leq 4 \rho_a(t,x)
				+ C^2 \overline \rho_a^2(t,x) \rho_a(t,x)
			\\&
			\leq (4 + C^2) \max\{1, \overline \rho_a^2(t,x)\} \rho_a(t,x).
	\end{split}
\ee
This, in addition to (\ref{e.c41909}), concludes the proof of~\eqref{e.c41905}.~$\Box$

\subsubsection{Proof of \Cref{prop.pp_O(1)}.(ii): the lower bound on $\btilde\urd$}

In this section, we work only with $\urd$. As there is no possibility of confusion,
we simply referring to $\urd$ as $u$ below. 
First, to shorten the proof, we use a preliminary bound 
obtained in Section 4.1 of \cite{Giletti}: for any $\eps>0$, we have
\be\label{e.Giletti}
	\tilde u(t, x)
		\geq \frac{C_\eps}{t^\eps} e^{-x - \frac{x^2}{Ct}}.
\ee
Note that this bound is off the optimal by a factor of $t^\eps$ but is useful as a first step.

The second ingredient is a preliminary upper bound on the shape defect
function: for all~$t \geq 1$ and any $x$, we have
\be\label{e.c41704}
	\tilde w(t,x)
		\leq \begin{cases}
				C \qquad &\text{ if } x < -\log t,\\
				\dfrac{C}{t} \big(x+1 + \log t \big) e^{-x - \frac{(x + \log t)^2}{Ct}} \qquad&\text{ if } x > - \log t.
			\end{cases}
\ee
Indeed, the uniform bound on $w$ over $(-\infty, - \log t)$ follows by parabolic regularity theory.  
To obtain the second bound in~\eqref{e.c41704}, we note that, under the present assumptions,
and since we are in the pushmi-pullyu case, the function $\eta_*(u)=u - A(u)$ is concave and
\[
\eta_*'(u)(2-\eta_*'(u))=1-(1-\eta_*'(u))^2\le 1. 
\]
Therefore, it follows from~\eqref{dec2120} that the (non-negative) shape defect function $w$ satisfies
\be\label{e.c41703}
	w_t - w_{xx}
		\leq w.
\ee
A lesson 
of \cite{HNRR} is that bounded functions starting from compactly supported (on the right) initial data that satisfy~\eqref{e.c41703} must also satisfy
\be
	w\big(t, x + 2t - \frac{3}{2}\log t\big)
		\leq C (x + 1) e^{-x - \frac{x^2}{Ct}},
\ee
from which the second bound in~\eqref{e.c41704} follows after changing variables.


The main estimate we need for the proof of \Cref{prop.pp_O(1)}.(ii)
is a uniform bound on the exponential moment of $\tilde u$
\be
I(t) := \int e^x \tilde u(t,x) dx.
\ee
This is provided by the following lemma:
\begin{lem}\label{lem.exponential_moment}
There exists $C>0$ such that, for sufficiently large $t$,
\be\label{jul802}
\frac{1}{C}\leq \frac{I(t)}{\sqrt t}\leq C.
\ee
\end{lem}
We note that the upper bound in \Cref{lem.exponential_moment} is not used in this paper; however, 
it comes `for free' in the proof, so we state it as well.  We postpone the proof 
of \Cref{lem.exponential_moment} momentarily and, first, show how to apply it to conclude 
the lower bound in \Cref{prop.pp_O(1)}.

\subsubsection*{\bf Proof \Cref{prop.pp_O(1)}.(ii)}  

Fix $N>0$ to be chosen.  We argue by contradiction, assuming, for $t$ fixed and sufficiently large,
that
\be
- \mu(t) < - N.
\ee
In this case~\eqref{e.u_asymptotics} implies that
\be\label{e.c7}
	\tilde u(t,x)
		\leq C e^{-x - N}.
\ee
Our goal is to deduce from (\ref{e.c7}) an upper bound on $I(t)$ that violates the lower bound of \Cref{lem.exponential_moment} if $N$ is too large.

Define the left, middle, and right domains as
\be
	L = \{x : x < - N\},
	\quad M = \{x : -N \leq x < N \sqrt t\},
	\quad\text{and}\quad
	R = \{x: N \sqrt t \leq x\},
\ee
with $I_L(t)$, $I_M(t)$, and $I_R(t)$ the decomposition of $I(t)$ into integrals on each respective domain.  
To bound $I_L$, we simply use that $\tilde u\leq 1$:
\be\label{e.c72801}
I_L(t)= \int_{-\infty}^{-N} e^x \tilde u(t,x) \dx
\leq e^{-N}.
\ee
For $I_M$, we use~\eqref{e.c7} to find
\be\label{e.c72802}
	I_M(t)
		= \int_{-N}^{N\sqrt t} e^{x} \tilde u(t,x) \dx
		\leq \int_{-N}^{N\sqrt t} C e^{-N} \dx
		= CN e^{-N} (1 + \sqrt t).
\ee
Finally, we estimate $I_R$ from above.  
We integrate by parts and use the shape defect function to 
obtain
\be
\begin{split}
I_R(t)&=\int_{N\sqrt t}^\infty e^x \tilde u \dx= - e^{N \sqrt t} \tilde u(t,N \sqrt t)
			- \int_{N\sqrt t}^\infty e^x \tilde u_x \dx
		\leq - \int_{N\sqrt t}^\infty e^x \tilde u_x \dx
		\\&
		= \int_{N \sqrt t}^\infty e^x(\tilde w - \eta_*(\tilde u)) \dx
		\leq \int_{N \sqrt t}^\infty e^x \tilde w \dx.
\end{split}
\ee
Applying the second bound in~\eqref{e.c41704} yields
\be\label{e.c72803}
I_R(t)\leq C\int_{N \sqrt t}^\infty  \frac{(x+\log t) }{ t} e^{-\frac{x^2}{Ct}} \dx
		\leq C  e^{-\frac{N^2}{C}}.
\ee
Putting together~\eqref{e.c72801},~\eqref{e.c72802}, and~\eqref{e.c72803} with~\eqref{e.c6}, we find
\be
	\frac{1}{C}
		\leq \frac{I(t)}{\sqrt t}
		\leq \frac{e^{-N}}{\sqrt t} + CN e^{-N} + \farc{C}{\sqrt{t}} e^{-\frac{N^2}{C}}.
\ee
This yields a contradiction for $t$ and $N$ sufficiently large. 
We conclude that there exists $\underbar L$ such that
\be
	\alpha(u(t, 2t - (1/2)\log(t) - \underbar L) )\geq 1/2,
\ee
which finishes the proof.~$\Box$


\subsubsection*{\bf Proof of \Cref{lem.exponential_moment}}

First, we compute the time derivative of $I(t)$, using integration by parts repeatedly 
and the definition of the shape defect function: 
\be\label{e.c4}\begin{split}
\dot I(t)&= \int e^x \Big(\tilde u_{xx} + f(\tilde u) + 
\Big(2 - \frac{1}{2(t+1)}\Big) \tilde u_x \Big) \dx
			\\&
			= \int e^x \Big(-\tilde u_x + \eta_*(\tilde u)(2 - \eta_*'(\tilde u)) + 
2 \tilde u_x + \frac{1}{2(t+1)} \tilde u \Big) \dx
			\\&
			=  \int e^x (\eta_*(\tilde u) + (\eta_*(\tilde u))_x -\tilde w(1 - \eta_*'(\tilde u))) \dx
				+ \frac{I}{2(t+1)}
			\\&= - \int e^x (1 - \eta_*'(\tilde u)) \tilde w \dx
				+ \frac{I}{2(t+1)}.
	\end{split}
\ee
In the pushmi-pullyu case we have
\[
	u - \eta_*(u) = A(u)
	\quad\text{ so that }\quad
	1- \eta_*'(u) = A'(u).
\]
Recall that $A(u)$ is an increasing $C^2$ function with $A(0) = A'(0) = 0$. It follows that 
\be\label{e.c41803}
	0 \leq \int e^x (1 - \eta_*'(\tilde u)) \tilde w \dx
		\leq C\int e^x \tilde u \tilde w \dx
		=: \cE(t).
\ee

The upper bound of $I$ is immediate from~\eqref{e.c4} and the fact that $\cE (t)\geq 0$
due to the positivity of~$w$.  Next, we prove the lower bound of $I$.
We claim that
\be\label{e.c42605}
	\cE (t)\leq \frac{C\log^2 (t+1)}{t+1}.
\ee

Before establishing~\eqref{e.c42605}, we show how to use it to conclude the lower bound of $I$.  
Using~\eqref{e.c42605} in~\eqref{e.c4} yields
\be
	\frac{d}{dt}\left(\frac{I(t)}{\sqrt{t+1}}\right)
		\geq -\frac{C \log^2(t+1)}{(t+1)^{3/2}}.
\ee
Fix $t_0>0$ to be chosen.  Integrating the above from $t_0$ to $t$, we find
\be\label{e.c5}
	\frac{I(t)}{\sqrt {t+1}}
		\geq \frac{I(t_0)}{\sqrt{t_0+1}} - C\int_{t_0}^t\frac{\log^2(s+1)}{(s+1)^{3/2}}\ds
		\geq \frac{I(t_0)}{\sqrt{t_0+1}} - C\frac{\log^2(t_0+1)}{\sqrt{t_0+1}}.
\ee
On the other hand, for any $\eps>0$, by~\eqref{e.Giletti}, we have
\be
	\frac{I(t_0)}{\sqrt{t_0+1}}
		\geq \frac{C}{(t_0+1)^{\eps}}.
\ee
Hence, fixing $t_0$ sufficiently large,  the right  side of~\eqref{e.c5} is positive.  
It follows that
\be\label{e.c6}
	\frac{I(t)}{\sqrt{t+1}}
		\geq C_0
			\qquad\text{ for all $t> t_0$,}
\ee
which yields the claimed lower bound in (\ref{jul802}).

To finish the proof, we establish~\eqref{e.c42605}.  We decompose $\cE(t)$ as
\be
	\cE(t) = \int_{-\infty}^{-\log (t+1)} e^x \tilde u \tilde w \dx
		+ \int_{-\log(t+1)}^0 e^x \tilde u \tilde w \dx
		+ \int_0^\infty e^x \tilde u \tilde w \dx
		= \cE_L(t) + \cE_M(t) + \cE_R(t).
\ee
For $\cE_L(t)$, we note that $0 \leq \tilde u, \tilde w \leq C$ in order to find
\be\label{e.c42606}
	\cE_L
		\leq C \int_{-\infty}^{-\log (t+1)} e^x \dx
		= C e^{-\log (t+1)}
		\leq \frac{C}{t+1}.
\ee
Next, on the domain of integration of $\cE_M(t)$,~\eqref{e.c41704} implies that 
\[
e^x\tilde w\le \frac{C\log(t+1)}{t+1}.
\]
This gives the bound
\be\label{e.c42607}
\cE_M(t)\leq 
\frac{C \log(t+1)^2}{t+1}.
\ee
Finally, using~\eqref{e.c41704} again, as well as the upper bound 
for $\tilde u$ in~\eqref{e.u_asymptotics}, yields
\be
\cE_R(t)
		\leq C \int_0^\infty \frac{x + \log(t+1)}{t+1} e^{-x - \mu(t) - \frac{(x+\log(t+1))^2}{C(t+1)}} \dx
		\leq \frac{C\log(t+1)e^{-\mu(t)}}{t+1}.
\ee
By the already proved \Cref{prop.pp_O(1)}.(i), we know that~$-\mu(t) \leq \overline L$, and, thus,
\be\label{e.c42608}
	\cE_R(t)
		\leq \frac{C\log(t+1)}{t+1}.
\ee
Putting together~\eqref{e.c42606},~\eqref{e.c42607}, and~\eqref{e.c42608}, 
we obtain~\eqref{e.c42605}, which completes the proof of \Cref{lem.exponential_moment}.~$\Box$

\subsection{The second step: the behavior of $\tilde u$ at $x=t^\gamma$}
\label{sec:pp_second}

We now identify the constant order term in the expansion of $m(t)$.  
This is done in the manner of~\cite{NRR1}.  The main step is  the following proposition.
\begin{prop}\label{prop.pp_o(1)}
	Let $u$ be the solution to either~\eqref{e.rde} or~\eqref{e.rcl} under the assumptions of \Cref{t.rde-intro} or \Cref{t.rcl-intro} with $\tilde u$ defined by~\eqref{e.moving_frame}.  Then there exists $\alpha_\infty>0$ such that, for any~$\gamma \in (0,1/2)$, we have
	\be
		\lim_{t\to\infty}
			e^{t^\gamma}\tilde u(t, t^\gamma)
			= \alpha_\infty.
	\ee
\end{prop}

As the proof of this is similar to that of \cite[Corollary~6.3]{AHR}, which is, in turn, 
based on the proof of \cite[Lemma~4.2]{NRR1}, we merely provide an outline of the main points, 
as well as a description of the changes needed to be made in the present setting.

As the proof is significantly simpler for solutions to the reaction-diffusion equation~\eqref{e.rde} 
than for solutions to the reactive conservation law~\eqref{e.rcl},
we begin with the proof in former case.

\subsubsection{The proof of \Cref{prop.pp_o(1)} for the reaction-diffusion equation (\ref{e.rde})}\label{s.pp_o(1)}

In this subsection we only use $\urd$, so we drop the ``rd'' subscript.   
We begin with the weighted Hopf-Cole transform
\be
	v(t,x)
		= e^\Gamma \tilde u(t,x)
		= \exp\Big\{x + \int_x^\infty \alpha(\tilde u(t,y)) dy\Big\} \tilde u(t,x),
\ee
which satisfies
\be\label{e.c42107}
	v_t + \frac{1}{2(t+1)} (v_x - v) - v_{xx} = G,
\ee
with 
\be\begin{split}
G
&= e^\Gamma \Big(f(\tilde u) - \tilde u - \alpha(\tilde u)^2 \tilde u + 2 \alpha(\tilde u) \tilde u_x\Big) +  u \Big( - \int_x^\infty \alpha''(\tilde u) \tilde u_y^2 dy + \int_x^\infty \alpha'(u)  f(\tilde u) \dy\Big).
	\end{split}
\ee
See the details of \Cref{prop-mar902} in order to see how~\eqref{e.c42107} is computed.  Due to \Cref{prop-mar902}, we have
\be\label{e.c42101}
	G \leq 0.
\ee
Furthermore, for all $\gamma >0$, $t\geq 1$, and $x \geq 0$, we have 
\be\label{e.c42102}
	G(t, x + t^\gamma)
		\geq - C \exp(-x - t^\gamma).
\ee
Here, we are using that, for $u$ small, $f(u) - u = O(u^2)$, and, for all $t\geq 1$ and $x\geq 0$,
we have 
\be\label{e.c42103}
	\tilde u(t,x) \leq C e^{-x}
	\quad\text{ and }\quad
	-\tilde u_x \leq C e^{-x}.
\ee
The former is due to the upper bounds in \Cref{prop.pp_O(1)} and~\eqref{e.u_asymptotics}, 
while the latter follows from the former by parabolic regularity theory.

Additionally, the argument of \cite[Lemma~6.1]{AHR} applies nearly verbatim, so we assert its conclusion without proof: for all $t\geq 1$,
\be\label{e.c42104}
	v(t,x) \leq C
		\quad\text{ for all } x \geq 0,
	\quad\text{ and }\quad
	v(t,x) \geq \frac{1}{C}
		\quad\text{ for all } x \leq \frac{\sqrt t}{C}.
\ee
Finally, we claim that, for all $t\geq 1$, we have
\be\label{e.c42105}
	v_x(t, t^\gamma)
		\leq 0
	\quad\text{ and }\quad
	|v_x(t, - t^\gamma)|
		\leq Ce^{ - \alpha'(1) t^\gamma}.
\ee
We see this as follows: using (\ref{jul702}), we write
\be\label{e.c50601}
	v_x(t,x) = (\tilde u_x(t,x) + \eta_*(\tilde u(t,x)))e^\Gamma=-w e^\Gamma\le 0,
\ee
giving the first inequality in~\eqref{e.c42105}.  

We now justify the second inequality in~\eqref{e.c42105}.  
Using the computation~\eqref{e.c50601}, there are three terms to bound: $\tilde u_x$, $\eta_*(\tilde u)$, and $\exp(\Gamma).$ 
First, observe that, by \Cref{prop.pp_O(1)} and \eqref{e.u_asymptotics}, we have 
\be\label{e.c42108}
	0 \leq 1- \tilde u(t,x)
		\leq C e^{\alpha'(1) x}.
\ee
Thus, by parabolic regularity theory, we have 
\be\label{e.c50602}
	|\tilde u_x(t,-t^\gamma)|
		\leq C e^{-\alpha'(1) t^\gamma}.
\ee
Next, from~\eqref{e.c42108} and a Taylor approximation (recall that $\alpha(1) = 1$), we find
\be
	0
		\leq 1 - \alpha(\tilde u(t,x))
		\leq C e^{\alpha'(1) x}.
\ee
This yields  
\be\label{e.c42109}
	\eta_*(\tilde u(t,-t^\gamma))
		= \tilde u(t,-t^\gamma) (1 - \alpha(\tilde u(t,-t^\gamma)))
		\leq (1 - \alpha(\tilde u(t,-t^\gamma)))
		\leq C e^{-\alpha'(1) t^\gamma}.
\ee
In addition, using~\eqref{e.u_asymptotics} and a Taylor approximation again, we find, for $x \leq 0$,
\be
	\begin{split}
		\exp\Big| x + &\int_x^\infty\alpha(\tilde u(t,y)) \dy\Big|
			= \exp\Big|\int_x^0(\alpha(\tilde u(t,y))-1)dy + \int_0^\infty \alpha(\tilde u(t,y)) dy\Big|
		\\&
			\leq \exp\Big|\int_x^0 C e^{\alpha'(1) y}dy\Big| \cdot \exp \Big|\int_0^\infty C e^{-y} dy\Big|
			\leq C.
	\end{split}
\ee
The combination of~\eqref{e.c50601},~\eqref{e.c50602} and~\eqref{e.c42109} yields~\eqref{e.c42105}. 

As discussed in detail in \cite[Section 6]{AHR},  these ingredients, that is,~\eqref{e.c42101},~\eqref{e.c42102},~\eqref{e.c42104}, and~\eqref{e.c42105}, are all that 
is needed to prove \Cref{prop.pp_o(1)}.  This concludes the proof of \Cref{prop.pp_o(1)} in the case of the reaction-diffusion equation~\eqref{e.rde}.

\subsubsection{The proof of \Cref{prop.pp_o(1)} for the reactive conservation law (\ref{e.rcl})}\label{sec:rcl-hopf}

We now drop the ``rcl'' subscript from $\urcl$ and denote by $u$ 
the solution to~\eqref{e.rcl}.
We begin by reviewing the key ingredients in \Cref{s.pp_o(1)} for the reaction-diffusion case.  
Clearly it is possible to establish~\eqref{e.c42104} and ~\eqref{e.c42105} verbatim in the reactive conservation law case. 
In order to check~\eqref{e.c42101} and~\eqref{e.c42102}, we need to 
consider the weighted Hopf-Cole transform 
\be
v(t,x)= e^\Gamma \tilde u(t,x)
= e^{x + \int_x^\infty \alpha(\tilde u(t,y)) dy} \tilde u(t,x).
\ee
Changing to the moving frame from~\eqref{e.rcl}, we have
\be\label{e.c42110}
	\tilde u_t - \Big(2 - \frac{1}{2(t+1)}\Big)\tilde u_x
		+ (\tilde u \alpha(\tilde u))_x
		= \tilde u_{xx} + \tilde u(1 - \alpha(\tilde u)).
\ee
This allows us to compute
\be
	\begin{split}
		v_t + &\frac{1}{2(t+1)} (v_x - v) - v_{xx}
			= \Gamma_t v
				+ e^\Gamma \tilde u_t
				+ \frac{1}{2(t+1)} (e^\Gamma \tilde u_x - \alpha v)
				- ((1 - \alpha) v + e^\Gamma \tilde u_x)_x
			\\&
			= \Gamma_t v
				+ e^\Gamma
					\Big(
					\Big(2 - \frac{1}{2(t+1)} \Big)\tilde u_x
					- (\tilde u \alpha(\tilde u))_x
					+ \tilde u_{xx}
					+ \tilde u(1 - \alpha(\tilde u))
					\Big)
				+ \frac{1}{2(t+1)} (e^\Gamma \tilde u_x - \alpha v)
			\\&\qquad
				- \left(
					- \alpha'(\tilde u) \tilde u_x v
					+ (1-\alpha)^2 v
					+ 2 (1-\alpha) e^\Gamma \tilde u_x
					+e^\Gamma \tilde u_{xx}
				\right)
			\\&
			= \Gamma_t v
				- \frac{\alpha}{2(t+1)} v
				- \alpha(\tilde u) e^\Gamma \tilde w.
	\end{split}
\ee
To compute $\Gamma_t$, we use again~\eqref{e.c42110}, as well as (\ref{jul702}) 
and  integration by parts, to find
\be
	\begin{split}
	\Gamma_t
		&= \int_x^\infty \alpha'(\tilde u)
			\Big(
					\Big(2 - \frac{1}{2(t+1)} \Big)\tilde u_x
					- (\tilde u \alpha(\tilde u))_x
					+ \tilde u_{xx}
					+ \tilde u(1 - \alpha(\tilde u))
					\Big)
 \dy
		\\&
		= \frac{\alpha}{2(t+1)}
			+ \tilde u \alpha \alpha'
			- \tilde u_x \alpha'(\tilde u)
			+ \int_x^\infty \big(2 \alpha' \tilde u_x
					+ \alpha''(\tilde u) \tilde u \alpha(\tilde u) \tilde u_x
					- \alpha''(\tilde u) \tilde u_x^2
					+ \alpha'(\tilde u) \eta_*(\tilde u)
				\big) \dy
		\\&
		= \frac{\alpha}{2(t+1)}
			+ \alpha' \tilde w
			+ \int_x^\infty \left(
					- \tilde u \alpha'' \tilde u_x
					+ \alpha' \tilde u_x
					+ \alpha''(\tilde u) \tilde u \alpha(\tilde u) \tilde u_x
					- \alpha''(\tilde u) \tilde u_x^2
					+ \alpha'(\tilde u) \eta_*(\tilde u)
				\right) \dy
		\\&
		= \frac{\alpha}{2(t+1)}
			+ \alpha' \tilde w
			- \int_x^\infty (\alpha' - \alpha'' \tilde u_x) \tilde w \dy.
	\end{split}
\ee
Combining both computations above yields
\be\label{e.c42111}
	v_t
		+ \frac{1}{2(t+1)} (v_x - v) - v_{xx}
		= e^\Gamma \Big( (\tilde u\alpha' - \alpha) \tilde w
				-\tilde u \int_x^\infty (\alpha' - \alpha'' \tilde u_x) \tilde w \dy
				\Big).
\ee
Unfortunately, due to the convexity of $\alpha(u)$, the first term in the right side 
of~\eqref{e.c42111} is positive.  Amazingly, it is zero for exactly one choice, $\alpha(u)=u$, 
which was by coincidence considered in \cite{AHR}.  As a result, the analogue of~\eqref{e.c42101} does not hold in the reactive conservation law case.

Let us now outline how to bypass this difficulty. 
Arguing exactly as in \Cref{s.pp_o(1)}, it is easy to check that the 
analogue of~\eqref{e.c42102} holds; that is,
\be
	v_t
		+ \frac{1}{2(t+1)} (v_x - v) - v_{xx}
			\geq - C e^{ -x}
				\qquad\text{ for } x \geq t^\gamma.
\ee
In addition, as we have observed, $v(t,x)$ still satisfies the first inequality in
(\ref{e.c42105}): $v_x(t,t^\gamma)\le 0$. Therefore, a subsolution for $v(t,x)$
can be found as the solution to 
\be
\underline v_t
		+ \frac{1}{2(t+1)} (\underline v_x - \underline v) - \underline v_{xx}
			= - C e^{-x}
				\quad\text{ for } x \geq t^\gamma,
\ee
with the Neumann boundary condition $\underline v_x(t,t^\gamma)=0$.  That is, if 
$v(T,x)=\underline v(T,x)$ for all $x\ge T^\gamma$ at some time $T$, then 
\be
v(t,x)\ge \underline v(t,x),~~~\hbox{ for all $t\ge T$ and $x\ge t^\gamma$.}
\ee
This part of the proof is unaffected and proceeds exactly as in~\cite{AHR}.

In order to construct a supersolution for $v(t,x)$, we note that 
\be
	e^\Gamma \Big( (\tilde u\alpha' - \alpha) \tilde w
				-\tilde u \int_x^\infty (\alpha' - \alpha'' \tilde u_x) \tilde w \dy
				\Big)
		\leq C v \tilde w.
\ee
Therefore, a supersolution for $v(t,x)$ is given by the solution to 
\be\label{e.c42113}
	\overline v_t
		+ \frac{1}{2(t+1)} (\overline v_x - \overline v) - \overline v_{xx}
			= C \overline v \tilde w
				\quad\text{ for } x \geq -t^\gamma,
\ee
with the boundary condition
\be\label{e.c42202}
	\overline v_x(t,- t^\gamma)
		= - C e^{-\alpha'(1) t^\gamma}.
\ee
Here, we took into account the second inequality in
(\ref{e.c42105}).
%
%

Changing to self-similar variables
\be
	\overline V(\tau, y)
		= \overline v( e^\tau, e^{\tau/2} y),
\ee
yields  
\be\label{e.c42114}
	\overline V_\tau
		+ \cL \overline V
		+ \frac{1}{2} e^{-\tau/2} \overline V_y
		= C \overline V e^\tau \tilde w(e^\tau, e^{\tau/2}y)
		\quad\text{ for } y \geq - e^{-\tau(1/2 - \gamma)}.
\ee
Here, $\cL$ is the linear operator associated to the heat equation in self-similar variables:
\be
	\cL
		:= - \partial_y^2 - \frac{y}{2} \partial_y - \frac{1}{2}.
\ee
A key feature required to make the proof strategy of~\cite[Corollary~6.3]{AHR} work is that, 
the right hand side of~\eqref{e.c42114} must be `suitably' bounded and tend to 
zero in a `suitable' fashion as $\tau\to+\infty$.  In particular, one needs to show that
\be\label{e.c42201}
	\int_0^\infty \|e^\tau \tilde w(e^\tau, e^{\tau/2}\cdot)\|_{L^2_{\smuG}} \, d\tau
		<\infty.
	\quad\text{ and } \quad
	\sup_{\tau\geq 0} \|e^\tau \tilde w(e^\tau,\cdot)\|_{L^\infty} < \infty
\ee
with the Gaussian-weighted measure 
\[
d\tmuG = e^{y^2/4} \dy.
\]
Roughly, this allows to treat the right side of~\eqref{e.c42114}
as an error term when pursuing $L^2_{\smuG}$ estimates.

The bounds in~\eqref{e.c42201} are a consequence of the following lemma, that we state here and prove in \Cref{s.w_bounds}.
\begin{lem}\label{l.w_bounds}
For all $t\geq 1$, we have 
	\be\label{jul1004}
		\tilde w(t,x)
			\leq 
				\begin{cases}
					 {C}{t}^{-1}
						\qquad&\text{ if } x < 1,\\
					 {C}{t}^{-1} x e^{-x - \frac{x^2}{5t}}
						\qquad&\text{ if } x \geq 1.
				\end{cases}
	\ee
\end{lem}

We now discuss how the above subsolution and supersolution can be used to prove 
the the convergence of $V(\tau, y) = v(e^\tau, e^{\tau/2}y)$ in a slightly more detail.  
First, for the sake of a simple discussion, we assume that both $\overline V$ and $\underline V$ (defined analogously using $\underline v$) satisfy Neumann boundary conditions -- it is 
easy to account for the the errors coming 
from the `approximate' Neumann boundary condition~\eqref{e.c42202} as they are exponentially small in $t^\gamma$.

Fix $T\gg1$ and impose initial conditions
\be
	\underline V(T,\cdot) = \overline V(T,\cdot) = V(T,\cdot).
\ee
Notice that, changing to self-similar variables, the equation for $\underline V$ becomes
\be\label{e.c42203}
	\underline V_\tau
		+ \cL \underline V
		+ \frac{1}{2} e^{-\tau/2} \underline V_y
		= - C e^{ - ye^{\tau/2}}
		\quad\text{ for }
		y \geq e^{-\tau (1/2 - \gamma)}.
\ee
Recall the equation for $\overline V$ is given by~\eqref{e.c42114}.   Let us take for granted that $\underline V(t,\cdot)$, and $\overline V(t,\cdot)$ are bounded in $L^2_{\smuG} \cap L^\infty$ uniformly in $t$ and $T$ (this is easily established for $V$ in exactly the manner of~\cite{AHR, HNRR, NRR1} and it is then inherited by $\underline V$ and $\overline V$).

We consider the projection onto the principal eigenfunction 
\be
	\psi(y)
		= \frac{1}{Z} e^{-y^2/4}.
\ee
Here, $Z$ is chosen so that $\| \psi\|_{L^2_{\smuG}} = 1$.  For $\overline V$, we find
\be
	\frac{d}{dt} \langle \overline V, \psi\rangle_{L^2_{\smuG}}
		= -\frac{1}{2} e^{-\tau/2}\overline V(\tau, 0)
			+ Ce^\tau \|\tilde w(e^\tau, e^{\tau/2}\cdot)\|_{L^2_{\smuG}},
\ee
and, for $\underline V$, we find
\be
	\frac{d}{dt} \langle \underline V, \psi\rangle_{L^2_{\smuG}}
		= -\frac{1}{2} e^{-\tau/2} \underline V(\tau, 0)
			- Ce^{-\tau/2 - e^{\gamma \tau}}.
\ee
Integrating both identities above and using the initial data, it follows that, for any $\tau > T$,
we have 
\be
	\begin{split}
	- C e^{-T/2}
		&\leq  \langle \underline V(\tau,\cdot), \psi\rangle_{L^2_{\smuG}}
			-  \langle V(T,\cdot), \psi\rangle_{L^2_{\smuG}}
		\leq \langle  V(\tau,\cdot), \psi\rangle_{L^2_{\smuG}}
			-  \langle V(T,\cdot), \psi\rangle_{L^2_{\smuG}}
		\\&\leq \langle \overline V(\tau,\cdot), \psi\rangle_{L^2_{\smuG}}
			-  \langle V(T,\cdot), \psi\rangle_{L^2_{\smuG}}
		\leq C\int_T^\tau e^{\tau}\|\tilde w(\tau',\cdot)\|_{L^2_{\smuG}} \, d\tau'.
	\end{split}
\ee
As the left and right sides tend to zero as $T\to\infty$,  due to~\eqref{e.c42201}, 
we deduce that $\langle V(\tau,\cdot), \psi\rangle_{L^2_{\smuG}}$ is a Cauchy sequence in $\tau$.  
It follows that there is $\alpha_\infty$ so that
\be
	\langle V, \psi\rangle_{L^2_{\smuG}}
		\to \alpha_\infty
		\qquad\text{ as } \tau \to \infty.
\ee
By~\eqref{e.c42104}, we also know that $\alpha_\infty > 0$.

Similar arguments using the spectral gap of $\cL$ 
show that
\be
	\|V - \langle V, \psi\rangle_{L^2_{\smuG}} \psi\|_{L^2_{\smuG}}
		\to 0
		\quad\text{ as } \tau\to\infty,
\ee
which implies that
\be
	\|V - \alpha_\infty \psi\|_{L^2_{\smuG}}
		\to 0
		\quad\text{ as } \tau\to\infty.
\ee
Then, using parabolic regularity theory and the boundedness of the right hand sides of~\eqref{e.c42203} and~\eqref{e.c42114} (here we are using the $L^\infty$ bound of $e^\tau\tilde w$~\eqref{e.c42201}), this can be upgraded to
\be
	\|V - \alpha_\infty \psi\|_{L^\infty_{\rm loc}}
		\to 0
		\quad\text{ as } \tau\to\infty.
\ee
After undoing the change of variables, this is precisely \Cref{prop.pp_o(1)}.  As the technical details are exactly those in~\cite{AHR, NRR1}, except for the small changes indicated above, we omit them.

\subsection{The third step: convergence to the wave}
\label{sec:pp_third}

We now conclude the proof of Theorems~\ref{t.rde-intro}.(ii) and \ref{t.rcl-intro}.  We note that the proof 
of the final step here is different from that in~\cite{AHR,NRR1,NRR2}.  
It is no extra effort to consider both cases at the same time, so
we use $\tilde u$ for both $\tilde\urd$ and $\tilde\urcl$.

Fix $\gamma \in (0,1/10)$.  We let $\alpha_\infty>0$ be the constant from \Cref{prop.pp_o(1)}, 
and choose $x_\infty$ be such that
\be\label{e.c42505}
\lim_{t\to\infty} e^{t^\gamma} U_\infty(x_\infty + t^\gamma) = \alpha_\infty.
\ee
Consider the difference 
\be
	s(t,x) = e^x (\tilde u(t,x) - U_*(x_\infty + \eps + x)).
\ee
We point out that, due to~\eqref{e.c42505} and \Cref{prop.pp_o(1)}, we have 
\be\label{e.c42510}
	0 < s(t,t^\gamma)
		< 2 \eps,
\ee
for $t$ sufficiently large. 
Using the identities
\be
	\tilde w
		= - \tilde u_x - \tilde u + A(\tilde u)
	\quad\text{ and }\quad
	0 = - U_*' - U_* + A(U_*),
\ee
we derive a first order ODE for $s$:
\be\label{e.c42508}
	s_x
		= e^x( \tilde u_x + \tilde u - U_*' - U_*)
		= e^x( - \tilde w + A(\tilde u) - A( U_*))
		= - e^x \tilde w + \xi s
\ee
where
\be
	\xi = \frac{A(\tilde u) - A(\tilde U_*)}{\tilde u - U_*}
		\geq 0.
\ee
As $\tilde w \geq 0$ and $s(t,t^\gamma)>0$,~\eqref{e.c42508} implies that
\be\label{e.c42511}
	s(t,x) > 0
		\qquad\text{ for all } x < t^\gamma.
\ee

Next, by~\eqref{e.c41704}, in the reaction diffusion case, and \Cref{l.w_bounds}, in the reactive conservation law case, we have
\be
	\tilde w(t,x)
		\leq \frac{C (x+\log t+1)}{t} e^{-x - \frac{(x+\log t)^2}{Ct}}
		\qquad \text{ for all } x > -\log t,
\ee
so that
\be
	s_x > - \frac{C(x + \log t+1)}{t} e^{- \frac{(x+\log t)^2}{Ct}},
	~~~\hbox{for all $x \in (-\log t, t^\gamma)$}. 
\ee
Integrating this, we find, for any $x > - \log(t)$, we have
\be
	s(t,x)
		- s(t,t^\gamma)
		< C \int_x^{t^\gamma}\frac{(y + \log t+1)}{t} e^{- \frac{(y + \log t)^2}{Ct}} \dy
		\leq \frac{C}{\sqrt t}.
\ee
In view of~\eqref{e.c42510}, we deduce that 
\be\label{e.c42512}
	s(t,x)
		< \frac{C}{\sqrt t} + 2 \eps.
\ee
Putting together~\eqref{e.c42511} and~\eqref{e.c42512} and undoing the change of variables, we obtain
\be
	0
		< \tilde u(t,x)
			- U_*(x_\infty + \eps + x)
		< e^{-x} \Big( \frac{C}{\sqrt t} + 2\eps\Big)
		\qquad\text{ for all } x  < t^\gamma.
\ee
Notice also that
\be
	\|U_*(x_\infty + \eps + \cdot) - U_*(x_\infty + \cdot)\|_{L^\infty}
		\leq \|U_*'\|_{L^\infty} \eps
		\leq C \eps.
\ee
As a consequence of the previous two inequalities, we obtain 
\be\label{e.c42601}
	\sup_{x \in [-\log\frac{1}{2\eps}, t^\gamma]}
		|\tilde u(t,x) - U_*(x_\infty + x)|
		\leq \frac{C}{\sqrt t} + C \sqrt \eps.
\ee
On the other hand, clearly we have 
\be\label{e.c42602bis}
	\begin{split}
	\sup_{x < - \log\frac{1}{2\eps}} |\tilde u(t,x) - U_*(x_\infty +x)|
		&\leq \sup_{x < - \log\frac{1}{2\eps}} |\tilde u(t,x) - 1|
			+ \sup_{x < - \log\frac{1}{2\eps}} |1 - U_*(x_\infty +x)|
		\\&
		\leq C e^{-\alpha'(1) \log\frac{1}{2\eps}}
		= C \eps^{\alpha'(1)/2}
	\end{split}
\ee
and
\be\label{e.c42603}
	\sup_{x > t^\gamma} |\tilde u(t,x) - U_*(x+ x_\infty)|
		\leq \sup_{x > t^\gamma} \tilde u(t,x)
			 + \sup_{x > t^\gamma} U_*(x+x_\infty)
		\leq C e^{-t^\gamma}.
\ee
These inequalities can be seen by~\eqref{e.u_asymptotics}, \Cref{prop.pp_O(1)} and~\eqref{e.U_asymptotics}.

The proof is now finished by combining~\eqref{e.c42601},~\eqref{e.c42602bis}, and~\eqref{e.c42603}.
~$\Box$

\subsection{The proof of Lemma~\ref{l.w_bounds}}\label{s.w_bounds}

For notational ease, we drop the `rcl' subscript for the remainder of the proof as we only work with the reactive conservation law in this section.

As the construction of the supersolution is somewhat complicated and opaque, let us first describe, roughly, why the claim holds.  
First, fix $T\gg 1$ and recall that $\tilde w$ satisfies
\be
	\tilde w_t
		- \Big(2 - \frac{1}{2(t+T)}\Big) \tilde w_x
		+\, A'(\tilde u) \tilde w_{x}
			\,\,=\,\, \tilde w_{xx}
				+ \tilde w(1- A'(\tilde u)).
\ee
Ignoring the $A'$ terms, momentarily, the work of \cite[equation (20)]{HNRR}, 
suggests that, for $x\geq 1$, we should have 
\be\label{jul1002}
	\tilde w(t,x)
		\leq \frac{C x}{t} e^{-x - \frac{x^2}{(4 + o(1)) t}}.
\ee
This is exactly the desire bound on $x \geq 1$.  
On the other hand, in order to obtain a bound on $x<1$, we must understand the reaction coefficient $1 - A'$.  An important observation is that $A'(1) > 1$, which follows from the assumptions:~$A$ is convex, $A'(0) = 0$, and $A(1)=1$.  
%
%
Hence, since $A'(\tilde u) \approx A'(1) > 1$ for $x < 1$ due to \Cref{prop.pp_O(1)} (up to a constant shift), the comparison principle immediately yields
\be
	\sup_{x < 1} \tilde w(t,x)
		\leq \max\{ e^{-t/C}, \tilde w(t,1)\}
		\leq \frac{C}{t}.
\ee
The above two inequalities give (\ref{jul1004}). 


We use the above heuristics  to construct a suitable supersolution for $\tilde w$.   
Additionally, slightly abusing the notation, we let
\be
	\tilde w(t,x) = w(t, x + 2t - \tfrac{1}{2}\log(t+T) - C)
	\quad\text{ and }\quad
	\tilde u(t,x) = u(t, x + 2t - \tfrac{1}{2}\log(t+T) - C).
\ee
Here, we choose $C$ so that, in view of \Cref{prop.pp_O(1)},
\be\label{e.c42402}
	\inf_{t \geq 0, x \leq 10}
		\tilde u(t,x)
			\geq u_*,
\ee
with $u_*\in(0,1)$ chosen so that 	$A'(u_*) > 1$.
Note that it suffices to prove (\ref{jul1004}) for this shift of~$\tilde w$.

The function $\tilde w$ satisfies
\be\label{e.c42401}
	\tilde w_t
		- \Big(2 - \frac{1}{2(t+T)} \Big) \tilde w_x
			+ A'(\tilde u) \tilde w_x
			= \tilde w_{xx} + \tilde w (1 - A'(\tilde u)).
\ee
Let us define two auxiliary functions:
\be
	\theta(t) = \frac{T}{t+T}
\ee
and, for $\kappa \in (0,1/4)$ to be chosen,
\be
	\begin{split}
		\omega(t,x)
			&= x \frac{T}{t+T} \exp\Big\{4 - 2\sqrt{\frac{T}{t+T}} - \frac{x^2}{4(t+T)}\Big(1 - \kappa\sqrt{\frac{T}{t+T}}\Big)\Big\}
			\\&
			= x \theta(t) \exp\Big\{4 - 2 \sqrt{\theta(t)} - \frac{x^2}{4(t+T)}(1 - \kappa \sqrt{\theta(t)}) \Big\}.
	\end{split}
\ee
We also set
\be
	\bar w(t,x)
		= \begin{cases}
			\frac{T}{t+T}
				\qquad &\text{ for all } x \leq 1,\\
			\min\big\{\frac{T}{t+T}, e^{-x} \omega(t,x)\big\}
				\qquad &\text{ for all } x > 1.
		\end{cases}
\ee
Clearly, it suffices to show that there is $A>0$ such that
\be
	w(t,x) \leq A \overline w(t,x)
		\qquad\text{ for all } t\geq 1, x \in \R.
\ee
To this end, by the comparison principle, it is enough to show the following:
\begin{enumerate}[(i),itemsep=-1.75pt,topsep=5pt]
	\item $w(1,\cdot) \leq A \overline w(1,\cdot)$,
	\item $\bar w = e^{-x} \omega$ on $[10,\infty)$,
	\item $e^{-1}\omega(t,1) > \theta(t)$ for all $t\geq 0$,
	\item $\theta$ is a supersolution to~\eqref{e.c42401} on $(-\infty, 10]$, and
	\item $e^{-x} \omega$ is a supersolution to~\eqref{e.c42401} on $(0,\infty)$.
\end{enumerate}
Before proceeding, let us discuss the purpose of (ii) and (iiii).  The former allows us to only check (iv) on the domain $(-\infty,10]$.  
This is crucial because $\theta$ is {\em not} a supersolution to~\eqref{e.c42401} on $\R$
but only when $\tilde u \approx 1$.  On the other hand, (iii) guarantees the continuity of $\bar w$.  As such, we see that $\bar w$ is a minimum of two supersolutions on $[1,10]$ and, thus, is itself a supersolution. 
%
%

We now check conditions (i)-(v). 
As $w$ satisfies a parabolic equation with $L^1$ initial data:
\be
	\int |w(0,x)| \dx
		= \int w(0,x) \dx
		= \int (-\partial_x u_0 - \eta_*(u_0))\dx
		\leq \int (-\partial_x u_0) \dx
		= 1,
\ee
it can easily be bounded by
\be
	w(1,x)
		\leq C e^{-{x^2}/{C}}.
\ee
Up to increasing $T$ and $A$, it is clear that
\be
	C e^{- {x^2}/{C}}
		\leq A \overline w,
\ee
whence (i) is established.

Both (ii) and (iii) are elementary, by increasing $T$ if necessary.  
Additionally, after increasing $T$ further and recalling~\eqref{e.c42402}, we find
\be
		\theta_t
		- \Big(2 - \frac{1}{2(t+T)}\Big) \theta_x
		+ A'(\tilde u) \theta_x
		- \theta_{xx}
		- \theta(1 - A'(\tilde u))
		= - \frac{\theta}{t+T} + \theta (A'(u_*) - 1)
		> 0,
\ee
for $x<10$, which implies  (iv).

Finally, we check (v).  This computation  
is made easier by noting that it suffices to check that
\be\label{e.c42403}
	\omega_t + \frac{1}{2(t+T)} (\omega_x - \omega) + A'(\tilde u) \omega_x - \omega_{xx} 
		\geq 0.
\ee
A direct computation yields
\be\label{e.c42003}
	\begin{split}
		\omega_t + &\frac{1}{2(t+T)} (\omega_x - \omega)
		+ A'(\tilde u) \omega_x
		- \omega_{xx}
%
%
%
		\\&
		= \omega \Bigg[\frac{\dot\theta}{\theta}
			- \frac{\dot \theta}{\sqrt\theta}
			+\frac{x^2}{4(t+T)^2} (1 - \kappa \sqrt \theta)
			+ \frac{x^2}{4 (t+T)}  \kappa \frac{\dot\theta}{2\sqrt \theta}
		\\&\qquad\quad
			+\Big(\frac{1}{2(t+T)} + A'(\tilde u)\Big)\Big( \frac{1}{x}
				- \frac{x}{2(t+T)} (1 - \kappa\sqrt\theta)
				\Big)
				- \frac{1}{2(t+T)}
		\\&\qquad\quad
			- \Big( \Big( \frac{1}{x}- \frac{x}{2(t+T)} (1 - \kappa \sqrt\theta)\Big)
^2 + \frac{1}{x^2} + \frac{1}{2(t+T)}\Big(1 - \kappa \sqrt\theta\Big)\Big)
			\Bigg]
		\\&=
			\omega \Bigg[
			- \frac{\dot \theta}{\sqrt\theta}
			+\frac{x^2}{4(t+T)^2} \kappa \sqrt\theta \Big(1 - \kappa \sqrt\theta + \frac{\dot \theta}{2\theta}\Big)
		\\&\qquad\quad
			+\left(\frac{1}{2(t+T)} + A'(\tilde u)\right)\left( \frac{1}{x}
				- \frac{x}{2(t+T)} (1 - \kappa \sqrt\theta)
				\right)
  - \frac{3}{2(t+T)}\kappa \sqrt\theta
			\Bigg].
	\end{split}
\ee
Up to increasing $T$, we have
\be
	-\frac{\dot \theta}{\sqrt \theta}
		= \frac{\sqrt \theta}{t+T}
		> \frac{3 \kappa \sqrt\theta}{t + T}
	\qquad\text{ and }\qquad
	1 - \kappa \sqrt \theta + \frac{\dot\theta}{2 \theta}
		\geq \frac{1}{2}
\ee 
(recall that $\kappa \in (0,1/4)$ and $\theta \leq 1$).
Using these inequalities and that $A' \geq 0$,~\eqref{e.c42003} becomes
\be\label{e.c42501}
	\begin{split}
		\omega_t + &\frac{1}{2(t+T)} (\omega_x - \omega)
		+ A'(\tilde u) \omega_x
		- \omega_{xx}
			\\&
			\geq \omega \Bigg[\frac{\sqrt \theta}{2(t+T)}
				+ \frac{x^2}{4(t+T)^2} \frac{\kappa \sqrt \theta}{2}
			+\left(\frac{1}{2(t+T)} + A'(\tilde u)\right)\left( \frac{1}{x}
				- \frac{x}{2(t+T)}
				\right)
			\Bigg].
	\end{split}
\ee
The first and second terms in the right side of~\eqref{e.c42501} are positive.  Next,
we show that they dominate the potentially negative last term.

When $x \leq \sqrt{2(t+T)}$, the final term in~\eqref{e.c42501} is positive as well.  We obtain
\be\label{e.c42503}
	\omega_t + \frac{1}{2(t+T)} (\omega_x - \omega)
		+ A'(\tilde u) \omega_x
		- \omega_{xx}
		\geq 0
			\quad\text{ when } x \leq \sqrt{2(t+T)}, 
\ee
as desired.  Hence, we need only consider the case $x > \sqrt{2 (t+T)}$.

In this case, we first show that $A'(\tilde u)$ is small when $x > \sqrt{2(t+T)}$.  
Indeed, using \Cref{prop.pp_O(1)} and (\ref{e.u_asymptotics}), 
we see that, up to increasing $T$, we have
\be
	\begin{split}
	\tilde u(t, x)
		&< \tilde u(t, 2 \sqrt{t+T})\leq C \exp(- \sqrt{t+T}).
	\end{split}
\ee
Due to the assumptions on $A$, we then deduce
\be\label{e.c42502}
	A'(\tilde u(t,x))
		\leq C \tilde u(t,x)
		\leq C \exp(-\sqrt{t+T})
		\leq \frac{1}{2(t+T)}.
\ee
Applying~\eqref{e.c42502} and Young's inequality, we have
\be
	\begin{split}
		\Big(\frac{1}{2(t+T)} &+ A'(\tilde u)\Big)\Big( \frac{1}{x}
				- \frac{x}{2(t+T)}
				\Big)
			\geq
				-\frac{x}{4(t+T)^2}
			\\&
			\geq -\frac{x^2}{8(t+T)^{5/2}}
				- \frac{1}{8(t+T)^{3/2}}
			= - \frac{1}{\sqrt T}\left(\frac{x^2}{8(t+T)^2} \sqrt \theta
				+ \frac{\sqrt\theta}{2(t+T)} \right).
	\end{split}
\ee
Clearly, up to increasing $T$, this, along with~\eqref{e.c42501} implies that
\be\label{e.c42504}
	\omega_t + \frac{1}{2(t+T)}(\omega_x - \omega) + A'(\tilde u) \omega_x - \omega_{xx}
		\geq 0
			\quad\text{ when } x \geq \sqrt{2 (t+T)}.
\ee

The combination of~\eqref{e.c42503} and~\eqref{e.c42504} yields~\eqref{e.c42403}.  Due to the equivalence of~\eqref{e.c42403}, this yields (v), which concludes the proof.~$\Box$

\section{Semi-FKPP fronts: proof of \Cref{t.rde-intro}.(i)}
 \label{sec:sfkpp}

In this section, we show that there is $x_0$, depending on the initial condition, 
such that solutions $u$ to~\eqref{jul1314} with $\chi \in[0, 1)$, satisfy
\be
	u(t,x+ 2t - \tfrac{3}{2}\log t + x_0) \to U_*(x),
		\qquad\text{ as } t\to\infty.
\ee
The general structure of this proof is similar to the previous section.  First, we obtain upper and lower bounds of $u$ that show that the front $2t - (3/2)\log t + O(1)$.  Second, we examine the behavior at $2t + t^\gamma$ to identify a candidate for the precise constant shift term.  Finally, we use the `pulled' nature of the front to show that closeness of the $u$ and the traveling wave $O(t^\gamma)$ ahead of the front yields convergence everywhere.

The first two steps proceed as in the analogous case ($\beta < 2$) in the Burgers-FKPP setting
of~\cite{AHR}, with the use of the weighted Hopf-Cole transform in \Cref{prop-mar902}. This is recalled
briefly in Section~\ref{sec:sfkpp_tgamma}, where Proposition~\ref{p.sfkpp_tgamma} is proved.

On the other hand, the final step, convergence to the traveling wave, is significantly more difficult than in~\cite{AHR}.  This is mostly due to the fact that the equation for the weighted Hopf-Cole transform of $u$ involves terms that did not appear in~\cite{AHR} where $\alpha(u) = u$. This part of the
proof is presented in Section~\ref{sec:sfkpp_conv}.
The key result here is 
the differential inequality in Lemma~\ref{l.s_eqn} that follows from yet another surprising and intricate set
of algebraic cancellations.

In this section, we use the tilde for a change to the moving frame~$x \mapsto x + 2t - (3/2)\log(t+1)$, such as
\be
	\tilde u(t,x)
		:= u(t, x + 2t - \frac{3}{2}\log(t+1)).
\ee

\subsection{The precise behavior at $2t + t^\gamma$}
\label{sec:sfkpp_tgamma}

The goal of this section is to establish the following proposition, the analogue of \cite[Lemma~4.1]{AHR}.
\begin{prop}\label{p.sfkpp_tgamma}
	Under the assumptions of \Cref{t.rde-intro}.(i), there exists $\tgammacoefficient_\infty>0$, depending on the initial condition, such that, for any $\gamma \in (0,1/2)$, we have
	\be
		\lim_{t\to\infty} \frac{e^{t^\gamma}}{t^\gamma} \tilde u(t, t^\gamma)
			= \tgammacoefficient_\infty.
	\ee
\end{prop}
We discuss below the essential ingredients of \cite[Lemma~4.1]{AHR} and show that they are present in our setting.
We begin with a preliminary upper bound that is required to show that the integral term in the weighted Hopf-Cole transform is not `too big.'  Notice that, by the comparison principle,
\be
	u \leq u_{\rm pp},
\ee
where $u_{\rm pp}$ is the pushmi-pullyu front associated to $\chi =1$, with the same
initial condition. It follows immediately from \Cref{t.rde-intro}.(ii) and~\eqref{e.u_asymptotics} that
\be\label{e.c42705}
	\tilde u(t,x)
		\leq \min(1, C\exp\{- x + \log (t+1)\})
		= \min(1, C (t+1) e^{-x}).
\ee
We deduce that
\be\label{e.c42701}
	\int_x^\infty \tilde u(t,y) \dy
		\leq \begin{cases}
				C + \log(t+1) - x
					\qquad&\text{ if } x \leq \log (t+1),\\
				C(t+1) e^{-x}
					\qquad&\text{ if } x \geq \log (t+1).
			\end{cases}
\ee

Next, we work with the weighted Hopf-Cole transform of $\tilde u$:
\be
	v(t,x)
		= \exp\Big( x + \sqrt \chi \int_x^\infty \alpha(\tilde u(t,y)) \dy\Big) \tilde u(t,x).
\ee
Following the computations in the proof of \Cref{prop-mar902}, we see that
\be\label{e.c42704}
	v_t - v_{xx} + \frac{3}{2(t+1)} (v_x - v) - \frac{g[\tilde u]}{\tilde u} v = 0,
\ee
where
\be\label{e.c42702}
	\begin{split}
	g[\tilde u]: =  &f(\tilde u) - \tilde u - \chi \alpha(\tilde u)^2 \tilde u + 2 \sqrt \chi \alpha(\tilde u) \tilde u_x
		\\&
		+\tilde u\Big(- \sqrt \chi \int_x^\infty \alpha''(\tilde u) \tilde u_y^2 \dy + \sqrt \chi \int_x^\infty\alpha'(\tilde u) f(\tilde u) \dy \Big)
		\leq 0.
	\end{split}
\ee
Using~\eqref{e.c42701}, we make two crucial observations:
\be\label{e.c42703}
	\begin{split}
		v(t, -t^\gamma)
			&\leq C \exp\big( - \tfrac{1}{2}(1-\sqrt \chi) t^\gamma\big) \qquad \text{ and }\\
		\frac{g[\tilde u](t, t^\gamma)}{\tilde u(t,t^\gamma)}
			&\geq - C \exp\big( - \tfrac{1}{2} t^\gamma\big),
	\end{split}
\ee
for any $\gamma > 0$.

The key ingredients in the proof of~\cite[Lemma~4.1]{AHR} are
\begin{enumerate}[(i),itemsep=-1.75pt,topsep=5pt,]
	\item (subsolution) $v$ satisfies
		\be
			v_t - v_{xx} + \frac{3}{2(t+1)} (v_x - v) \leq 0;
		\ee
	\item (not too far from supersolution) for any $\gamma \in (0,1/2)$ and $x > t^\gamma$,
		\be
			v_t - v_{xx} + \frac{3}{2(t+1)} (v_x - v) \geq - c e^{- t^\gamma/C};
		\ee
	\item (approximate Dirichlet boundary condition) for any $\gamma>0$, 
	\be\label{e.c42923}
		v(t, -t^\gamma) \leq C e^{- t^\gamma/C}.
	\ee
\end{enumerate}
It is clear that (i) follows from~\eqref{e.c42704}-\eqref{e.c42703} while (ii) and (iii) follow from~\eqref{e.c42703}.

From this point, the proof of \Cref{p.sfkpp_tgamma} is verbatim the same as \cite[Lemma~4.1]{AHR} 
and is omitted.

%
%
%

\subsection{Convergence to the wave}
\label{sec:sfkpp_conv}

\subsubsection{The setup}

We now use \Cref{p.sfkpp_tgamma} to finish the proof of \Cref{t.rde-intro}.(i).  
While we follow the approach of~\cite{AHR}, there is a technical complication due 
to extra terms that appear in the general setting, such as the nontiviality of $\alpha''$.

To begin, we make two observations.  Firstly, the combinations of \Cref{p.sfkpp_tgamma} and~\eqref{e.c42705} yields
\be\label{e.v_tgamma}
	\lim_{t\to\infty} \frac{v(t,t^\gamma)}{t^\gamma} = \tgammacoefficient_\infty.
\ee
In fact, this is the bound we use, not \Cref{p.sfkpp_tgamma}.   Secondly, since
\be
	\tilde  u v_x
		= v \tilde u_x +   (1 -\sqrt{\chi}\alpha(\tilde u))\tilde uv,
\ee
we can rewrite the equation for $v$ as
\be\label{e.c42710}
	v_t
		- v_{xx}
		+ \frac{3}{2 t} (v_x - v)
		- 2 \sqrt \chi \alpha(\tilde u) v_x
		- G[\tilde u] v
		= 0.
\ee
Here, for a suitably smooth and decaying function $\rho$, we define
\be
	\begin{split}
		G[\rho]
			&:= \frac{g[\rho]}{\rho} - \frac{2 \sqrt \chi \alpha(\rho) \rho_x}{\rho} - 2 \sqrt \chi \alpha(\rho) (1 - \sqrt \chi \alpha(\rho)) 
			\\&
			= \frac{f(\rho)}{\rho} - 1
				- \chi \alpha(\rho)^2
				- 2 \sqrt\chi \alpha(\rho) (1 - \sqrt \chi\alpha(\rho))
		\\&
		\quad - \sqrt \chi \int_x^\infty \alpha''(\rho) \rho_y^2 \dy + \sqrt \chi \int_x^\infty\alpha'(\rho) f(\rho) \dy.
	\end{split}
\ee 
With some arithmetic, this can be simplified to the form we use below:
\be\label{e.c42922}
	\begin{split}
	G[\rho]
		= &- (1+\chi) \alpha(\rho)
			- 2 \sqrt \chi(1 - \sqrt \chi) \alpha(\rho)
			+ \chi \rho \alpha'(\rho)(1-\alpha(\rho))
		\\&
			- \sqrt \chi \int_x^\infty \alpha''(\rho) \rho_y^2 \dy + \sqrt \chi \int_x^\infty\alpha'(\rho) f(\rho) \dy
	\end{split}
\ee
The re-writing of~\eqref{e.c42710} in terms of $G$, at the expense of introducing
an extra drift term in the left side of ~\eqref{e.c42710} compared to (\ref{e.c42704}), 
is crucial in the proof of 
the main  \Cref{l.sfkpp_s}, below.  This algebraic ingredient is a major difference with the proofs of \cite{AHR,NRR1,NRR2}.   

We now take the weighted Hopf-Cole transform of a suitable shift of the traveling wave.  
To this end, for $\tgammacoefficient \in (\tgammacoefficient_\infty/2, 2\tgammacoefficient_\infty)$, we let
\be\label{jul2802}
	\vphi_\tgammacoefficient(t,x)= U_*(x+\zeta_\tgammacoefficient(t)),
\ee
with $\zeta_\tgammacoefficient$ to be chosen.   
This satisfies
\be\label{e.c6295}
	\begin{aligned}
		&\partial_t \varphi_\tgammacoefficient - \partial_x^2\varphi_\tgammacoefficient - \left(2 - \frac{3}{2(t+1)}\right) \partial_x \varphi_\tgammacoefficient 
			- f(\varphi_\tgammacoefficient)
			= \left( \frac{3}{2(t+1)} + \dot \zeta_\tgammacoefficient\right) \partial_x \varphi_\tgammacoefficient.
	\end{aligned}
\ee
Next, we define its Hopf-Cole transform:
\be\label{e.c714}
	\psi_\tgammacoefficient(t,x)
		= \exp\Big(x + \sqrt{\chi}\int_x^\infty \alpha( \vphi_\tgammacoefficient(t,y))\dy\Big) \vphi_\tgammacoefficient(t,x),
\ee
and fix the shift $\zeta_\tgammacoefficient(t)$ by the
normalization 
\be\label{e.c6291}
	\psi_\tgammacoefficient(t,t^\gamma)
	= \tgammacoefficient t^\gamma.
\ee
Recall that there are constants $D>0$ and $B\in \R$ so that the traveling wave has the asymptotics
\be\label{e.c42708}
U_*(x)= D x e^{-x} + Be^{-x} + o(e^{-x})
\qquad\text{ for } x \gg 1.
\ee
Thus, we have 
\be\label{e.zeta}
\zeta_\tgammacoefficient(t)
= - \log\big(\frac{\tgammacoefficient}{D}\big) + \farc{1}{t^\gamma}
\Big(\frac{B}{D} - \log\big(\frac{\tgammacoefficient}{D}\big)\Big)
			+ o(t^{-\gamma}),
\ee
and 
\be\label{e.c6297}
	|\dot \zeta_\tgammacoefficient(t)|\leq \frac{C}{t^{1 + \gamma}}.
\ee
The constant $C$ is independent of $\tgammacoefficient$ due to the restriction $\tgammacoefficient \in (\tgammacoefficient_\infty/2, 2\tgammacoefficient_\infty)$.  From~\eqref{e.c6295} and~\eqref{e.c6297}, we find
\be\label{e.c6298}
	\Big|\partial_t \psi_\tgammacoefficient
		- \partial_x^2\psi_\tgammacoefficient
		+ \frac{3}{2 t}(\partial_x \psi_\tgammacoefficient-\psi_\tgammacoefficient)
		- 2 \sqrt \chi \alpha(\vphi_\tgammacoefficient) \partial_x \psi_\tgammacoefficient
		- G[\vphi_\tgammacoefficient] \psi_\tgammacoefficient\Big|
    	\leq \frac{C}{t^{1-\gamma}}
		\quad \text{ for } |x| < t^\gamma.
\ee

\subsubsection{The main lemma and heuristic description of the proof of \Cref{t.rde-intro}.(i)}

The main step in the proof of \Cref{t.rde-intro}.(i) is 
the following analogue of \cite[Lemma~4.3]{AHR} bounding the difference between the
weighted Hopf-Cole transforms. 
\begin{lem}\label{l.sfkpp_s}
Under the assumptions of \Cref{t.rde-intro}.(i), if $\eps \in (0, \tgammacoefficient_\infty/2)$ 
and if $\gamma < 1/3$, then there is $\lambda > 0$ so that
\be\label{jul1102}
	v(t,x) - \psi_{\tgammacoefficient_\infty+\eps}(t,x)
			\leq \frac{C}{t^\lambda}
		\qquad\text{ for all } t > T_\eps
			\text{ and } |x| < t^\gamma.
	\ee
\end{lem}
The proof of \Cref{l.sfkpp_s} is substantially more intricate than \cite[Lemma~4.3]{AHR}.
We postpone it  to \Cref{s.main_lemma}.

A simple consequence of \Cref{l.sfkpp_s}  useful in the proof of \Cref{t.rde-intro}.(i) is that, 
for any $L >0$, we have
\be\label{e.c42903}
	\sup_{|x| \leq L} v(t,x) \leq C.
\ee


\subsubsection*{\bf Proof of \Cref{t.rde-intro}.(i)}

While the proof of \Cref{t.rde-intro}.(i) is very similar to that of \cite[Theorem 1.1 for $\beta < 2$]{AHR}, we give the complete details due to the subtle but important 
difference between \Cref{l.sfkpp_s} and \cite[Lemma 4.3]{AHR}; that is, that we use $v$ directly, here, whereas~\cite{AHR} uses an auxiliary function.

The two key inequalities that we prove are: for any $L, \eps>0$,
\be\label{e.c42908}
	\inf_{x< t^\gamma}(\tilde u - \vphi_{\tgammacoefficient_\infty - \eps}) \geq 0
	\quad\text{ and }\quad
	\lim_{t\to\infty}\sup_{x \in [-L,L]}(\tilde u - \vphi_{\tgammacoefficient_\infty + \eps})
		\leq C \eps.
\ee
Before establishing these, we show how to use them to conclude the proof.

We begin with an observation:  
\be\label{e.c42706}
	\Big\|\frac{e^{x_+}}{x_+ + 1}(\vphi_{\tgammacoefficient_\infty} - \vphi_{\tgammacoefficient_\infty-\eps})\Big\|_{L^\infty}
		\leq C \eps,
\ee
due to the regularity of the traveling wave $U_*$, the definition of $\varphi_\omega$~\eqref{jul2802}, and the expansion of $\zeta_\omega$~\eqref{e.zeta}.

We immediately see from~\eqref{e.c42908} and~\eqref{e.c42706} that
\be\label{e.c72701}
	\lim_{t\to\infty} \|\tilde u - \varphi_{\omega_\infty} \|_{L^\infty([-L,L])}
		= 0.
\ee
To handle the domain $(-\infty,-L)$, it suffices to use both inequalities in~\eqref{e.c42706} to find: for all $x \in (-\infty, -L)$ and $t$ sufficiently large,
\be\label{e.c72703}
	\begin{split}
	|\varphi_{\omega_\infty}(t,x) - \tilde u(t,x)|
		&\leq \|\varphi_{\omega_\infty} - \varphi_{\omega_\infty - \eps}\|_{L^\infty}
			+ |\varphi_{\omega_\infty-\eps}(t,x) - \tilde u(t,x)|
		\\& \leq  C \eps + |\varphi_{\omega_\infty-\eps}(t,x) - \tilde u(t,x)|
		= C \eps + \tilde u(t,x) - \varphi_{\omega_\infty-\eps}(t,x)
		\\&
		\leq C \eps + 1 - \varphi_{\omega_\infty-\eps}(t,x)
		\leq C \eps + C e^{-\alpha'(1) L}.
	\end{split}
\ee
To handle the domain $(L,\infty)$, we argue similarly to find: for all $x \in (L, +\infty)$ and $t$ sufficiently large,
\be\label{e.c72702}
	|\varphi_{\omega_\infty}(t,x) - \tilde u(t,x)|
		\leq C\eps + C e^{-L}.
\ee

The claim is established after putting together~\eqref{e.c72701},~\eqref{e.c72702}, and~\eqref{e.c72703} and then sequentially taking $t$ to infinity, $\eps$ to zero, and $L$ to infinity.  Hence, our goal is now to prove~\eqref{e.c42908}.

We begin with the lower bound in~\eqref{e.c42908}.    For $t$ sufficiently large, it 
follows from~(\ref{e.c6291}) and the smallness of the integrals in the respective
definitions of the 
weighted Hopf-Cole transforms that
\be
	\tilde u(t,t^\gamma) > \vphi_{\tgammacoefficient_\infty-\eps}(t,t^\gamma).
\ee
Since $\tilde u$ is steeper than $U_*$ and, thus, than $\vphi_{\tgammacoefficient_\infty-\eps}$, 
we deduce that
\be\label{e.c42904}
	\tilde u(t,x) > \vphi_{\tgammacoefficient_\infty-\eps}(t,x)
		\qquad\text{ for all } x < t^\gamma,
\ee
which concludes the proof of the first inequality in~\eqref{e.c42908}.

Next,  we consider the second inequality in~\eqref{e.c42908}. 
We decompose the difference as 
\be\label{e.c42720}
	\begin{split}
		\tilde u&(t,x) - \vphi_{\tgammacoefficient_\infty + \eps}(t,x)
			\\&
			= \exp\Big\{- x - \sqrt \chi \int_x^\infty \alpha(\vphi_{\tgammacoefficient_\infty + \eps}(t,y))\dy\Big\}
				\big[ v(t,x) - \psi_{\tgammacoefficient_\infty+\eps}(t,x)\big]
			\\&
			\quad + e^{-x} v(t,x)
				\Big[ \exp\Big\{- \sqrt \chi\int_x^\infty \alpha( \tilde u(t,y)) \dy\Big\}
				- \exp\Big\{- \sqrt \chi\int_x^\infty \alpha(\vphi_{\tgammacoefficient_\infty +\eps}(t,y)) \dy\Big\}\Big].
	\end{split}
\ee
The first term is easy to estimate using \Cref{l.sfkpp_s}:
\be\label{e.c42721}
	\exp\Big\{- x - \sqrt \chi \int_x^\infty \alpha(\vphi_{\tgammacoefficient_\infty + \eps}(t,y))\dy\Big\}
				\big[ v(t,x) - \psi_{\tgammacoefficient_\infty+\eps}(t,x)\big]
		\leq e^L \frac{C}{t^\lambda},
\ee
as desired.

Next, we consider the second term in~\eqref{e.c42720}.  Using~\eqref{e.c42705}, it is clear that, up to an error term, we can restrict the limits of integration slightly:
\be
	\begin{split}
	e^{-x} v(t,x)
		&\Big[ \exp\Big\{- \sqrt \chi\int_x^\infty \alpha( \tilde u(t,y)) \dy\Big\}
		- \exp\Big\{- \sqrt \chi\int_x^\infty \alpha(\vphi_{\tgammacoefficient_\infty +\eps}(t,y)) \dy\Big\}\Big]
		\\&
		\leq C e^{- \frac{1}{2} t^\gamma} + e^{-x} v(t,x)
			\Big[ \exp\Big\{- \sqrt \chi\int_x^{t^\gamma} \alpha( \tilde u(t,y)) \dy\Big\}
		\\&\qquad
		- \exp\Big\{- \sqrt \chi\int_x^{t^\gamma} \alpha(\vphi_{\tgammacoefficient_\infty + \eps}(t,y)) \dy\Big\}\Big].
	\end{split}
\ee
We further decompose the second term above
\be
	\begin{split}
		&\exp\Big\{- \sqrt \chi\int_x^{t^\gamma} \alpha( \tilde u(t,y)) \dy\Big\}
		- \exp\Big\{- \sqrt \chi\int_x^{t^\gamma} \alpha(\vphi_{\tgammacoefficient_\infty + \eps}(t,y)) \dy\Big\}
		\\&\quad
		=
			\Big[ \exp\Big\{- \sqrt \chi\int_x^{t^\gamma} \alpha( \tilde u(t,y)) \dy\Big\}
				-  \exp\Big\{- \sqrt \chi\int_x^{t^\gamma} \alpha(\vphi_{\tgammacoefficient_\infty - \eps}(t,y)) \dy\Big\}\Big]
		\\&\qquad\quad
			+ \Big[\exp\Big\{- \sqrt \chi\int_x^{t^\gamma} \alpha(\vphi_{\tgammacoefficient_\infty - \eps}(t,y)) \dy\Big\}
		- \exp\Big\{- \sqrt \chi\int_x^{t^\gamma} \alpha(\vphi_{\tgammacoefficient_\infty + \eps}(t,y)) \dy\Big\}\Big].
	\end{split}
\ee
The first term in the right side is nonpositive because $\alpha$ is increasing and $\tilde u \geq \varphi_{\tgammacoefficient_\infty - \eps}$ (recall~\eqref{e.c42904}). 
For the second term, we use a Taylor expansion and~\eqref{e.c42706} to conclude that
\be
	\begin{split}
	&\exp\Big\{- \sqrt \chi\int_x^{t^\gamma} \alpha(\vphi_{\tgammacoefficient_\infty - \eps}(t,y)) \dy\Big\}
		- \exp\Big\{- \sqrt \chi\int_x^{t^\gamma} \alpha(\vphi_{\tgammacoefficient_\infty + \eps}(t,y)) \dy\Big\}
		\\&\quad
		=	e^{- \sqrt \chi\int_x^{t^\gamma} \alpha(\vphi_{\tgammacoefficient_\infty - \eps}(t,y)) \dy}\Big[1
		- \exp\Big\{- \sqrt \chi\int_x^{t^\gamma} (\alpha(\vphi_{\tgammacoefficient_\infty + \eps}(t,y)) -\alpha(\vphi_{\tgammacoefficient_\infty - \eps}(t,y)) )\dy\Big\}\Big]
		\\&\quad
		\leq 1
		- \exp\Big\{- \sqrt \chi\int_x^{t^\gamma} (\alpha(\vphi_{\tgammacoefficient_\infty + \eps}(t,y)) -\alpha(\vphi_{\tgammacoefficient_\infty - \eps}(t,y)) )\dy\Big\}
		\\&\quad
		\leq \sqrt \chi\int_x^{t^\gamma} \big(\alpha(\vphi_{\tgammacoefficient_\infty + \eps}(t,y)) -\alpha(\vphi_{\tgammacoefficient_\infty - \eps}(t,y))\big)\dy
		\leq C \eps.
	\end{split}
\ee
Using the boundedness of $e^{-x}v$~\eqref{e.c42903} as well as the work above, we find
\be\label{e.c42907}
	\begin{split}
	e^{-x} v(t,x)
		&\Big[ \exp\Big\{- \sqrt \chi\int_x^\infty \alpha( \tilde u(t,y)) \dy\Big\}
		- \exp\Big\{- \sqrt \chi\int_x^\infty \alpha(\vphi_{\tgammacoefficient_\infty +\eps}(t,y)) \dy\Big\}\Big]
		\\&\quad
		\leq C e^{-\frac{1}{2}t^\gamma} + \frac{C}{t^\lambda} + C \eps.
	\end{split}
\ee
By applying~\eqref{e.c42721} and~\eqref{e.c42907} in~\eqref{e.c42720}, we obtain
\be\label{e.c42722}
	\tilde u(t,x)
		- \varphi_{\tgammacoefficient_\infty + \eps}
		\leq C e^{-\frac{1}{2}t^\gamma}+\frac{C}{t^\lambda} + C \eps.
\ee
This yields exactly the second inequality in~\eqref{e.c42908}, which completes the proof.~$\Box$

\subsubsection{Proof of Lemma~\ref{l.sfkpp_s}}
\label{s.main_lemma}

A key step in proving Lemma~\ref{l.sfkpp_s} is the derivation of a differential inequality for
\be
	s(t,x) = v(t,x) - \psi_{\tgammacoefficient_\infty + \eps}(t,x).
\ee
As the computation is somewhat complicated, we state it here as a further lemma and prove it after we show how it is used to prove \Cref{l.sfkpp_s}.
\begin{lem}\label{l.s_eqn}
	In the setting of \Cref{l.sfkpp_s}, whenever $s \geq 0$, we have 
	\be\label{e.c42711}
		s_t - s_{xx} + \frac{3}{2t} (s_x - s)
			- 2 \sqrt \chi \alpha (\tilde u) s_x
			\leq \frac{C}{t^{1-\gamma}},
			\qquad\text{ for } |x|\leq t^\gamma.
	\ee
\end{lem}

\subsubsection*{Proof of \Cref{l.sfkpp_s}}  

We prove this via the construction of a supersolution for 
(\ref{e.c42711}) on $[-t^\gamma, t^\gamma]$.  First, we check the boundary conditions for $s$.  Using~\eqref{e.c42923}, we find
\be\label{e.c42713}
	s(t, -t^\gamma)
		\leq C \exp\big( - \tfrac{1}{C} t^\gamma\big).
\ee
Additionally, using~\eqref{e.v_tgamma} and~\eqref{e.c6291} and the choice of $\gamma = \tgammacoefficient_\infty+\eps$ with $\eps >0$, we find
\be\label{e.c42714}
	s(t,t^\gamma)
		\leq 0.
\ee

We now define a supersolution to~\eqref{e.c42711}. Fix any $\lambda \in (0,1 - 4\gamma)$ and $Q\geq 1$, and let
\be
	\overline s(t,x)
		= \frac{Q}{t^\lambda} \cos \Big( \frac{x + t^\gamma}{2 t^\gamma}\Big)
		= \frac{Q}{t^\lambda} \cos \Big( \frac{x}{2 t^\gamma} + \frac{1}{2}\Big).
\ee
We first notice that $\overline s_x \leq 0$ on $[-t^\gamma, t^\gamma]$.  This allows us to drop a 
positive term in the first step below, after which we directly compute:
\be
	\begin{split}
		\overline s_t - &\overline s_{xx}
			+ \frac{3}{2t} (\overline s_x - \overline s)
			- 2\sqrt\chi \alpha(\tilde u) \overline s_x
			\geq \overline s_t - \overline s_{xx}
			+ \frac{3}{2t} (\overline s_x - \overline s)
			\\&
			= -\frac{\lambda}{t} \overline s
				+ \frac{\gamma Q x}{2 t^{1 + \lambda + \gamma}}  \sin\Big( \frac{x}{2 t^\gamma} + \frac{1}{2}\Big)
				+ \frac{1}{4t^{2\gamma}} \overline s
				+ \frac{3}{2t} \Big(-\frac{xQ}{2t^{\gamma+\lambda}}\sin \Big( \frac{x}{2 t^\gamma} + \frac{1}{2}\Big)
					- \overline s \Big).
	\end{split}
\ee
Taking $t$ sufficiently large and then using that, for $|x| \leq t^\gamma$, we have 
$\overline s(t,x) \geq Q\cos(1)/t^\lambda$, we see
\be
	\begin{split}
		\overline s_t - &\overline s_{xx}
			+ \frac{3}{2t} (\overline s_x - \overline s)
			- 2\sqrt\chi \alpha(\tilde u) \overline s_x
			\geq \frac{1}{8 t^{2\gamma}} \overline s
				- \frac{CQ}{t^{1+\lambda}}
			\geq \frac{Q\cos(1)}{8 t^{2\gamma + \lambda}}
				- \frac{CQ}{t^{1+\lambda}}
			\geq \frac{Q\cos(1)}{16 t^{2\gamma + \lambda}}.
	\end{split}
\ee
The last inequality follows because $2\gamma + \lambda < 1 +\lambda$ (recall that $\gamma < 1/3$).  Recalling that $Q\geq 1$ and decreasing $\lambda$ so that $2\gamma + \lambda < 1 -\gamma$ yields
\be\label{e.c42910}
	\overline s_t - \overline s_{xx}
			+ \frac{3}{2t} (\overline s_x - \overline s)
			- 2\sqrt\chi \alpha(\tilde u) \overline s_x
		\geq \frac{C}{t^{1-\gamma}}
\ee
for $t$ sufficiently large. 
In other words, $\overline s$ is a supersolution to~\eqref{e.c42711} for $t$ sufficiently large (independent of $Q$) and $|x| < t^\gamma$.  We let $T_0$ be such that~\eqref{e.c42910} holds for $t\geq T_0$.

By \Cref{l.s_eqn} and the comparison principle, we have
\be\label{e.c42921}
	s(t,x) \leq \overline s(t,x)
		\qquad\text{ on } |x| \leq t^\gamma
\ee
for all $t\geq T_0$ sufficiently large as long as we can verify that 
$s \leq \overline s$ on the parabolic boundary.

We first check the portion of the parabolic boundary $t=T_0$.  There, it is easy to see that we can choose $Q$ such that $s(T_0,x) \leq \overline s(T_0,x)$ for all $|x| < T_0^\gamma$; indeed, $s(T_0,\cdot)$ is bounded above and $\overline s/Q$ is bounded below on $[-t^\gamma,t^\gamma]$.

Next we check the portion of the parabolic boundary $|x| = t^\gamma$.  Notice that
\be
	\overline s(t, \pm t^\gamma)
		\geq \frac{Q \cos(1)}{t^\lambda}.
\ee
From~\eqref{e.c42713} and~\eqref{e.c42714}, we see that
\be
	s(t, \pm t^\gamma)
		\leq C e^{-t^\gamma /C}.
\ee
Clearly $s(t, \pm t^\gamma) \leq \overline s(t,\pm t^\gamma)$ for $t\geq T_1$ for some $T_1$.  Hence, up to increasing $Q$ to handle the range~$t \in [T_0, \max\{T_0, T_1\}]$, we have
\be\
	s(t, \pm t^\gamma) \leq \overline s(t,\pm t^\gamma)
		\qquad \text{ for all } t \geq T_0.
\ee

The previous two paragraphs show that $s \leq \overline s$ on the parabolic boundary.  The inequality~\eqref{e.c42921} follows, from which the conclusion of \Cref{l.sfkpp_s} follows.~$\Box$

\subsubsection*{Proof of \Cref{l.s_eqn}}

We now establish the differential inequality~\eqref{e.c42711}.  As $\tgammacoefficient$ remains 
fixed throughout the proof below, we drop it notationally, for ease.
We start by obtaining a few preliminary results that help us establish~\eqref{e.c42711}.  First, note that
\be\label{e.c42802}
	\tilde u > \vphi
		\quad\text{ whenever }\quad
	s > 0.
\ee
Indeed, take $x_t$ to be the `farthest right' point where $s(x_t) = 0$ (it is clearly negative at $+\infty$, so this is well defined).  Then, we have
\be
	\tilde u(t,x_t) \exp\Big(x_t +  \sqrt{\chi}\int_{x_t}^\infty \alpha(\tilde u(t,y)) \dy\Big)
		= \vphi(t,x_t) \exp\Big(x_t + \sqrt{\chi}\int_{x_t}^\infty \alpha(\vphi(t,y)) \dy\Big).
\ee
Suppose that $\tilde u(t,x_t) < \vphi(t,x_t)$.  In this case, 
as $\tilde u$ is steeper than $\vphi$, it follows that
$\tilde u < \vphi$ on~$(x_t,\infty)$, so we find, using monotonicity of $\alpha(u)$,
\be
	\tilde u(t,x_t)
		< \vphi(t,x_t)
		= \tilde u(t,x_t) \exp\Big( -\sqrt{\chi} \int_{x_t}^\infty (\alpha(\vphi(t,y)) - \alpha(\tilde u(t,y))) \dy\Big)
		< \tilde u(t,x_t),
\ee
which is clearly a contradiction. It follows that $\tilde u(t,x_t) \ge\vphi(t,x_t)$, and 
another use of the steepness
comparison together with the choice of $x_t$ implies that
\be\label{e.c72704}
	\tilde u(t,x) > \varphi(t,x)
	\quad\text{ for all $x < x_t$}.
\ee
Now, by construction of $x_t$, if $s(t,x)>0$ it must be that $x< x_t$.  Thus, by~\eqref{e.c72704}, $\tilde u(t,x) > \varphi(t,x)$, as desired.  This completes the proof of the claim (\ref{e.c42802}).

Next, since $u(t,\cdot)$ and $U_*$ are strictly decreasing, we can find the function
$\eta(t,u)$ 
such that
\be
	- u_x(t,x)
		= \eta(t,u(t,x)).
\ee
One can see immediately that $\eta(t,u(t,x)) = w(t,x) + \eta_*(u(t,x))$, and the positivity of $w$ implies that, 
for any $z\in(0,1)$,
\be\label{e.c42902}
	\eta(t,z) > \eta_*(z).
\ee
In addition, Proposition~\ref{dec28-prop06} and Lemma~\ref{lem-jun1402}
yield 
\be\label{e.phase_plane}
	f(z) = \eta_*(z) (2 - \eta_*'(z))
		\quad\text{ and }\quad
	\sqrt \chi (z - A(z))	
		\leq \eta_*(z)
		\leq z - A(z).
\ee

We are now ready to prove~\eqref{e.c42711}.  By combining~\eqref{e.c42710} and~\eqref{e.c6298}, we see that
\be\label{e.c42801}
	\begin{split}
	s_t
		- s_{xx}
		&+ \frac{3}{2 t} (s_x - s)
		- 2 \sqrt \chi \alpha(\tilde u) s_x
		\\&
		\leq \frac{C}{t^{1-\gamma}}
			+ 2 \sqrt \chi (\alpha(\tilde u) - \alpha(\vphi)) \psi_x
			+ G[\tilde u] v
			- G[\vphi]
			\psi
		\\&
		\leq \frac{C}{t^{1-\gamma}}
			+ 2 \sqrt \chi (\alpha(\tilde u) - \alpha(\vphi)) \psi_x
			+ G[\tilde u] s
			+ \Big( G[\tilde u]
			- G[\vphi]\Big)
			\psi.
	\end{split}
\ee
Thus, it suffices to prove
\be\label{e.c42901}
		2 \sqrt \chi (\alpha(\tilde u) - \alpha(\vphi)) \psi_x
			+ G[\tilde u] s
			+ \Big( G[\tilde u]
			- G[\vphi]\Big)
			\psi
			\leq 0.
\ee
This is the reason for the change from using $g$ to $G$, cf.~\eqref{e.c42704} and~\eqref{e.c42710}.  
As we see below, the key observation in the proof of~\eqref{e.c42901}, and, as a consequence,
a key step in the proof of Theorem~\ref{t.rde-intro}.(i),
is that 
\be\label{jul1104}
G[\tilde u] - G[\vphi] \leq 0.
\ee
The analogous term had we not made the `swap' is $g[\tilde u] - g[\vphi]$; however, we are unable to prove this is nonpositive.  On the other hand, the cost of making the `swap' is the $- 2 \chi \alpha(\tilde u) s_x$ term in the left hand side of~\eqref{e.c42711}. 
As we have seen in the proof of \Cref{l.sfkpp_s},  this term, while it may not have a 
definite sign, did not pose an issue with constructing a supersolution.

The easiest term in (\ref{e.c42901}) is $\psi_x$.  
Here, we simply use its relationship to $\vphi$ and~\eqref{e.phase_plane} to find:
\be
	\psi_x
		= \frac{\psi}{\vphi}( \vphi_x + (1 - \sqrt \chi \alpha(\vphi)) \vphi)
		\leq \frac{\psi}{\vphi}( - \sqrt \chi \vphi(1 - \alpha(\vphi)) + (1 - \sqrt \chi \alpha(\vphi)) \vphi)
		= (1-\sqrt{\chi})\psi .
\ee
Now, we prove (\ref{jul1104}).
We estimate integral terms in the definition of $G$ using~\eqref{e.c42902}:
\be
	\begin{split}
	-\sqrt \chi\int_x^\infty &\left(\alpha''(\tilde u) \tilde u_y^2 
			- \alpha'(\tilde u) f(\tilde u)\right) \dy
		+ \sqrt \chi\int_x^\infty \left(\alpha''(\vphi) \partial_y \vphi^2 
			- \alpha'(\vphi) f(\vphi)\right) \dy
		\\&
		= -\sqrt \chi\int_0^{\tilde u} \Big(\alpha''(z) \eta(t,z)
			- \frac{\alpha'(z) f(z)}{\eta(t,z)}\Big) \dz
			+ \sqrt\chi\int_0^{\vphi} \Big( \alpha''(z) \eta_*(z) - \frac{\alpha'(z) f(z)}{\eta_*(z)} \Big) \dz
		\\&
		\leq -\sqrt\chi \int_0^{\tilde u} \Big(\alpha''(z) \eta_*(z)
			- \frac{\alpha'(z) f(z)}{\eta_*(z)}\Big) \dz
			+ \sqrt \chi \int_0^{\vphi} \Big( \alpha''(z) \eta_*(z) - \frac{\alpha'(z) f(z)}{\eta_*(z)} \Big) \dz
		\\&
		=  - \sqrt \chi\int_{\vphi}^{\tilde u} \Big( \alpha''(z) \eta_*(z) - \frac{\alpha'(z) f(z)}{\eta_*(z)} \Big) \dz.
	\end{split}
\ee
On the other hand, we can convert some of the non-integral terms in $G$ to integral terms:
\be
	\begin{split}
	\chi \tilde u &\alpha'(\tilde u) (1 - \alpha(\tilde u))
		- \chi \vphi \alpha'(\vphi) (1- \alpha(\vphi))
		\\&
		= \chi \int_{\vphi}^{\tilde u} (z\alpha'(z)(1- \alpha(z))' \dz
		= \chi \int_{\vphi}^{\tilde u} (\alpha'(z)(z - A(z))' \dz
		\\&
		= \chi \int_{\vphi}^{\tilde u} \alpha''(z)(z- A(z)) \dz
			+ \chi \int_{\vphi}^{\tilde u} \alpha'(z)(1- A'(z)) \dz
		\\&
		= \chi \int_{\vphi}^{\tilde u} \alpha''(z)(z- A(z)) \dz
			+ \chi \int_{\vphi}^{\tilde u} \alpha'(z)\Big(1- \frac{1}{\chi}\Big(\frac{f(z)}{z - A(z)}-1\Big)\Big) \dz.
	\end{split}
\ee
Using then~\eqref{e.phase_plane}, we find
\be
	\begin{split}
	\chi \tilde u &\alpha'(\tilde u) (1 - \alpha(\tilde u))
		- \chi \vphi \alpha'(\vphi) (1- \alpha(\vphi))
		\\&
		\leq \sqrt\chi \int_{\vphi}^{\tilde u} \Big(\alpha''(z) \eta_*(z)
			+ \chi \int_{\vphi}^{\tilde u} \alpha'(z)\Big(1- \frac{1}{\chi}\Big(\sqrt \chi\frac{f(z)}{\eta_*(z)}-1\Big)\Big) \dz.
	\end{split}
\ee
Hence, we have 
\be
	\begin{split}
	&-\sqrt \chi \int_x^\infty \left(\alpha''(\tilde u) \tilde u_y^2 
			- \alpha'(\tilde u) f(\tilde u)\right) \dy
		+ \sqrt \chi \int_x^\infty \left(\alpha''(\vphi) \partial_y \vphi^2 
			- \alpha'(\vphi) f(\vphi)\right) \dy
		\\&+ \chi \tilde u \alpha'(\tilde u) (1 - \alpha(\tilde u))
			- \chi \vphi \alpha'(\vphi) (1- \alpha(\vphi))
		\\&\qquad
		\leq \chi \int_{\vphi}^{\tilde u} \alpha'(z) \Big(1 + \frac{1}{\chi}\Big) \dz
		= (\chi+1) (\alpha(\tilde u) - \alpha(\vphi)).
	\end{split}
\ee
To summarize, we have so far arrived at:
\be
	\begin{split}
	G[\tilde u]	- G[\vphi]
		\leq &(\chi+1) (\alpha(\tilde u) - \alpha(\vphi))
			-(1+\chi) (\alpha(\tilde u) - \alpha(\vphi))
			\\&
			- 2 \sqrt \chi (\alpha(\tilde u) - \alpha(\vphi)) (1 - \sqrt \chi)
			= - 2 \sqrt \chi(1-\sqrt \chi)(\alpha(\tilde u) - \alpha(\vphi)).
	\end{split}
\ee
A similar  computation yields
\be
	\begin{split}
		G[\tilde u]
			\leq &- 2 \sqrt \chi(1-\sqrt \chi)\alpha(\tilde u) .
	\end{split}
\ee
Combining all of the above, we have
\be
	\begin{split}
		2\sqrt \chi &(\alpha(\tilde u) - \alpha(\vphi))\psi_x
			+ G[\tilde u] s + (G[\tilde u]	- G[\vphi]) \psi
		\\&
		\leq 
			2 \sqrt \chi (\alpha(\tilde u) - \alpha(\vphi))(\vphi_x + (1 - \sqrt \chi \alpha(\vphi)) \varphi) \frac{\psi}{\vphi}
			\\&\qquad
			- 2\sqrt \chi(1-\sqrt \chi)\alpha(\tilde u) s
			- 2 \sqrt \chi(1-\sqrt \chi) (\alpha(\tilde u) - \alpha(\vphi)) \psi.
	\end{split}
\ee
Using again~\eqref{e.phase_plane} to estimate $\vphi_x$, in addition to the fact that $\tilde u > \vphi$ (recall~\eqref{e.c42802}) and the nonegativity of $\psi$, we find
\be
	\begin{split}
	2 \sqrt \chi &(\alpha(\tilde u) - \alpha(\vphi))\psi_x
			+ G[\tilde u] s + (G[\tilde u]	- G[\vphi]) \psi
		\\&
		\leq 
			2 \sqrt \chi(\alpha(\tilde u) - \alpha(\vphi))(-\sqrt \chi \vphi(1-\alpha(\vphi)) + (1 - \sqrt \chi \alpha(\vphi)) \varphi) \frac{\psi}{\vphi}
			\\&\qquad
			- 2\sqrt \chi(1-\sqrt \chi)\alpha(\tilde u) s
			- 2 \sqrt \chi(1-\sqrt \chi) (\alpha(\tilde u) - \alpha(\vphi)) \psi
		\\&
		=	- 2\sqrt \chi(1-\sqrt \chi)\alpha(\tilde u) s \leq 0,			
	\end{split}
\ee
where the last inequality follows due to the assumed nonnegativity of $s$.

Thus,~\eqref{e.c42901} is established, completing the proof.~$\Box$

\appendix

\section{Auxiliary results on nonlinearities and traveling waves}\label{sec:appendix} 

In this appendix, we discuss some auxiliary elementary facts related to the
models introduced in this paper.

\subsection{Non-triviality of the semi-FKPP range}

We first show that the semi-FKPP range is nontrivial.  
\begin{lem}\label{l.semiFKPP_range}
Suppose that $A \in C^2([0,1])$ satisfies~\eqref{mar718} and $\chi_{FKPP}$ is defined by~\eqref{mar804}.  Then, we have $\chi_{FKPP} \leq 1/2$.
\end{lem}
We note that a sharper bound than $1/2$ could be obtained by taking into account the order of vanishing where $A$ ``lifts away'' from~$0$.  If the order is $n-1$, one obtains $1/n$.  We opt for a simpler proof since sharpness is not our goal here.

\smallskip
\noindent
{\bf Proof.}
If $A''(0)>0$, one obtains this from L'Hopital's rule:
\be
	\begin{split}
		\chi_{FKPP}
			&\leq \lim_{u \to 0} \frac{A(u)}{A'(u)(u-A(u))}
			= \lim_{u \to 0} \frac{A'(u)}{A'(u)(1-A'(u)) + A''(u)(u-A(u))}
			\\&
			= \lim_{u \to 0} \frac{A''(u)}{A''(u)(1-A'(u)) - A'(u) A''(u) + A'''(u)(u-A(u)) + A''(u) (1 - A'(u))}
			\\&
			= \frac{A''(0)}{ A''(0) + A''(0)}
			= \frac{1}{2}.
		\end{split}
\ee
This concludes the proof in this case.

If $A''(0) = 0$, 
let $u_\eps$ be the smallest $u>0$ so that $A''(u) = \eps$:
\be
	u_\eps = \inf\{u : A''(u) = \eps\}.
\ee
The existence of such a $u_\eps$ follows from~\eqref{mar718} for $\eps$ sufficiently small.  
Clearly, $u_\eps$ is increasing in $\eps$.  Hence, there is $u_0\in [0,1]$ such that $u_\eps \to u_0$ as $\eps \to 0$.  By construction, we have 
\be
A(u_0) = A'(u_0) = A''(u_0) = 0.
\ee
The definition of $u_\eps$ and (\ref{mar718}) imply that 
\be
A'(u_\eps)=\int_0^{u_\eps}A''(u)du \leq \eps u_\eps,
\ee
and, more generally, we have
\be
A'(u)=\int_0^{u}A''(v)dv\le \eps u,~~\hbox{ for all $0\le u\le u_\eps$.}
\ee
This gives
\be
A(u_\eps) =\int_0^{u_\eps}A'(u)du\leq\eps u_\eps^2/2.
\ee
Arguing again by L'Hopital's rule, we find
\be
	\begin{split}
		&\chi_{FKPP}
			\leq \lim_{u \to u_0} \frac{A(u)}{A'(u)(u-A(u))}
			= \lim_{\eps \to 0} \frac{A(u_\eps)}{A'(u_\eps)(u_\eps - A(u_\eps)}
			\\&~~~
			= \lim_{\eps \to 0} \frac{A'(u_\eps)}{A'(u_\eps)(1-A'(u_\eps)) + A''(u_\eps)(u_\eps-A(u_\eps))}
			= 	\lim_{\eps \searrow 0} \frac{1}{1-A'(u_\eps) + \frac{A''(u_\eps)}{A'(u_\eps)}(u_\eps-A(u_\eps))}
			\\&~~~
			\leq \lim_{\eps \to 0} \frac{1}{1 - \eps u_\eps + \frac{\eps}{\eps u_\eps}(u_\eps - \eps u_\eps^2/2)}
			= 	\lim_{\eps\to 0} \frac{1}{2 - 3\eps u_\eps/2}
			= \frac{1}{2}.
	\end{split}
\ee
Thus, the proof is complete.
~$\Box$

\subsection{The pushed case}

We have mentioned in the introduction  that waves are pushed when $\chi > 1$ in~\eqref{jul1314}.  Although the pushed regime is not 
a focus of this paper, we provide a short proof of this fact for the sake of completeness.
\begin{prop}\label{p.large_chi-explicit}
Suppose that $\chi > 1$, the function $A(u)$ satisfies~\eqref{jul1802}, and
\be\label{e.c72201}
f(u)=f_\chi(u):= f'(0) (u - A(u)) (1 + \chi A'(u)).
\ee
Then, the minimal speed $c_*$ is given by
\be\label{aug314}
c_*=\sqrt{f'(0)}\Big(\frac{1}{\sqrt \chi} + \sqrt \chi\Big) > 2 \sqrt{f'(0)},
\ee
and the minimal speed traveling wave profile function is
\be\label{jul2408}
\eta_*(u)=\sqrt{\chi}(u-A(u)).
\ee
As a consequence, traveling wave solutions of~\eqref{e.rde} are pushed and have a purely exponential decay, as in (\ref{jul1310}). 
\end{prop}
{\bf Proof.}
We may assume that $f'(0) = 1$ without loss of generality.
To obtain an upper bound on~$c_*$, first observe that
\be\label{jul2404}
	f(u) = p_f(u) (c_\chi - p_f'(u))
	\qquad\text{where } ~ 
	c_\chi = \sqrt \chi + \frac{1}{\sqrt \chi} ~
	\text{ and }~ p_f(u) = \sqrt \chi( u- A(u)).
\ee
Then applying~\eqref{dec2808bis}, we obtain the upper bound
\be
	c_*
		\leq \sup_{u\in[0,1]} \Big( p_f'(u) + \frac{f(u)}{p_f(u)}\Big)
		=  c_\chi.
\ee
 
To prove the lower bound 
\be\label{aug316}
c_*\ge c_\chi,
\ee
we argue by contradiction, assuming that~$c_* < c_\chi$.  Let $\eta_*$ be the traveling wave profile function.   It follows from~\eqref{dec2828},~\eqref{dec2820}, and the remark in between
that, since $c_*<c_\chi$, we have  
\be\label{e.c72302}
	\eta_*'(0)
		= \frac{c_* + \sqrt{c_*^2 - 4}}{2}
		< \frac{c_\chi + \sqrt{c_\chi^2 - 4}}{2}
		= \sqrt \chi
		= p_f'(0).
\ee
As we see below, the contradiction comes from the fact that (\ref{e.c72302}) says that the traveling wave moving with the speed $c_\chi>c_*$ decays
faster that the minimal speed wave $U_*$.

\Cref{prop-dec2802} implies that
\be\label{e.c72301}
	p_f(u) (c_\chi - p_f'(u)) = f(u) = \eta_*(u)(c_* - \eta_*'(u)).
\ee
Integrating the identity~\eqref{e.c72301} from $0$ to $u$, letting $P$ and $N$ be the anti-derivatives of $p_f$ and $\eta_*$, 
respectively, and rearranging terms yields
\be\label{e.c72202}
	c_\chi P(u) - c_* N(u) = \frac{p_f^2(u)}{2} - \frac{\eta_*^2(u)}{2}.
\ee
Observe that, due to~\eqref{e.c72302}, we know that $p_f(u) > \eta_*(u)$ 
for all $u$ that are positive and sufficiently small. Hence, we have~$c_\chi P(u) > c_*N(u)$ for small $u$ as well.  Let $u_0$ be the
smallest $u>0$ so that~$p_f(u_0) = \eta_*(u_0)$.  This must exist because $p_f(1) = \eta_*(1)$.

From the construction of $u_0$, we have that
\be
	p_f(u) > \eta_*(u)
		\qquad\text{ for all } u \in (0,u_0).
\ee
Recalling that $c_*< c_\chi$, by assumption, it follows that
\be\label{e.c72203}
	c_\chi P(u_0) > c_* P(u_0) = c_* \int_0^{u_0} p_f(u) du
		> c_*\int_0^{u_0} \eta_*(u) du
		= c_*N(u_0).
\ee
On the other hand, by~\eqref{e.c72202} and the fact that $p_f(u_0) = \eta_*(u_0)$, we find
\be
	c_\chi P(u_0) = c_* N(u_0) + \frac{p_f(u_0)^2}{2} - \frac{\eta_*(u_0)^2}{2}
		= c_*N(u_0).
\ee
This contradicts~\eqref{e.c72203}.  We conclude that $c_* \geq c_\chi$, finishing the proof.
~$\Box$

As a consequence of Proposition~\ref{p.large_chi-explicit} and the comparison principle, we have the following. 
\begin{cor}\label{p.large_chi}
Suppose that $\chi > 1$, the function $A(u)$ satisfies~\eqref{jul1802}, $f(u)$ satisfies~\eqref{dec2804bis}, and
	\be\label{e.c72201bis}
		f(u) \geq f'(0) (u - A(u)) (1 + \chi A'(u)).
	\ee
Then, the minimal speed $c_*$ obeys a lower bound
	\be
		c_* \geq \sqrt{f'(0)}\Big(\frac{1}{\sqrt \chi} + \sqrt \chi\Big) > 2 \sqrt{f'(0)}.
	\ee
	As a consequence, traveling wave solutions of~\eqref{e.rde} are pushed.
\end{cor}

%

\subsection{The generality of the class of nonlinearities}\label{s.A}

We discuss here how restrictive the assumptions~\eqref{jul1314}-\eqref{jul1802} on the nonlinearity $f(u)$are.  
The main point is that, while there is always at least one pair $(\chi,A)$ such that the form~\eqref{jul1314} holds for a given $f(u)$
satisfying assumptions (\ref{dec2804bis}), it is not clear when $A$ satisfies the conditions in~\eqref{jul1802}.  For simplicity, we assume that $f'(0) = 1$.
	
We first show the existence of a pair $(\chi,A)$.  Indeed, letting $\chi$ be the larger solution of
\be
c_*= \frac{1}{\sqrt \chi} + \sqrt \chi,
\ee
and using the traveling wave profile function $\eta_*$, we can define,  as in (\ref{jul2408}), 
\be\label{e.c72504}
	A(u) =u-\frac{1}{\sqrt{\chi}} \eta_*(u). 
\ee
Then, from \Cref{prop-dec2804}, we have:
\be
	\begin{split}
	f(u)
		&= \eta_*(u)( c_* - \eta_*'(u))
		= \sqrt{\chi } (u-A(u)) ( c_* - \sqrt{ \chi }(1 - A'(u)))
		\\&
		= \sqrt{\chi} (u-A(u)) \left(\left(c_* - \sqrt{\chi}\right) + \sqrt{ \chi } A'(u)\right)
		\\&
		= \sqrt{\chi } (u-A(u)) \bigg(\frac{1}{\sqrt\chi} + \sqrt{ \chi} A'(u)\bigg)
		= (u-A(u))(1 + \chi A'(u)).
	\end{split}
\ee
Hence $f(u)$ is in the form of~\eqref{jul1314}.  However, the function $A(u)$ defined in (\ref{e.c72504}) need not satisfy the conditions in~\eqref{jul1802}.
Actually, when the minimal speed traveling wave has an extra linear factor in the exponential asymptotics (\ref{jul1312}), 
we do know from Proposition~\ref{dec28-prop06}
that $\eta_*(u)$ is not $C^2[0,1]$. Thus, in that case this construction 
does not give a function $A(u)$ satisfying the regularity assumptions in~\eqref{jul1802}.  

Neither does $A(u)$ given by  (\ref{e.c72504}) have to be convex.
Indeed, take any non-convex, smooth, positive $\eta(u)$ satisfying $\eta'(0) = 1$ and $\eta(0) = 0 = \eta(1)$, for example, 
\be
	\eta(u) = \Big( \frac{1}{2} u + \frac{1}{2} \frac{\sin(10 u)}{10}\Big) (1-u).
\ee
Then, defining $A(u)$ by~\eqref{e.c72504} with $\chi = 1$, and setting
\be\label{e.c72502}
	f(u) =  (u - A(u))(1 + A'(u))
\ee
we obtain a nonlinearity of the form~\eqref{jul1314} and satisfying~\eqref{e.rde}-\eqref{dec2804bis}, but 
for which the choice of~$A(u)$ in (\ref{e.c72504})  is not convex.


Next we discuss the uniqueness of $\chi$ and $A$.   Due to \Cref{p.large_chi-explicit}, 
uniqueness of $\chi>1$ and~$A(u)$ satisfying \eqref{jul1802} holds for nonlinearities of the form (\ref{e.c72201}). 
 On the other hand, uniqueness is not necessarily true in the pulled case.  Indeed, consider the following simple example: for any $n \geq 3$, take the nonlinearity
\be
	f(u) = u - u^n.
\ee
Immediately, we see that, with $\chi = 0$ and $A(u) = u^n$, this has the form of~\eqref{jul1314}:
\be
	f(u) = (u - u^n)(1 + 0 \cdot (n u^{n-1}))
	 = (u- A(u)) (1 + 0\cdot A'(u)).
\ee
There is, however, another decomposition of the form~\eqref{jul1314}: choosing $m = (n+1)/2$, $\tilde \chi = 1/(m-1)$, and $\tilde A(u) = u^m$, 
we can write 
\be
	\begin{split}
	f(u)
		&= u - u^m + u^m - u^n
		= u - u^m + u^{m-1}(u - u^m)
		= (u - u^m)(1 + u^{m-1})
		\\&
		= (u - \tilde A(u))(1 + \tilde \chi \tilde A'(u)).
	\end{split}
\ee
Hence, the choice of $\chi$ and $A(u)$ is not necessarily unique in the pulled case.

Let us comment that we saw above an example of $\chi$ and $A(u)$ for which $(u)$ has the form~\eqref{jul1314} but for which $A(u)$ is not convex.  
However, this example was pulled, so  we do not necessarily have uniqueness of this decomposition.  
It is, thus, possible that $f(u)$ has another decomposition of the form~\eqref{jul1314} with a new $A(u)$ that is convex.
We leave it as an open question to determine conditions on $f(u)$ that guarantee the existence of a {\em suitable} $\chi$ and $A(u)$ such that~\eqref{jul1314}-\eqref{jul1802} hold.

\subsection{Bounds for the traveling wave profile}

\begin{lem}\label{lem-jun1402}
Let $f(u)$ be of the form (\ref{jul1314}) with $\chi \in [0,1]$ and $A(u)$  that satisfies~\eqref{jul1802}.
Then, the minimal speed traveling wave profile satisfies
\begin{equation}\label{jun1402}
\sqrt{\chi}(z-A(z))\leq \eta_*(z)\leq z-A(z).
\end{equation}
\end{lem}
{\bf Proof.} 
The minimal speed  traveling wave $U(x)=U_*(x)$ satisfies 
\be\label{jun1404}
-c_*U' = U'' +f'(0) (U-A(U))(1+\chi A'(U)).
\ee
Introducing~$V=-U'>0$,  we write (\ref{jun1404}) as 
\begin{equation}\label{sep0310}
\frac{dU}{dx} = -V,~~\frac{dV}{dx} =-c_*V+ f'(0)(U-A(U))(1+\chi A'(U)).
\end{equation}
This leads to 
\begin{align}\label{ratio}
    \frac{dV}{dU} = c_* -\frac{f'(0)(U-A(U))(1+\chi A'(U))}{V}.
\end{align}
We claim that when $0\le \chi < 1$, all trajectories are trapped in the region $D_1$ bounded by the curves 
\[
\ell_1 = \big\{V =  {\sqrt{\chi f'(0)}}(U-A(U))\big\},
\]
and $\ell_2 = \{V = \sqrt{f'(0)}(U-A(U))\}$. The reason is that along $\ell_1$, we have
\be
\frac{dV}{dU} = c_*-\frac{\sqrt{f'(0)}}{\sqrt{\chi}}(1+\chi A'(U)),
\ee
and the slope of the curve $\ell_1$ itself is $\sqrt{\chi f'(0)}(1-A'(U))$. Since $U$ is decreasing along the trajectory,
the trajectories point into the region~$D_1$ along $\ell_1$  if
\be
c_*-\frac{\sqrt{f'(0)}}{\sqrt{\chi}}(1+\chi A'(U))\le \sqrt{\chi f'(0)}(1-A'(U)),\quad \text{for all} ~U\in[0,1].
\ee
This condition is satisfied, as 
\[
c_*=2\sqrt{f'(0)}\le \frac{\sqrt{f'(0)}}{\sqrt{\chi}}+\sqrt{\chi f'(0)}.
\]
On the other hand, the trajectories point into the region $D_1$ along $\ell_2$  if
\be
c_*-\sqrt{f'(0)}(1+\chi A'(U))\ge \sqrt{f'(0)}(1-A'(U)),\quad \text{for all} ~U\in[0,1].
\ee
This condition is also satisfied, as 
\[
c_*=2\sqrt{f'(0)}\ge 2\sqrt{f'(0)}+\sqrt{f'(0)}(\chi-1)A'(U),
\]
for all $U\in[0,1]$ and $0\le\chi\le1$, because $A(U)$ is increasing. 
Therefore, the minimal speed traveling wave trajectory lies in the region $D_1$.
Thus,  for any $0<z<1$ such that $z=U_*(x)$, we have 
\begin{equation}
\sqrt{f'(0)} \sqrt{\chi}(z-A(z)) \leq \eta_*(z)\leq \sqrt{f'(0)}(z-A(z)),
\end{equation}
which is~(\ref{jun1402}).~$\Box$

\bibliographystyle{plain}

\end{document}